\documentclass[11pt]{article}

\usepackage{geometry}

\usepackage{etoolbox}
\usepackage[utf8]{inputenc}
\usepackage[T1]{fontenc}
\usepackage{textcomp}

\usepackage{amsmath}
\usepackage{amsthm} 
\usepackage{mathtools}

\usepackage{amssymb}

\usepackage{pgfplots}
\pgfplotsset{compat=1.14}

\usepackage{beton}

\usepackage{eucal}
\usepackage[euler-digits,euler-hat-accent]{eulervm}
\usepackage{mathdots}

\DeclareSymbolFont{eulerup}{U}{zeur}{m}{n}
\DeclareMathSymbol{\uppi}{\mathord}{eulerup}{"19}
\DeclareMathSymbol{\upi}{\mathord}{eulerup}{"69}

\usepackage[numbers,sort&compress,square]{natbib}%

\mathtoolsset{mathic}

\usepackage{graphics}
\usepackage{newunicodechar}

\usepackage{breqn}
\DeclareFlexSymbol{\Gamma}  {Var}{latin}{00}
\DeclareFlexSymbol{\Delta}  {Var}{latin}{01}
\DeclareFlexSymbol{\Theta}  {Var}{latin}{02}
\DeclareFlexSymbol{\Lambda} {Var}{latin}{03}
\DeclareFlexSymbol{\Xi}     {Var}{latin}{04}
\DeclareFlexSymbol{\Pi}     {Var}{latin}{05}
\DeclareFlexSymbol{\Sigma}  {Var}{latin}{06}
\DeclareFlexSymbol{\Upsilon}{Var}{latin}{07}
\DeclareFlexSymbol{\Phi}    {Var}{latin}{08}
\DeclareFlexSymbol{\Psi}    {Var}{latin}{09}
\DeclareFlexSymbol{\Omega}  {Var}{latin}{0A}
\DeclareFlexSymbol{0}{Var}{latin}{30}
\DeclareFlexSymbol{1}{Var}{latin}{31}
\DeclareFlexSymbol{2}{Var}{latin}{32}
\DeclareFlexSymbol{3}{Var}{latin}{33}
\DeclareFlexSymbol{4}{Var}{latin}{34}
\DeclareFlexSymbol{5}{Var}{latin}{35}
\DeclareFlexSymbol{6}{Var}{latin}{36}
\DeclareFlexSymbol{7}{Var}{latin}{37}
\DeclareFlexSymbol{8}{Var}{latin}{38}
\DeclareFlexSymbol{9}{Var}{latin}{39}

\usepackage{hyphenat}
\usepackage{tabu} 

\usepackage[dutch,french,german,italian,english]{babel}

\newcommand{\dutch}[1]{\foreignlanguage{dutch}{#1}}

\newcommand{\german}[1]{\foreignlanguage{german}{#1}}

\usepackage{xspace}
\newcommand*\elide{\textup{[\,\dots]}\xspace}

\DeclareUnicodeCharacter{03A3}{\ensuremath{\Sigma}}
\DeclareUnicodeCharacter{03A6}{\ensuremath{\Phi}}
\DeclareUnicodeCharacter{03A8}{\ensuremath{\Psi}}
\DeclareUnicodeCharacter{03A0}{\ensuremath{\Pi}}
\DeclareUnicodeCharacter{03B1}{\ensuremath{\alpha}}
\DeclareUnicodeCharacter{03B2}{\ensuremath{\beta}}
\DeclareUnicodeCharacter{03B3}{\ensuremath{\gamma}}
\DeclareUnicodeCharacter{03B4}{\ensuremath{\delta}}
\DeclareUnicodeCharacter{03B5}{\ensuremath{\epsilon}}
\DeclareUnicodeCharacter{03B6}{\ensuremath{\zeta}}
\DeclareUnicodeCharacter{03B7}{\ensuremath{\eta}}
\DeclareUnicodeCharacter{03BB}{\ensuremath{\lambda}}
\DeclareUnicodeCharacter{03BC}{\ensuremath{\mu}}
\DeclareUnicodeCharacter{03BD}{\ensuremath{\nu}}
\DeclareUnicodeCharacter{03BE}{\ensuremath{\xi}}
\DeclareUnicodeCharacter{03C1}{\ensuremath{\rho}}
\DeclareUnicodeCharacter{03C3}{\ensuremath{\sigma}}
\DeclareUnicodeCharacter{03C6}{\ensuremath{\phi}}
\DeclareUnicodeCharacter{03C7}{\ensuremath{\chi}}
\DeclareUnicodeCharacter{03C9}{\ensuremath{\omega}}

\newcommand{\apart}{\ensuremath\mathrel{\#}}
\newcommand{\absval}[1]{\ensuremath \lvert #1\rvert}

\newcommand{\numberset}[1]{\mathbb{#1}} 

\newcommand\numnat{\numberset{N}}
\newcommand\numrat{\numberset{Q}}
\newcommand\numreal{\numberset{R}}

\newcommand{\known}[1]{\ensuremath{\mathord{K}#1}}

\makeatletter
\newcommand*{\IfItTF}{%
  \ifx\f@shape\my@test@it
    \expandafter\@firstoftwo
  \else
    \expandafter\@secondoftwo
  \fi
}
\newcommand*{\my@test@it}{it}
\makeatother

\newcommand{\cs}[2]{%
\if\relax\detokenize{#1}\relax%
CS%
\else
	\if\relax\detokenize{#2}\relax%
	\mbox{CS\IfItTF{\/}{}\(^{#1}\)}%
	\else%
	\mbox{CS\IfItTF{\/}{}\(^{#1}\text{\IfItTF{\!}{}#2}\)}%
	\fi	
\fi
}

\newcommand{\bks}[2]{%
\if\relax\detokenize{#1}\relax%
BKS%
\else
	\if\relax\detokenize{#2}\relax%
	\mbox{BKS\IfItTF{\/}{}\(^{#1}\)}%
	\else%
	\mbox{BKS\IfItTF{\/}{}\(^{#1}\text{\IfItTF{\!}{}#2}\)}%
	\fi	
\fi
}

\newcommand{\csop}{\Box}
\newcommand{\gro}{\mathrel{\raisebox{0.75pt}{\scalebox{0.7}{\(\circ\)}}\mkern-9mu>}}
\newcommand{\klo}{\mathrel{<\mkern-9mu\raisebox{0.75pt}{\scalebox{0.7}{\(\circ\)}}}}
\newcommand{\leftabstord}{\mathrel{\sqsubset}}

\DeclareMathOperator{\length}{length}

\DeclarePairedDelimiter{\abs}{\lvert}{\rvert}

\newcommand{\inverse}[1]{#1^{-1}}

\let\oldexists\exists
\renewcommand*{\exists}{\mathop{}\!\oldexists}
\let\oldforall\forall
\renewcommand*{\forall}{\mathop{}\!\oldforall}

\newcommand{\breakingcomma}{%
  \begingroup\lccode `~=`,
 \lowercase{\endgroup\expandafter\def\expandafter~\expandafter{~\penalty0 }}}

\newcommand{\mathlist}[1] {\ensuremath{\breakingcomma #1}}
\newcommand{\seq}[1]{\ensuremath{\breakingcomma\langle\mathlist{#1}\rangle}}

\theoremstyle{definition}
\newtheorem{thm}{Theorem}
\newtheorem{dfn}[thm]{Definition}
\newtheorem{crl}[thm]{Corollary}
\newtheorem{lem}[thm]{Lemma}
\newtheorem{pargmt}[thm]{Plausibility argument}
\newtheorem{prf}[thm]{Proof}
\newtheorem{pri}[thm]{Principle}
\newtheorem{wce}[thm]{Weak counterexample}

\makeatletter
\gdef\emailauthor#1#2{\stepcounter{ead}%
     \g@addto@macro\@elseads{\raggedright%
      \let\corref\@gobble
      \eadsep\texttt{#1}\def\eadsep{\unskip,\space}}%
}
\def\urlauthor#1#2{\g@addto@macro\@elsuads{\let\corref\@gobble%
    \raggedright\eadsep\texttt{#1}%
    \def\eadsep{\unskip,\space}}%
}
\makeatother

\usepackage[pdfencoding=auto, psdextra]{hyperref}

\usepackage{refcount}

\newcounter{marginruletmp}
\newcommand\tikzmark[1]{\tikz[overlay,remember picture] \node (#1) {};}

\newlength{\marginlinexshift}
\AtBeginEnvironment{quote}{\setlength{\marginlinexshift}{-1.5em}} 

\newcommand\StartMarginRule{%
\stepcounter{marginruletmp}%
\tikzmark{a}\label{a\themarginruletmp}%
  \ifnum\getpagerefnumber{a\themarginruletmp}=\getpagerefnumber{b\themarginruletmp} 
\else
\begin{tikzpicture}[overlay, remember picture]
    \draw[ultra thick,black]
      ([xshift=\marginlinexshift,yshift=1.5ex]a-|current page text area.east) --  ([xshift=\marginlinexshift]current page text area.south east);  
  \end{tikzpicture}%
  \fi%
}

\newcommand\EndMarginRule{%
\tikzmark{b}\label{b\themarginruletmp}%
  \ifnum\getpagerefnumber{a\themarginruletmp}=\getpagerefnumber{b\themarginruletmp}
  \begin{tikzpicture}[overlay, remember picture]
    \draw [ultra thick,black]
      ([xshift=\marginlinexshift,yshift=1.5ex]a-|current page text area.east) -- ([xshift=\marginlinexshift]b-|current page text area.east);
  \end{tikzpicture}%
  \else
  \begin{tikzpicture}[overlay, remember picture]
    \draw [ultra thick,black]
      ([xshift=\marginlinexshift]current page text area.north east) -- ([xshift=\marginlinexshift]b-|current page text area.east);
  \end{tikzpicture}%
  \fi%
}

\begin{document}


\title{The Creating Subject, the Brouwer-Kripke Schema, and infinite proofs}
\author{Mark van Atten\thanks{Sciences, Normes, Décision (CNRS/Paris 4), 1 rue Victor Cousin, 75005 Paris, France. As of September 1, 2018: Archives Husserl (CNRS/ENS), 45 rue d'Ulm, 75005 Paris, France. \url{vanattenmark@gmail.com}}}

\maketitle

\begin{tikzpicture}[overlay, remember picture]
      \node[text=black,
      execute at begin node={\begin{varwidth}{\textheight}},
      execute at end node={\end{varwidth}},
      rotate=90]
      at ([xshift=20mm,yshift=0mm]current page.west)
      {Forthcoming in \textit{Indagationes Mathematicae}.};
    \end{tikzpicture}

\begin{abstract}
Kripke’s Schema
(better the Brouwer-Kripke Schema)
and
the Kreisel-Troelstra Theory of the Creating Subject
were introduced 
around the same time
for
the same purpose,
that of analysing Brouwer’s
‘Creating Subject arguments’;
other applications have been found since.
I first look in detail at
a representative choice of 
Brouwer’s arguments.
Then I discuss
the original use 
of the Schema and the Theory,
their justification
from a Brouwerian perspective,
and 
instances of the Schema
that can in fact be found in Brouwer's own writings.
Finally,
I defend the Schema and the Theory
against
a number of objections
that have been made.
\end{abstract}





\begin{flushright}
\textit{Brouwer’s views may be wrong or crazy
(e.g.~self-contradictory),\\
but one will never find out without looking
at their more dubious aspects.}\\
Kreisel~\citep{Kreisel1967b}\vspace{-\baselineskip}
\end{flushright}

\tableofcontents

\section{Notation}

Unless noticed otherwise,
the following notation is used:

\begin{center}
\begin{tabu} to \textwidth {lX<{\strut}}
\(\seq{a_n}_n\) & a sequence of elements each indexed by \(n\) \\
\(i, j, k, m, n, p, v, w\) & variables ranging over the natural numbers \\
\(f,g,h\) & variables ranging over functions \\
\(j(x,y)\) & a pairing function \(\numnat \times \numnat  \rightarrow \numnat\)\\
\(r,s,t\) & variables ranging over real numbers\\
\(x,y,z\) & variables whose range depends on the context\\
\(A, B, C\) & variables ranging over propositions\\
\(P, Q, R\) & variables ranging over predicates\\
\(X, Y, Z\) & variables ranging over species\\
\(α, β,γ,ξ\) & variables ranging over choice sequences\\
\(\bar{α}m\) & the initial segment of α \(\seq{α(1), α(2), \dots, α(m)}\)\\
\end{tabu}
\end{center}

In quotations,
notation has been left unchanged.

\section{Introduction}

Can mathematical arguments
depend 
not only on mathematical objects and their properties, 
but also
on the temporal order in which some ideal mathematician 
establishes theorems about them?
Brouwer affirmed this,
and it led him to devise a number of 
reasonings
now known as
‘Creating Subject arguments’.
The first known example
occurs in the 1927 Berlin lectures~\citep{Brouwer1992}.%
\footnote{Brouwer says in the opening of the 1948 paper 
‘Essentially negative properties’
that he had given an example of
such reasoning ‘now and then in courses and lectures since 1927’~[\citealp[p.963]{Brouwer1948A};
trl.~\citealp[p.478]{Brouwer1975}].
Known places are 
the second Vienna lecture (see further on in the text),
the 1933 Groningen lectures,
and the 1934 Geneva lectures.
The latter two have remained unpublished,
but 
Niekus~[\citealp[section 7]{Niekus2010};~\citealp[section 10]{Niekus2017}]
makes some comments on them.
There is no Creating Subject argument in the 1932 lecture
‘Will, knowledge, and speech’~\citep{Brouwer1933A2}.}
Let
\(\numreal\)
be understood
as the species of intuitionistic real numbers
(convergent choice sequences).
Brouwer considers the possibility of
an order relation on 
\(\numreal\)
based on 
‘the naive “before” and “after”\,’
according to which
\(r < s\)
if and only if 
on the intuitive continuum
(seen with the mind’s eye),
\(r\)
appears to the left
of
\(s\).%
\footnote{In the second Vienna lecture,
this is called the
‘natural order’~\citep[p.8]{Brouwer1930A}
and in the Cambridge lectures 
the ‘intuitive order’~\citep[p.43]{Brouwer1981A}.}

\begin{wce}[{{\citep[p.31–32]{Brouwer1992}}}]%
\label{L136}%
Let
\(<\)
be the naive order relation.
Then there is no hope of showing  that 
\begin{equation}
\forall x\in \numreal \forall y \in \numreal (x \neq y \rightarrow x < y \vee y < x)
\end{equation}
\end{wce}

\begin{pargmt}%
\label{L117}
Let
\(P\)
be a unary predicate
and
\(e\)
a mathematical object in its domain,
such that at present
neither
\(\neg P(e)\)
nor
\(\neg \neg P(e)\)
have become evident.

Define a real number
\(r\) 
as a convergent choice sequence
of rationals
\(r(n)\):
\begin{itemize}
\item As long as,
by the choice of 
\(r(n)\),
one
has obtained evidence neither of 
\(P(e)\)
nor of 
\(\neg P(e)\),
\(r(n)\) is chosen to be 0.

\item If between the choice of 
\(r(m-1)\) and
\(r(m)\),
one
has obtained evidence
of 
\(P(e)\),
\(r(n)\)
for all
\(n \geq m\) 
is
chosen to be
\(2^{-m}\).

\item If between the choice of 
\(r(m-1)\) and
\(r(m)\),
one
has obtained evidence
of 
\(\neg P(e)\),
\(r(n)\)
for all
\(n \geq m\) 
is
chosen to be
\(-2^{-m}\).
\end{itemize}
Then it is true that
\(r \neq 0\)
but,
as long as
neither
\(\neg P(e)\)
nor
\(\neg \neg P(e)\)
has become evident,
neither
\(r < 0\)
nor
\(r > 0\)
is true.
As one may always expect  to be able to find
such
\(P\)
and
\(e\),%
\footnote{See on this point also the remark on the relation between completely open problems and untested propositions on p.\pageref{L130} below.}
there is no hope of showing 
\(\forall x\in \numreal \forall y \in \numreal (x \neq y \rightarrow x < y \vee y < x)\).
\end{pargmt}

The justification is a plausibility argument
(for the conclusion that 
the given proposition
will never be demonstrated),
not a proof,
because it depends on an expectation,
albeit one that is highly likely to be fulfillable.
Thus,
the justification contains not only mathematical reasoning,
but also a value judgement.
This is the defining characteristic of weak counterexamples.%
\footnote{As the referee emphasised, 
Brouwer never calls his weak counterexamples 
‘theorems’.}

In the Berlin lectures,
this weak counterexample serves as a motivation 
for introducing
the notion of
\textit{virtual order},
and it is shown that that is the closest one can come to
an order on 
\(\numreal\);
we will come back to that in subsections~\ref{L114}
and~\ref{L018}.
For here,
the salient feature of the argument is that
the choices in the sequence
\(r\)
depend not only on properties
of 
\(P\)
and
\(e\),
but also on
\textit{the moment}
at which a certain proof about them becomes available
to the mathematician constructing the sequence.
It will be clear that
this sequence
\(r\)
is therefore not determined by a law,
to the extent that one accepts the idea
that the mathematician’s activity
depends on its free choices in directing its efforts.
The role of freedom in intuitionistic mathematics will be discussed further in
subsections~\ref{L034}
and~\ref{L077}.

‘The mathematician’
here is evidently an idealised one,
who in particular is always there to make 
the
\(n\)-th choice,
however large
\(n\)
may become.
Brouwer baptised this idealised mathematician
‘\dutch{het scheppende subject}’
in Dutch in 1948~\citep[p.963]{Brouwer1948A}
and
‘the Creating Subject’\label{L115}
in English
in 1949~\citep[p.1246]{Brouwer1949C}.%
\footnote{%
The basic meaning of
‘\dutch{scheppend(e)}’,
the participial adjective
of the verb
‘\dutch{scheppen}’,
here is
‘bringing something into existence’.
According to the historical dictionary
\textit{\dutch{Woordenboek der Nederlandsche Taal}},
the alternative
‘\dutch{creatief}’
had already been introduced in Dutch when Brouwer chose  
‘\dutch{scheppend}’,
but
often takes on a restricted sense expressing that this bringing about happens
in an imaginative, 
original,
or artistic
way.
The
\textit{Oxford English Dictionary}
shows the same relation between
‘creating’
and
‘creative’.
To reflect the Subject’s ontological responsibility in its full generality,
then,
Brouwer’s chosen terms in Dutch and English are more apt
than the often seen
‘\dutch{creatief subject}’
and
‘creative subject’.}
This idealisation will be further explained 
in subsection~\ref{L105},
and an objection to the claim that Brouwer is making it will be discussed in
subsection~\ref{L065}.

Brouwer had given a weak counterexample with a very similar conclusion
a few years before,
in 1924:
\begin{wce}[{{\citep[p.3]{Brouwer1924N}}}]%
\label{L137}
There is no hope of showing that 
\begin{equation}
\forall x\in \numreal \forall y \in \numreal (x = y \vee x < y \vee y < x)
\end{equation}
\end{wce}

But the way he had argued for it then had been very different.

\begin{pargmt}%
\label{L118}
\mbox{}
\begin{quote}
Let
\(d_v\)
be the
\(v\)-th digit after the decimal point
in the decimal expansion of
\(\pi\)
and
\(m=k_n\),
when
in the ongoing decimal expansion
of
\(\pi\)
at
\(d_m\)
it happens for the
\(n\)-th time
that the part
\({d_m}{d_{m+1}}\dots d_{m+9}\)
of this decimal expansion
forms a sequence
\(0123456789\).
Let furthermore
\(c_v={(-1/2)}^{k_1}\)
if 
\(v \geq k_1\),
otherwise
\(c_v={(-1/2)}^v\);
then the infinite sequence
\mathlist{c_1, c_2, c_3, \dots}
defines a real number
\(r\),
for which neither
\(r=0\),
nor
\(r>0\),
nor
\(r<0\)
holds.~\citep[p.3]{Brouwer1924N}
\end{quote}
\end{pargmt}
Note that
the choices in the sequence
\(c_v\)
are not defined in terms of 
the moment at which
the Creating Subject 
comes to know the truth of a certain proposition,
as they would in 
Plausibility argument~\ref{L117}. 
As the decimal expansion of
\(\pi\)
is lawlike,
and for every
\(v\)
it is decidable whether
\(v \geq k_1\),
\(r\)
is itself lawlike.

When Brouwer wrote this,
no value  
\(k_1\)
was known,
but in the meantime it has been found
that
\(k_1=17,387,594,880\)~\citep{Borwein1998}; 
so \(r>0\).
But as Brouwer remarks in a footnote~\citep[p.3n4]{Brouwer1924N},
\(r\)
can be defined using any other decidable 
property of natural numbers
of which one knows neither that there is an instance,
nor that there cannot be one.
He came to call such a property 
a ‘fleeing property’.	

\begin{dfn}[{[\citealp[p.161]{Brouwer1929A};~\citealp[p.6–7]{Brouwer1981A}]}]%
\label{L113}
A 
\textit{fleeing property} 
is defined by a decidable predicate 
\(P\)
on the natural numbers
such that
at present
there is evidence of neither
\(\exists n P(n)\) 
nor 
\(\forall n \neg P(n)\).
The 
\textit{critical number} 
\(k_P\)
of
\(P\)
is the as yet hypothetical smallest number
\(k\)
such that
\(P(k)\).
\end{dfn}

This notion is used in the second of the known Creating Subject arguments,
which is made 
in the second Vienna lecture from 1928.
It was
published in 1930
and,
as observed by Dirk van Dalen~\citep{Dalen1999a},
that makes it the first occurrence of a Creating Subject argument in print. 
Its conclusion is that of Weak counterexample~\ref{L136}:

\begin{pargmt}
\mbox{}
\begin{quote}
That the continuum is 
\emph{not ordered} 
by the sequence of its elements as derived from intuition
[\german{durch die der Anschauung entnommene Reihenfolge ihrer Elemente}]
is shown by an element
\(p\)
determined by the convergent sequence
\mathlist{c_1, c_2, \dots},
for which I%
\footnote{%
\label{L135}Brouwer’s use of
‘I’
here is conventional,
and does not mean that this is no Creating Subject argument.
See subsections~\ref{L105} 
and~\ref{L065},
and compare footnote~\ref{L134}.} 
choose
\(c_1\)
to be the zero point
and every
\(c_{ν+1} = c_{ν}\)
with one exception:
As soon as I find a critical number
\(λ_f\)
of a certain fleeing property
\(f\),
I choose the next
\(c_{ν}\)
to be equal to
\(-2^{-ν-1}\),
and as soon as I find a proof of the absurdity
of such a critical number,
I choose the next
\(c_{ν}\)
to be equal to
\(2^{-ν-1}\).
This element 
\(p\)
is different from zero,
and yet it is neither smaller nor greater than zero.~[\citealp[p.7–8]{Brouwer1929A};
trl.~\citep[p.59]{Mancosu1998}, modified]
\end{quote}
\end{pargmt}
Like the 1924 argument,
this one depends on a fleeing property,
but here,
if a critical number is found,
it influences the choices in the sequence
not through its value
but through the moment at which that value was found.
Thus,
\(p\)
is not a lawlike real number.

Neither in the Berlin nor in the Vienna lecture
Brouwer paused 
to isolate the notion of Creating Subject,
even though it is clearly operative in the arguments,
and
the audiences will have easily missed this novelty~[\citealp[p.308–310]{Dalen1999a}; \citealp[p.516–517]{Dalen2013}].
It is 
in  a series of publications beginning in 1948
that
Creating Subject arguments
appear in central position.

Wider attention to that series,
and indeed to Creating Subject arguments as such,
was first drawn
by
Heyting’s discussion in his book
\textit{Intuitionism}
from 1956~\citep{Heyting1956}.
In the early comments\label{L132} 
by 
Heyting himself,
Van Dantzig,
and Kleene
(who had spent January-June 1950 at the University of Amsterdam)
the Creating Subject arguments
were considered
controversial 
or not mathematical at all~[\citealp{Dantzig1949};~\citealp[title of chapter 8]{Heyting1956};~\citealp[p.175]{Kleene.Vesley1965}],
because 
these arguments 
were considered to let in an empirical or quasi-empirical element.
On the other hand,
Kripke and Kreisel
in the mid-1960s
encouraged further discussion of Brouwer’s arguments
by presenting explicitations
in terms of 
what have become known as Kripke’s Schema
and the Theory of the Creative Subject.
(For reasons given on p.\pageref{L139}
and
p.\pageref{L138},
it is more appropriate to speak of
‘the Brouwer-Kripke Schema’
and 
‘the Theory of the Creating Subject’.)
And whereas Brouwer had used Creating Subject arguments
only in weak and strong counterexamples to classical principles,
other uses of the idea have been found,
mostly using the Brouwer-Kripke Schema.

The purpose of this paper is
to present Brouwer’s
Creating Subject arguments
in their original context,
to give a critical survey of the debate of these arguments,
and to argue for their intuitionistic correctness
via a defense of the intuitionistic validity
of the Brouwer-Kripke Schema and the correctness of the Theory of the Creating Subject.
Given the historical and 
philosophical concern with what was,
or might reasonably be inferred to have been,
Brouwer’s thinking when introducing and employing the Creating Subject arguments,
more recent alternatives to his proofs and notions,
or intended outright replacements thereof,
are taken to be of additional,
but not of primary interest.

\section{Three Creating Subject arguments}\label{L112}

This section is an exposition of three of Brouwer’s post-war
Creating Subject arguments
that,
together,
bring out the mathematical
and historical aspects relevant to the later debate.
Two more of Brouwer’s Creating Subject arguments 
will be discussed in subsection~\ref{L095}.

\subsection{Ordering relations, testability, judgeability, decidability}\label{L114}

In the counterexamples that will be discussed in detail,
the following relations and notions are used.

\begin{dfn}[{[\citealp[p.8–9]{Brouwer1930A};~\citealp[p.1246n]{Brouwer1949C}]}]
Let
\(x\),
and
\(y\)
be two real numbers
given by convergent choice sequences.%
\footnote{I read Brouwer’s
characterisations of choice sequences~[\citealp[p.245n3]{Brouwer1925A};~\citealp[p.323]{Brouwer1942A}]
in such a way that
they include 
that of a lawlike sequence
as a limiting case.}
Define
‘%
\(x\) 
\textit{coincides} 
with 
\(y\)’
by
\begin{equation}
x = y
\equiv
\forall n\exists m\forall i(i \geq m \rightarrow 
\absval{x(i)-y(i)} < 2^{-n})
\end{equation}
Define
‘%
\(x\) is 
\textit{measurably smaller} 
than 
\(y\)’%
\footnote{The English terminology is Brouwer’s own~\citep[p.3]{Brouwer1951}. 
In Dutch, Brouwer used ‘\dutch{tastbaar kleiner}’~\citep[lecture 11]{Brouwer1933},
literally ‘tangibly smaller’.
(Vesley~\citep[p.143]{Kleene.Vesley1965}
uses 
‘%
\(\geq 2^{-n}\)’.)}
by
\begin{equation}
x\klo y
\equiv
\exists m \exists n\forall i(i \geq m \rightarrow y(i)-x(i) > 2^{-n})
\end{equation}
Correspondingly,
\(y \gro x\)
means that
\(y\) 
is 
‘%
\textit{measurably greater}’
than
\(x\).
Further,
\begin{align}
& x \leq y \quad \equiv  \quad\neg(x \gro γ)\\
& x < y \quad \equiv \quad \neg(x \gro y) \wedge \neg(x=y) \label{L116}
\end{align}
\end{dfn}

The relation
\(x < y\)
is negative%
\footnote{Brouwer~\citep[p.963]{Brouwer1948A} 
defines
‘a simply negative property’
as 
‘the absurdity of a constructive property’;
here we have a negative property that is a conjunction of two simple ones.}
and
\(x \klo y\)
positive 
(an existence statement).
Brouwer calls
\(\klo\)
the 
‘constructive’~\citep{Brouwer1949A}
or
‘natural measurable’
order~\citep{Brouwer1951}
on the continuum,
because
\(x \klo y\)
exactly if,
on the intuitive continuum
on which a scale has been placed such that
\(0\)
is to the left of
\(1\),
\(x\)
is to the left of
\(y\).
Thus,
where it applies,
this order mathematically captures
the ‘naive’ order
\(<\)
that we saw in the introduction.

Brouwer calls
\(<\)
as defined in~\eqref{L116}
the ‘virtual’~\citep[e.g.][]{Brouwer1951}
or ‘negative’
order~\citep{Brouwer1949A}
on the continuum.
It is different from
the order
\(<\)
in the 
Creating Subject argument
in 
the introduction.%
\footnote{There are more places where
the order
\(<\)
is not the virtual one~[\citealp[p.8]{Brouwer1930A};~\citealp[section 8]{Brouwer1933};~\citealp[lecture 2]{Brouwer1934};~\citealp[p.40-41]{Brouwer1981A}].}
Brouwer’s recasting of that argument based on
the new definition is discussed 
in subsection~\ref{L007}
below.

From the natural measurable order,
Brouwer also defined 
‘apartness’:
\begin{equation}
x  \apart y 
\equiv
x \klo γ \vee x \gro y
\end{equation}
where 
the octothorpe
is pronounced ‘is apart from’.%
\footnote{The term ‘apart’ (Dutch
‘\dutch{verwijderd}’,
‘\dutch{plaatselijk verschillend}’~\citep[section 2]{Brouwer1923C1};
German  
‘\german{entfernt}’~\citep[p.254]{Brouwer1925E},
‘\german{örtlich verschieden}’~\citep[p.3]{Brouwer1919A}) 
and the definition are Brouwer’s~\citep[p.1246]{Brouwer1949C},
the notation was introduced by  Heyting~\citep[p.20]{Heyting1934}, 
who earlier had used 
\(ω\)~\citep[section 6]{Heyting1925}.}
We have~\citep[p.254]{Brouwer1925E}
\begin{equation}\label{L040}
\neg(x \apart y)  \leftrightarrow x=y
\end{equation}
and thereby
\begin{equation}\label{L030}
x \neq y  \leftrightarrow \neg \neg x \apart y
\end{equation}

Properties of
\(<\)
and
\(\klo\)
that are immediate
from the definitions are
\begin{equation}
x \gro 0 \rightarrow x > 0
\end{equation}
and also~\citep[p.1245 and 1248]{Brouwer1949C}
\begin{align}
\neg\neg x > 0 & \leftrightarrow x > 0\\
\neg\neg x > 0 & \leftrightarrow \neg\neg x \gro 0
\end{align}
and hence
\begin{equation}
x > 0 \leftrightarrow \neg\neg x \gro 0
\end{equation}
But,
as shown by a Creating Subject argument that will be discussed in subsection~\ref{L007},
we do not have
\begin{equation}
x > 0 \rightarrow  x \gro 0 
\end{equation}
and therefore neither
\begin{equation}
\neg\neg x \gro 0 \rightarrow  x \gro 0 
\end{equation}
This means that 
\(<\)
is the weaker order relation.
As we will see,
Brouwer used that Creating Subject argument also to show that
\(<\)
is not a complete order.

\begin{dfn}
A proposition 
\(A\) 
is
\textit{decidable}
if a method is known 
to  prove 
\(A \vee \neg A\).
Instead of
‘decidable’,
Brouwer wrote
‘judgeable’~\citep[e.g.][p.114]{Brouwer1955}.
Testability of 
\(A\)
is 
decidability of
\(\neg A\):
a proposition 
\(A\) 
is
\textit{testable}
if a method is known 
to prove 
\(\neg A \vee \neg \neg A\).%
\footnote{As a schema,
\(\neg A \vee \neg\neg A\)
has been called
‘weak excluded middle’.
The intermediate logic obtained
by adding it to intuitionistic logic
was first
studied by Yankov~\citep{Yankov1968}.}
\end{dfn}

While decidability implies testability,
if only a proof of
\(\neg \neg A\)
is known
then
\(A\)
has been tested without having been decided.
No proposition can be untestable or undecidable in an absolute sense,
as both
\(\neg(\neg A \vee \neg \neg A)\)
and
\(\neg(A \vee \neg A)\)
are contradictory.

\subsection{An argument from 1948}\label{L007}
In
‘Essentially negative properties’,
Brouwer presented
\begin{wce}[{{\citep{Brouwer1948A}}}]%
\label{L013}
There is no hope of showing that
\begin{equation}\label{L001}
\forall x\in \numreal (x \neq 0 \rightarrow x < 0 \vee x > 0)
\end{equation}
let alone of
\begin{equation}\label{L002}
\forall x\in \numreal  (x \neq 0 \rightarrow x \klo 0 \vee x \gro 0)
\end{equation}
\end{wce}
with the common notational variant  
\begin{equation}
\forall x\in \numreal (x \neq 0 \rightarrow x \apart 0)
\end{equation}

\begin{pargmt}\label{L005}
Let 
\(A\) be a proposition 
that is at present not
testable.
The Creating Subject constructs a  
choice
sequence 
\(r\) of rational numbers 
\(r(n)\):
\begin{itemize}
\item As long as,
by the choice of 
\(r(n)\),
the Creating Subject
has obtained evidence neither of 
\(A\)
nor of 
\(\neg  A\),
\(r(n)\) is chosen to be 0.

\item If between the choice of 
\(r(m-1)\) and
\(r(m)\),
the Creating Subject
has obtained evidence
of 
\(A\),
\(r(n)\) 
for all 
\(n \geq m\) 
is chosen to be
\(2^{-m}\).

\item If between the choice of 
\(r(m-1)\) and
\(r(m)\),
the Creating Subject
has obtained evidence
of 
\(\neg  A\),
\(r(n)\) 
for all 
\(n \geq m\) 
is chosen to be
\(-2^{-m}\).
\end{itemize}
The choice sequence 
\(r\) 
converges,%
\footnote{Let 
\(ε\)
be given,
and determine an
\(n\)
such that
\(2^{-n} < ε\).
Construct the sequence
\(r\)
up to
\(r(n)\),
which can be done as each choice is decidable.
If
\(r(n)=0\),
all further choices will be in the interval
\([-2^{-(n+1)},2^{-(n+1)}]\)
and hence within
\(ε\)
from one another.
If
\(r(n) \neq 0\),
then the choices in
\(r\)
have already been fixed,
and hence within
\(ε\)
from one another.}
hence 
\(r\) 
is a real number.
By its definition,
\begin{equation}
r = 0 \leftrightarrow \neg A \wedge \neg\neg A
\end{equation}
and so the Creating Subject knows that 
\(r \neq 0\).
But it also reasons
\begin{equation}\label{L009}
\begin{aligned}
& r < 0 & \text{(assumption)}\\
& \neg(r > 0) & \text{(def.~\(<\))}\\
& \text{There never is evidence of 
\(A\)} & \text{(def.~\(r\))}\\
& \neg A & \text{(meaning of 
\(\neg\))}\\
& \text{\(A\) has been tested} & \\
\end{aligned}
\end{equation}
and,
symmetrically,
\begin{equation}\label{L010}
\begin{aligned}
& r > 0 & \text{(assumption)}\\
& \neg(r < 0) & \text{(def.~\(>\))}\\
& \text{There never is evidence of 
\(\neg A\)} & \text{(def.~\(r\))}\\
& \neg\neg A & \text{(meaning of 
\(\neg\))}\\
& \text{\(A\) has been tested} & \\
\end{aligned}
\end{equation}

\noindent
In both cases,
the Creating Subject arrives at a conclusion that contradicts
the hypothesis that
\(A\) cannot yet be tested.
Hence,
as long as 
\(A\) cannot be tested,
\(r < 0 \vee r > 0\) cannot be proved,
and,
by implication, 
neither
can
\(r \klo 0 \vee r \gro 0\).
\end{pargmt}

This reasoning\label{L051}
could not be reproduced  
starting from the hypothesis that
\(A\)
is at present
undecidable,
as that would allow for
\(A\)
having been tested
(the Creating Subject may have a proof of
\(\neg\neg A\)),
in which case no contradiction arises.
On the other hand,
it can be reproduced starting 
from the hypothesis that
\(A\)
is 
untestable
but replacing
\(>\)
by
\(\gro\);
this strengthens the assumptions towards a contradiction
in~\eqref{L009}
and~\eqref{L010},
but yields a correspondingly weaker result,
namely,
a weak counterexample to~\eqref{L002},
without implying one to~\eqref{L001}.

For a weak counterexample only to~\eqref{L002},
from either an untestable or an undecidable
proposition
\(A\),
there is an alternative argument
(not presented by Brouwer):

\begin{pargmt}\label{L021}
As above,
but instead of~\eqref{L009}
and~\eqref{L010},
the reasoning is
\begin{equation}\label{L011}
\begin{aligned}
& r \klo 0 & \text{(assumption)}\\
& \exists n(r(n) < 0) & \text{(def.~\(\klo\))}\\
& \text{\parbox[t]{19em}{By the choice of 
\(r(n)\), the Creating Subject\\ \mbox{}\hspace{1em} has obtained evidence of 
\(\neg A\)}} & \text{(def.~\(r\))}\\
& \text{\(A\) has been decided (and hence tested)} & \\
\end{aligned}
\end{equation}
and
\begin{equation}\label{L012}
\begin{aligned}
& r \gro 0 & \text{(assumption)}\\
& \exists n(r(n) > 0) & \text{(def.~\(\gro\))}\\
& \text{\parbox[t]{19em}{By the choice of 
\(r(n)\), the Creating Subject\\ \mbox{}\hspace{1em} has obtained evidence of 
\(A\)}} & \text{(def.~\(r\))}\\
& \text{\(A\) has been decided (and hence tested)} & \\
\end{aligned}
\end{equation}
\end{pargmt}

The reason why this shorter argument does not
also provide 
a weak counterexample to~\eqref{L001}
if
\(\klo\)
is replaced by
\(<\),
which would weaken the assumption
in~\eqref{L011}
and~\eqref{L012}
from which to work towards a contradiction,
is that
although
by definition of
\(<\) 
we have
\begin{equation}
r < 0   \rightarrow \neg\neg\exists n(r(n) < 0)
\end{equation}
we do not have
\begin{equation}
r < 0   \rightarrow \exists n(r(n) < 0)
\end{equation}
To get from the former to the latter,
Markov’s Principle\label{L006}
(MP)%
\footnote{%
\label{L124}%
Markov 
introduced this principle in lectures 1952–1953~\citep[p.44]{Dragalin1988}
and called it
‘the Leningrad Principle’~\citep{Markov1956},
after the place where he worked before moving to Moscow in 1955,
and continued to use it much later~\citep[p.284n20]{Margenstern1995}. 
Another name he used was
‘the principle of constructive glean’~\citep[p.137]{UspenskySemenov1993}.
Finally,
it is also known as 
‘the Principle of Constructive Choice’~\citep[p.44]{Dragalin1988}.
Philosophical and historical aspects of
Markov’s engagement with the work of Brouwer and Heyting
have been analysed by Vandoulakis~\citep{Vandoulakis2015}.} 
should hold for non-recursive sequences.
MP comes in several forms~[\citealp[section 4]{Smorynski1973};~\citealp{Troelstra1975};~\citealp[section 4.5]{Troelstra.Dalen1988}];
the one relevant here is
\begin{equation}
\neg\neg\exists n (α(n) = 1) \rightarrow  \exists n (α(n) = 1)
\end{equation}
and was originally formulated and defended for 
recursive sequences.
I will not here go into the question whether for such
sequences it is valid.%
\footnote{It seems very likely to me that it is;
I give a plausibility argument
elsewhere~\citep[section 2]{Atten-forthcomingB}.}
But it would certainly not be correct
to apply MP to a sequence such as
\(r\).
Briefly,
the Creating Subject is free to go about its activity,
and hence its knowing,
from the hypothesis 
\(r < 0\), 
that it is not the case that there is no stage
at which it establishes
\(\neg A\),
does not suffice for knowing 
a stage at which it does.
For more on the non-recursive nature
of the Subject’s activity,
see subsection~\ref{L034}.

Brouwer 
(in effect) 
once used a form of MP himself,
in a paper from 1918;
this is discussed below in Appendix~\ref{L131}.

In the same paper from 1948 as that in which
Weak counterexample~\ref{L013} appears,
Brouwer also observes that 
a slight modification to the definition
of
\(r\)
yields the following  weak counterexample:
\begin{wce}[{{\citep{Brouwer1948A}}}]%
\label{L028}
There is no hope of showing that
\begin{equation}
\forall x\in \numreal (x > 0 \rightarrow x \gro 0)
\end{equation}
\end{wce}

\begin{pargmt}\label{L029}
As Plausibility argument~\ref{L005},
modifying the clause in the definition of
\(r\) 
for
the case that evidence has been obtained of 
\(\neg  A\)
by changing 
the number chosen from
\(-2^{-n}\)
to
\(2^{-n}\).
By the argument above,
\(r \neq 0\)
and now also
\(r > 0\),
but not,
as long as the proposition
\(A\)
has not been tested,
\(r \gro 0\).
\end{pargmt}

This also presents a weak counterexample to MP 
((mis)applied to non-recursive sequences)
that
predates the introduction of the latter.%
\footnote{See footnote~\ref{L124}.}

As Brouwer mentions at the beginning
of 
‘Essentially negative properties’~\citep{Brouwer1948A},
he was moved to publish these weak counterexamples
by Griss’ and Van Dantzig’s contention
that negation is not acceptable
in intuitionistic mathematics~[\citealp{Griss1944,Griss1946};~\citealp{Dantzig1947};~\citealp[section 6.2.2]{Atten2017}].
The philosophical conclusion that Brouwer draws from
his counterexamples is that
\(\neq\)
and
\(>\)
are 
\textit{essentially} 
negative
in the sense that there is no hope
of ever proving that they
are equivalent
to positive relations,
for 
he believes that there will always be untestable propositions.
Krivtsov~\citep[p.167]{Krivtsov2000}
and also Veldman (in conversation) have pointed out that
that conclusion is of a greater generality than
is supported by
Brouwer’s actual arguments.
While Brouwer shows that
\(\neq\)
does not coincide with
\(\apart\),
nor 
\(>\)
with
\(\gro\),
he does not address the question whether
there may not be other possibilities.
That question,
while crucial to Brouwer’s discussion
with Griss and Van Dantzig,
goes beyond that of the acceptability of
Brouwer’s counterexamples themselves,
and will be of no further concern here.

\subsection{A refinement by Heyting}\label{L018}

Because the Creating Subject does not know
in advance by which choice 
\(r(n)\)
it will have obtained evidence of
\(A\)
or,
as the case may be,
of
\(\neg A\),
it can not know in advance any properties of 
\(n\)
to know which would require knowing its exact value,.
In other words,
if the Subject does know such a property of 
\(n\),
it must have
decided
\(A\)
already.
The definition of a sequence
\(r\)
may be made to depend on such properties.

This idea was used by Heyting in 1956
to argue
that
we cannot expect to prove that
Brouwer’s order on the continuum
\(<\),
which is a virtual order,
is also a pseudo-order.
The notion of a virtual order was introduced 
by Brouwer in print in 1926~\citep[p.453]{Brouwer1926A},
that of a pseudo-order by Heyting 
in his dissertation~\citep[p.8]{Heyting1925}.
These order concepts have in common the axioms for partial order:
\begin{gather}
x \leftabstord y \rightarrow \neg(y \leftabstord x) \wedge \neg(x=y)\label{L023}\\
x=y \wedge y \leftabstord z \rightarrow x \leftabstord  z\\
x \leftabstord y \wedge y=z \rightarrow x \leftabstord z\\
x \leftabstord y \wedge y \leftabstord z \rightarrow x \leftabstord z\label{L024}
\end{gather}
A pseudo-order satisfies,
in addition, 
\begin{gather}
\neg (x \leftabstord y) \wedge  \neg(y \leftabstord x) \rightarrow x=y\label{L025}\\
x \leftabstord y \rightarrow \forall z(x \leftabstord z \vee z \leftabstord y)\label{L027}
\end{gather}
whereas a virtual order satisfies,
in addition to~\eqref{L023}–\eqref{L024} and~\eqref{L025},
\begin{gather}
\neg x \leftabstord y \wedge  \neg(x=y) \rightarrow y \leftabstord x \label{L026}
\end{gather}
Both notions
are weaker than that of a complete order,
which  satisfies
\begin{gather}
x = y \vee x \leftabstord y \vee y \leftabstord x
\end{gather}

As we saw in
subsection~\ref{L007},
Brouwer had shown 
that his order on the continuum 
\(<\),
which is a virtual order,
is weaker than
his order on the continuum
\(\klo\),
which is an example of Heyting’s pseudo-orders.
Heyting’s argument then makes it plausible that 
no other choice of pseudo-order
would have yielded an equivalence:
\begin{wce}[{{\citep[p.117]{Heyting1956}}}]\label{L046}
There is no hope of showing that Brouwer’s virtual order on the continuum is a pseudo-order.%
\footnote{In fact,
the Continuity Principle
(Principle~\ref{L076} below)
can be used to obtain the stronger result that
the proposition stating the equivalence of these
two orders is absurd.
(I thank the referee for pointing this out.)}
\end{wce}

\begin{pargmt}
Let 
\(A\) be a proposition 
that has not been tested. 
The Creating Subject constructs two 
real numbers,
\(r\) and 
\(s\):
\begin{itemize}
\item As long as,
by the choice of 
\(r(n)\),
the Creating Subject
has not tested 
\(A\),
\(r(n)\) is chosen to be 
\(2^{-n}\).

\item If between the choice of 
\(r(m-1)\) and
\(r(m)\),
the Creating Subject
has tested 
\(A\),
\(r(n)\)
for all
\(n \geq m\)
is 
chosen to be
\(2^{-m}\).
\end{itemize}
The sequence
\(s\)
depends on 
\(r\):
\begin{itemize}
\item As long as,
by the choice of 
\(s(n)\),
the Creating Subject
has not tested 
\(A\),
\(s(n)\) is chosen to be 
\(2^{-n}\).

\item If between the choice of 
\(r(m-1)\) and
\(r(m)\),
the Creating Subject
has tested 
\(A\)
and
\(m\)
is odd,
\(s(n)\)
for all
\(n \geq m\)
is 
chosen to be
\(2^{-n}\).

\item If between the choice of 
\(r(m-1)\) and
\(r(m)\),
the Creating Subject
has tested 
\(A\)
and
\(m\)
is even,
\(s(n)\) 
for all
\(n \geq m\)
is chosen to be
\(2^{-m}\).
\end{itemize}
Then
\(0 < r\).
By definition,
\(0 < s   \rightarrow s \neq 0\);
and
if
\(s \apart 0\),
the third clause in the definition of
\(s\)
must have applied,
so 
\(s \apart 0 \rightarrow s=r\),
and therefore,
by
\(s < r   \rightarrow s \neq r\)
and~\eqref{L040},
\(s < r   \rightarrow s = 0\).
Combining the two gives
\(0 < s \vee s < r    \rightarrow s \neq 0 \vee s=0\).
Truth of
\(s \neq 0\)
means that 
it is not possible that
\(A\)
will not be tested between the choices
\(r(m-1)\)
and
\(r(m)\)
for some  even
\(m\),
and truth of
\(s=0\)
means that 
\(A\)
will not be tested 
between the choices
\(r(m-1)\)
and
\(r(m)\)
for some even
\(m\).
Either alternative cannot be known before
\(A\)
has,
in fact,
been tested.
On the hypothesis that
\(A\)
has not yet been tested,
the Creating Subject knows that
\(0 < r\),
but cannot yet know that
\(0 < s \vee s < r \).
Therefore,
axiom~\eqref{L027}
does not hold for
\(<\).
\end{pargmt}

Brouwer proved 
that any virtual order 
is an inextensible order
(and vice versa),
in the sense that any consistent addition to it is
already contained in it~[\citealp[ch.3]{Brouwer1992};~\citealp{Brouwer1950B}].%
\footnote{Martino~\citep{Martino1988}
analyses Brouwer’s
(changing)
argument.}
Van Dalen~\citep[p.12]{Brouwer1992}
observes that for Brouwer,
what with his drive for generality,\label{L052}
this seems to have resulted in a fondness for 
\(<\);
that would explain why Brouwer,
unlike Heyting,
never isolated the notion of pseudo-order.%
\footnote{In the 1934 Geneva lectures,
Brouwer precedes the definition of
\(\klo\)
with the remark
‘This will be just an auxiliary relation for the virtual order’~\citep[lecture 2, trl.~MvA]{Brouwer1934}.}
But Heyting
and later
constructivists
have
valued the pseudo-order
\(\klo\)
because its positivity makes it more practical,
and for 
(notational)
simplicity write it as
\(<\)~[\citealp[p.107]{Heyting1956};~\citealp[p.29n1]{Brouwer1992}].

\subsection{Drifts and checking numbers}\label{L071}

To be able to discuss weak counterexamples
in general terms,
Brouwer introduces the following definitions~[\citealp{Brouwer1948B};~\citealp{Brouwer1949A};~\citealp[p.1246]{Brouwer1949C}].

\begin{dfn}\label{L070} 
A 
\textit{drift}
\(γ\)
is the  union%
\footnote{The term is of course to be taken in Brouwer’s
species-theoretical sense,
not the 
set-theoretical sense:
the union of two species 
\(M\)
and
\(N\)
is the species of objects that have either been proved to be an element of the species
\(M\),
or been proved to be an element of the species
\(N\).~\citep[p.4]{Brouwer1918B}.
Thinking set-theoretically,
one would construe a drift as an ordered pair.} 
of
the real numbers 
in
a converging sequence
\(\seq{c_n(γ)}_n\)
and
the real number
\(c(γ)\)
it
converges to,
where
the 
\(c_n(γ)\)
are all apart from one another and
from 
\(c(γ)\).%
\footnote{The metaphor Brouwer has in mind is that of 
a sand drift;
‘drift’
translates
‘\dutch{aanstuiving}’~[\citealp[p.1239]{Brouwer1948B};~\citealp[p.122]{Brouwer1949A}]
which originally refers 
to the wind’s blowing sand to a given place,
thereby forming a dune.
The related term
‘checking number’ 
in Definition~\ref{L093}
is his translation of 
‘\dutch{dempingsgetal}’,
derived from the verb ‘\dutch{dempen}’,
one of whose meanings is 
‘to check’,
as in
‘to check a sand drift’.
(The historical dictionary
\textit{\dutch{Woordenboek der Nederlandsche Taal}}
gives the example phrase
‘\dutch{Het dempen der zandverstuivingen onder Eerbeek wordt op kleine schaal voortgezet}’.)
In the inland area where Brouwer lived,
sand drifts and the resulting dunes were,
and are,
well known.} 
The 
\(c_n(γ)\)
are 
the 
\textit{counting numbers} 
of the drift
and
\(c(γ)\)
its 
\textit{kernel}.
A counting number
\(c_n(γ) \klo c(γ)\)
is a
\textit{left counting number},
and a 
counting number
\(c_n(γ) \gro c(γ)\)
a
\textit{right counting number}.
If all counting numbers of a drift
are left counting numbers,
the drift is 
\textit{left-winged};
if
all are right counting numbers,
\textit{right-winged}.
If the sequence
\(\seq{c_n(γ)}_n\)
is the union of 
an infinite sequence
\(\seq{l_n(γ)}_n\)
of 
left counting numbers
and
an infinite sequence
\(\seq{d_n(γ)}_n\)
of 
right counting numbers,
the drift is 
\textit{two-winged}.
\end{dfn}

\begin{dfn}\label{L093}
Let
\(A\)
be a proposition
and 
\(γ\)
a drift.
The
\textit{direct checking number
of
\(γ\)
through
\(A\)}
is the real number 
\(D(γ,A)\),
constructed as follows.

The Creating Subject
constructs
a choice sequence
\(R(γ,A)\)
of real numbers
\mathlist{c_1(γ,A),c_2(γ,A),\dots}:
\begin{itemize}
\item As long as,
by the choice of 
\(c_n(γ,A)\),
the Creating Subject
has obtained evidence neither of 
\(A\)
nor of 
\(\neg  A\),
\(c_n(γ,A)\) is chosen to be 
\(c(γ)\).

\item 
\begin{tabular}[t]{@{}l@{ }p{25em}}
If 
\(n=1\) & and  
the Creating Subject
has obtained evidence
either of 
\(A\) or of 
\(\neg  A\)
before making the choice of 
\(c_{1}(γ,A)\),\\
or 
\(n>1\) &
and
the Creating Subject
has obtained evidence
either of 
\(A\) or of 
\(\neg  A\)
between the choice of 
\(c_{m-1}(γ,A)\) and
that of 
\(c_{m}(γ,A)\)
for some
\(m \leq n\),\\
then &
\(c_{n}(γ,A)\) 
is chosen to be
\(c_m(γ)\).
\end{tabular}
\end{itemize}
\(D(γ,A)\)
is the 
real number
to which
the sequence
\(R(γ,A)\)
converges.%
\footnote{Let 
\(ε\)
be given.
First consider the sequence
\(c_n(γ)\).
By hypothesis,
it converges to
\(c(γ)\),
so 
a number
\(n_0\)
can be constructed 
such that all values from place
\(n_0\)
onward
are within
\(ε/2\)
from one another
and from
\(c(γ)\).
With this
\(n_0\)
in hand,
turn to the sequence
\(R(γ,A)\)
and
assume that two values
\(c_k(γ,A)\)
and
\(c_m(γ,A)\)
have been constructed,
\(n_0 \leq k < m\).
Then it is decidable whether, 
by the time of constructing
\(c_k(γ,A)\),
\(A\) has been proved.
If it has,
then 
\(R(γ,A)\) has become constant,
and 
\(c_k(γ,A)\) and 
\(c_m(γ,A)\)
certainly lie within 
\(ε\) from one another.
If 
\(A\) has not been proved 
by the time of constructing 
\(c_k(γ,A)\),
but it has been proved by the time of constructing
\(c_m(γ,A)\),
then the values
at places 
\(k\) to 
\(m\) in 
\(R(γ,A)\)
are all within
\(ε/2\)
from
\(c(γ)\),
and 
hence within
\(ε\)
from one another.
Finally,
if 
\(A\) has not yet been proved by the time of
constructing
\(c_m(γ,A)\),
then the values
at places 
\(k\) to 
\(m\) in 
\(R(γ,A)\)
are just 
\(c(γ)\),
and hence
within 
\(ε\) from one another.
So in the sequence
\(R(γ,A)\),
all values from place
\(n_0\)
onward
lie within
\(ε\)
from one another.}

From this definition,
it follows that
\(D(γ,A) \neq c(γ)\),
that
\(\neg\neg\exists n(D(γ,A)=c_n(γ))\),
and that
\(\exists n(D(γ,A)=c_n(γ))\)
cannot be proved until
\(A\)
has been decided.
\end{dfn}

For the discussion
of BKS further on,
it should be noted that in Definition~\ref{L093}
(and Definitions~\ref{L110} and~\ref{L111})
Brouwer poses no requirement
that 
\(A\)
be untestable or undecidable.
Indeed,
the possibility that 
\(A\)
has been decided before the construction
of the direct checking number begins
is explicitly taken into account,
and in a footnote Brouwer mentions that
untested 
\(A\)
form 
(implicitly, not the general but)
a special case.

\begin{dfn}\label{L110}
Let
\(A\)
be a proposition
and 
\(γ\)
a drift.
The
\textit{conditional checking number
of
\(γ\)
through
\(A\)}
is the real number
\(C(γ,A)\)
constructed as follows.

The Creating Subject
constructs
a choice sequence
\(Q(γ,p)\)
of real numbers
\mathlist{c_1(γ,A),c_2(γ,A),\dots}:
\begin{itemize}
\item As long as,
by the choice of 
\(c_n(γ,A)\),
the Creating Subject
has not obtained evidence of 
\(A\),
\(c_n(γ,A)\) is chosen to be 
\(c(γ)\).

\item 
\begin{tabular}[t]{@{}l@{ }p{25em}}
If 
\(n=1\) & and  
the Creating Subject
has obtained evidence of 
\(A\) 
before making the choice of 
\(c_{1}(γ,A)\),\\
or 
\(n>1\) &
and
the Creating Subject
has obtained evidence
of 
\(A\) 
between the choice of 
\(c_{m-1}(γ,A)\) and
that of 
\(c_{m}(γ,A)\)
for some
\(m \leq n\),\\
then &
\(c_{n}(γ,A)\)
is
chosen to be
\(c_m(γ)\).
\end{tabular}
\end{itemize}
\(C(γ,A)\)
is the real number to which
the sequence
\(Q(γ,A)\)
converges.
\end{dfn}

\begin{dfn}\label{L111}
Let
\(A\)
be a proposition
and 
\(γ\)
a two-winged drift
with
\(\seq{c_n(γ)}_n\)
the union of
\(\seq{l_n(γ)}_n\)
and
\(\seq{d_n(γ)}_n\).
The
\textit{two-sided checking 
number
of
\(γ\)
through
\(A\)}
is the real number
\(E(γ,A)\)
constructed as follows.

The Creating Subject
constructs
a choice sequence
\(S(γ,A)\)
of real numbers
\mathlist{c_1(γ,A),c_2(γ,A),\dots}:
\begin{itemize}
\item As long as,
by the choice of 
\(c_n(γ,A)\),
the Creating Subject
has obtained evidence neither of 
\(A\)
nor of 
\(\neg  A\),
\(c_n(γ,A)\) is chosen to be 
\(c(γ)\).

\item 
\begin{tabular}[t]{@{}l@{ }p{25em}}
If 
\(n=1\) & and  
the Creating Subject
has obtained evidence
of 
\(A\)
before making the choice of 
\(c_{1}(γ,A)\),\\
or 
\(n>1\) &
and
the Creating Subject
has obtained evidence of
\(A\)
between the choice of 
\(c_{m-1}(γ,A)\) and
that of 
\(c_{m}(γ,A)\)
for some
\(m \leq n\),\\
then &
\(c_{n}(γ,A)\)
is chosen to be
\(d_m(γ)\).
\end{tabular}

\item 
\begin{tabular}[t]{@{}l@{ }p{25em}}
If 
\(n=1\) & and  
the Creating Subject
has obtained evidence
of 
\(\neg  A\)
before making the choice of 
\(c_{1}(γ,A)\),\\
or 
\(n>1\) &
and
the Creating Subject
has obtained evidence
of 
\(\neg  A\)
between the choice of 
\(c_{m-1}(γ,A)\) and
that of 
\(c_{m}(γ,A)\)
for some
\(m \leq n\),\\
then &
\(c_{n}(γ,A)\)
is
chosen to be
\(l_m(γ)\).
\end{tabular}
\end{itemize}
\(E(γ,A)\)
is the real number
to which the sequence
\(S(γ,A)\)
converges.
\end{dfn}

Plausibility argument~\ref{L005} in subsection~\ref{L007}
can be rephrased in terms of a 
two-sided checking number,
based on a
two-winged drift
\(γ\)
with 
kernel 
\(0\)
and
counting numbers
\(\seq{l_n(γ)}_n=\seq{-2^{-1}, -2^{-2},\dots}\)
and
\(\seq{d_n(γ)}_n=\seq{2^{-1}, 2^{-2},\dots}\)
An argument involving a
conditional checking number
is discussed in the next subsection.

\subsection{An argument from 1949}\label{L066}
In 1949 
Brouwer 
improved on the 1948 result
and devised
a strong counterexample:
\begin{thm}[{{\citep{Brouwer1949A}}}]%
\label{L004}
\begin{equation}
\neg\forall  x\in \numreal(x > 0 \rightarrow x \gro 0) 
\end{equation}
\end{thm}

\begin{crl}
Not stated by Brouwer,
but immediate from his definitions:
\begin{equation}
\neg\forall  x\in \numreal(x \neq 0 \rightarrow x \apart 0) 
\end{equation}
\end{crl}

\begin{crl}[{{\citep[p.205–206]{Troelstra.Dalen1988}}}]\label{L119}
Not stated by Brouwer, who never isolated MP:
\begin{equation}
\neg\text{MP}
\end{equation}
\end{crl}

\begin{prf}[of Corollary~\ref{L119}]
Assume MP\@.
Let
\(r \in \numreal\)
be an arbitrary real number,
given as a choice sequence of rationals, 
for which
\(r > 0\).
Without loss of generality,
we may assume that
\(\forall n(\absval{r(n)-r} < 2^{-n})\).
The assumption 
\(r > 0\)
implies
\(\neg\forall n(\absval{r(n)} < 2^{-n})\)
and hence
\(\neg\neg \exists n (\absval{r(n)}) \geq 2^{-n})\).
Define 
\(α\)
by
\(α(n) = 0 \text{ if } \absval{r(n)} < 2^{-n}\)
and
\(α(n) = 1 \text{ if } \absval{r(n)}  \geq 2^{-n}\).
Then
\(\neg\neg \exists n (α(n)=1)\)
and,
by MP, 
\(\exists n (α(n)=1)\).
Hence
for some
\(n\),
\(\absval{r(n)} \geq 2^{-n}\).
As, 
by assumption,
\(\absval{r(n)-r}  < 2^{-n}\),
also,
for some
\(m\),
\(\absval{r(n)-r} < 2^{-n} - 2^{-m}\),
and
therefore
\(r > 2^{-m}\).
By arbitrariness of
\(r\),
we conclude
\(\forall  x\in \numreal(x > 0 \rightarrow x \gro 0)\),
which contradicts
Theorem~\ref{L004}.
\end{prf}

Brouwer’s proof of Theorem~\ref{L004}
contains a mistake,
as
was pointed out by Myhill,
who also showed how to repair it.
We will look at Brouwer’s argument in this subsection,
and at Myhill’s reaction in the next.

It is not possible
to establish Theorem~\ref{L004}
by first proving 
\begin{equation}\label{L128}
\exists x\in \numreal\neg(x > 0 \rightarrow x \gro 0) 
\end{equation}
for as
\(
\neg\neg x > 0 \leftrightarrow x > 0
\)
and
\(
\neg\neg x > 0 \leftrightarrow \neg\neg x \gro 0
\),
the formula between the brackets is equivalent to
\(
\neg \neg x \gro 0  \rightarrow x \gro 0
\),
which is equivalent to an instance of PEM~\citep[p.252]{Brouwer1925E}
which is consistent~\citep{Brouwer1908C},
and so~\eqref{L128} is contradictory.%
\footnote{Likewise,
\(
\exists x\in \numreal\neg(x \neq 0 \rightarrow x\apart 0) 
\)
is contradictory,
because of~\eqref{L030}
and the consistency of
\(
\neg \neg β  \apart γ   \rightarrow β  \apart γ 
\).}
A strong counterexample,
if there is one, 
should therefore involve a property specific to
universal quantification over the real numbers.
For this,
Brouwer had the principle of Weak Continuity for Numbers%
\footnote{Thus named by Troelstra and Van Dalen~\citep[p.208–209]{Troelstra.Dalen1988}; 
in Brouwer’s writings it remained nameless.}~[\citealp[p.13]{Brouwer1918B}; \citealp[p.189]{Brouwer1924D2}; \citealp[p.253]{Brouwer1925A};~\citealp[p.63]{Brouwer1927B}; \citealp[p.15]{Brouwer1954A}]
and the Fan Theorem.

\begin{pri}[Weak Continuity for Numbers]\label{L076}
\begin{equation}\tag{WC-N}
\forall α\exists x\in \numnat A(α,x)
\rightarrow
\forall α\exists m \in \numnat\exists x\in \numnat\forall β(\bar{α}m=\bar{β}m \rightarrow  A(β,x))
\end{equation}
\end{pri}
where
\(\bar{α}m\)
stands for the initial segment
\(\seq{α(1), α(2), \dots, α(m)}\).

Brouwer never gave an explicit justification of WC-N,
which he must have had;
I refer to~\citep{Atten.Dalen2002}
and the later~\citep[ch.7]{Atten2007}
for discussion,
a justification,
and further
references.

WC-N can be strengthened to
\begin{pri}[Continuity for Numbers]\label{L140}
\begin{equation}\tag{C-N}
\forall α \exists x \in \numnat A(α, x) 
\rightarrow
\exists F \in K_0\forall α  A(α,F(α))
\end{equation}
\end{pri}
An informal justification runs as follows.%
\footnote{For details,
see Kleene and Vesley~\citep[p.71–73]{Kleene.Vesley1965},
whose term is
‘Brouwer’s Principle for Numbers’,
and Troelstra and Van Dalen~\citep[p.211–212]{Troelstra.Dalen1988},
who coined the term 
‘(strong) continuity for numbers’.
One also finds the principle referred to
as
‘\(\forall α\exists x\)-continuity’
(e.g., by Myhill in~\citealp[p.175]{Kreisel1967b}),
which is confusing 
(to some)
because that is the quantifier combination
in the antecedent of WC-N as well.
But
C-N may be held to lay greater claim to the term,
because it gives fuller expression
to the intuitionistic meaning of that combination.}
The quantifier combination
\(\forall α \exists x\)
entails that 
\(x\)
can be constructed from
\(α\)
by a method,
whose existence can be made explicit.
By WC-N, 
for a given
\(α\)
a value for
\(x\)
can be constructed from an initial segment of
\(α\).
It can be assumed that the method will
assign the same 
\(x\) 
to extensions of that segment, 
as additional information should not change its outcome;
and the constructivity of the method
should entail that it is decidable 
whether a given initial segment is long
enough for the method to produce its output. 
The methods for which these 
assumptions hold can be represented by a class 
of continuous functionals
\(K_0\) 
(which can be defined
inductively).
Alternatively,
C-N
can be proved
from WC-N,
monotone bar induction,
and
AC-NF~[\citealp[p.65–66]{Dummett2000}].%
\footnote{Thanks to Joan Moschovakis for reminding me of this fact
and for the reference to Dummett.
See also~\citep[p.333–334, 5.6.3(ii)]{Kreisel.Troelstra1970}.
For AC-NF,
see Principle~\ref{L141}
below.
For a statement of monotone bar induction,
see,
e.g.,~\citep[p.63]{Dummett2000}.
The intuitions that justify any of the forms of bar induction are not
simpler than those appealed to in the intuitive justification
of C-N from WC-N.}

WC-N is used in Brouwer’s  proofs of the Fan Theorem from 1927.%
\footnote{\label{L142}In ‘\german{Über Definitionsbereiche von Funktionen}’
from 1927~\citep{Brouwer1927B},
Brouwer presents two proofs,
in both of which the Fan Theorem is a corollary
of the Bar Theorem.
A 
\textit{bar}
is a set of nodes in a tree such that every path through the tree
intersects it.
The Bar Theorem states that
if a tree contains a bar,
then it contains a bar that admits of a well-ordering.
One proof of the Bar Theorem is based
on a general induction principle indicated in footnote 7 of Brouwer’s paper,
the other,
in section 2 of the main text, 
on the insight that  proofs
of the hypothesis of the Bar Theorem,
when considered as mental objects,
can be put into a canonical form.
In both the bar is defined by an application of 
WC-N.
Brouwer should have included that in the statement of the theorem
(or decidability, 
uniqueness, 
or 
monotonicity of the bar);
for otherwise,
as Kleene has shown,
the theorem is false~\citep[p.87–88]{Kleene.Vesley1965}.
(On this account,
in Brouwer’s presentation of the Bar Theorem in 1954~\citep{Brouwer1954A},
where the bars are not taken to arise by applications of WC-N,
there is a gap.)
Kleene gives various correct formulations~\citep[p.54–55]{Kleene.Vesley1965}.
There is ample discussion
of Brouwer’s argument for the Bar Theorem 
based on canonical proofs~[\citealp{Parsons1967};~\citealp{Martino.Giaretta1981};~\citealp[section 3.2]{Dummett2000b};~\citealp[ch.4]{Atten2004};~\citealp{Veldman2006a}].}

\begin{dfn}[{[\citealp[p.4]{Brouwer1923A}; \citealp[p.143]{Brouwer1952B}]}]
A 
\textit{fan} 
is a finitely branching tree.
An infinite path through a fan is called an 
\textit{element} 
of it.
\end{dfn}

When considered as a spread,
the unit continuum 
\([0,1]\)
is itself constructed as an infinitely branching tree,
the paths through which are choice sequences.
The following theorem shows that
it can in a sense be represented by a fan:

\begin{thm}[{{\citep[p.192]{Brouwer1924D2}}}]%
\label{L003}
There is a fan 
\(J\)
that coincides with the unit continuum 
\([0,1]\)
in that
every real number 
\(r \in J\) 
coincides with a real number 
\(s\in [0,1]\),
and every 
\(s \in [0,1]\) coincides with a real number 
\(r \in J\).
\end{thm}

\begin{thm}[{Fan Theorem [{\citealp[p.192]{Brouwer1924D2};~\citealp[p.66]{Brouwer1927B}}]}]
Let
\(F\)
be a method that assigns to every element 
\(e\)
of a fan
a 
(not necessarily unique)
number
\(F(e) \in \numnat\).
Then there exists
an
\(m \in \numnat\),
dependent on
\(F\)
but not on
\(e\),
such that
\(F(e)\)
depends on only the first
\(m\)
chosen values in 
\(e\).%
\footnote{Kleene~\citep[p.59]{Kleene.Vesley1965}
observes that classically
this theorem is false,
because of its dependence on WC-N.
He states a version
that is also classically true:
if there is a decidable bar in a fan,
then there is a uniform bound on the depth of the paths to the bar
(contraposition of König’s Lemma,
which itself is intuitionistically incorrect).
Current discussion of the Fan Theorem in non-intuitionistic constructivism is
concerned with variants of this latter version~[\citealp[section 3.2]{Dummett2000};~\citealp[section 4.1]{BridgesPalmgren2013};~\citealp{Veldman2014}].}
\end{thm}

The conclusion of the Fan Theorem 
should be seen not only in light of the contrast between
the presence and absence of a
(uniform)
bound on the determination of 
\(F(e)\),
but also of the contrast between
a determination of
\(F(e)\)
from only values in
\(e\)
and
a determination of
\(F(e)\)
from,
or also from,
restrictions that the Creating Subject may have imposed on how these values
are chosen.

\begin{prf}[of Theorem~\ref{L004}]\label{L074}
This is Brouwer’s incorrect proof of
1949;
in this presentation we will use mostly Brouwer’s own notation.
Let
\(γ\)
be a drift with kernel 
\(0\)
and
counting numbers
\mathlist{2^{-1}, 2^{-2}, 2^{-3},} \dots,
and
let 
\(f\)
be an arbitrary point of
the fan 
\(J\)
of Theorem~\ref{L003}.
A general condition imposed on the choices in
\(f\)
is that
\(f_n\)
must be chosen
after
the choice of
\(c_{n}(γ)\)
but before that of
\(c_{n+1}(γ)\).
This is to ensure that the Creating Subject
continues the construction of 
the direct checking number
as it continues constructing
\(f\);%
\footnote{Without such a condition
the Subject
would not be
obliged to do so,
on account of
its creative freedom.}
it also ensures that the choice of
\(c_{n+1}(γ)\)
is always informed by 
an initial segment
of 
\(f\)
of length
exactly
\(n\),
which means that,
in appropriate cases,
from a hypothesis about
\(c_{n+1}(γ)\)
we may infer a property
of
that
initial segment.

Let
\(α_f\)
be the proposition
\(f\in\numrat\),
\(ρ\)
the species
of the sequences
\(R(γ,α_f)\)
and
\(δ\)
the species
of the direct checking numbers
of
\(γ\)
through
\(α_f\),
\(D(γ,α_f)\).

Assume,
towards a contradiction,
that 
\(x > y \leftrightarrow x \gro y\).

It follows from
\(\forall f(\neg(\neg α_f \wedge \neg\neg α_f))\)
and 
the definition of a direct checking number
that
\(\forall f(D(γ,α_f) > 0)\).
By the assumption,
then also
\(\forall f(D(γ,α_f) \gro 0)\),
i.e.,
\(\forall e \in δ \exists v \in \numnat (e \gro 2^{-v})\).
Call
such a
\(v\)
depending on
\(e\),
which will not be unique,
\(n(e)\).
It follows that for each
\(f\),
\(D(γ,α_f) \in\numrat\)
can be proved
before 
\(n(e)\)
choices
have been made in 
\(D(γ,α_f)\), 
hence 
the proposition
\(f\in\numrat\)
can be decided
before 
\(n(e)\)
choices
have been made in 
\(R(γ,α_f)\),
and
hence,
by the general condition on the choices in 
\(f\),
the proposition
\(f\in\numrat\)
can be decided
before
\(n(e)\)
choices
have been made in 
\(f\).

This situation is reflected in the fan
\(J\).
From every 
\(p = R(γ,α_f) \in ρ\)
a real number
\(D(γ,α_f)\in [0,1]\)
can be obtained;
by Theorem~\ref{L003}, 
an element
\(f(p) \in J\)
can then be found that
coincides with the latter.
By the above,
the proposition
\(D(γ,α_f) \in \numrat\)
can be proved before
\(n(e)\)
choices have been made in
\(p\),
so there correspondingly is an
\(n(p) \in \numnat\)
such that 
the proposition
\(f(p) \in\numrat\)
can be proved before
\(n(p)\)
choices have been made in 
\(f(p)\).
(While coinciding,
\(p\)
and
\(f(p)\) 
may 
converge in different manners,
so that
if
for some predicate
\(P\)
one can prove 
\(P(p)\)
as soon as
\(\seq{p_1, \dots, p_k}\)
have been chosen,
the 
number of
choices in
\(f(p)\)
needed
to prove
\(P(f(p))\)
may be 
different.)
The number
\(n(p)\)
depends,
through
\(n(e)\),
on
\(f\),
and
the proposition
\(f(p) \in\numrat\)
is equivalent to
\(f \in\numrat \vee f \not\in\numrat\).
Thus,
for each
\(f \in J\),
an
\(n(f) \in \numnat\)
can be found 
such that
the proposition
\(f \in\numrat \vee f \not\in\numrat\)
can be proved
before 
\(n(f)\)
choices have been made in 
\(f\).

By the Fan Theorem,
the species
\(\{n(f) \mid f \in J \}\)
has a maximum element
\(m\).
Among the sequences admitted to the fan
\(J\)
are
sequences 
\(g\)
which up to
choice
\(m\)
have been chosen freely,
subject 
to the general
restriction on elements
of
\(J\)
and no others.
Therefore,
\(g \in \numrat \vee g \not\in \numrat\)
cannot be proved by stage
\(m\).
But the conclusion of the previous paragraph,
which depends on
the assumption
that
\(x > y \rightarrow x \gro y\),
implies that it can,
and hence this assumption
has led to a contradiction.
\end{prf}

\subsection{Myhill’s objection to, and repair of, the argument from 1949}\label{L067}

In July 1966,
John Myhill wrote a letter to Brouwer claiming
that the appeal to the Fan Theorem in Proof~\ref{L074}
is not correct~\citep[p.465–466]{Dalen2011}.%
\footnote{This objection was not yet made in Myhill’s contribution
to the discussion printed after
Kreisel’s seminal paper~\citep{Kreisel1967b};
the preface to the volume in which it appeared~\citep{Lakatos1967}
is dated July 1966 as well.}
No 
(draft)
reply from Brouwer
is known.%
\footnote{Brouwer in his note to that letter errs in taking it to refer to ‘Points and spaces’,
which contains the argument discussed in subsection~\ref{L069} below, 
but not one using the Fan Theorem.}
Myhill remarked on the problem in print in 1967~\citep[p.296]{Myhill1967}.

\begin{quotation}
In your version,
if I understand it correctly,
the proof runs as follows.
For each real number
\(α \in [0,1]\),
the real number
\(φ(α)\)
is defined as follows:
as long as the creating subject has not judged the proposition
‘\(α\) is rational’,
let
\(φ(α)(n)=1/{2^n}\);
if at the 
\(k\)th step
(after
\(k\)
choices for 
\(α\))
he decides that 
\(α\) is rational
or irrational,
let
\(φ(α)(k+q)=1/{2^k}\)
for all
\(q\).
Then
\(φ(α)\)
cannot be
\(0\),
for if it were
\(α\)
could be neither rational nor not rational.
All this is quite clear.
The difficulty lies in the second half of the proof.

\elide

It is the application of the fan theorem which I question here.
As I understood it,
the fan theorem applies only to those cases
in which to every free choice sequence
\(α\) 
belonging to the finitary spread
\(F\)
we can assign a natural number
\(n_{α}\)
\emph{using only the values}
\mathlist{α(0), α(1), α(2),} \dots. 
The proof of the fan theorem,
it seems to me,
depends essentially on this condition
(which is met in the usual mathematical cases:
for instance in the theorem,
that I used above,
that
\([0,1]\)
has no detachable subspecies.)
But it is not met in the situation to which you apply the fan-theorem
here,
because in computing the 
\(n\) 
from the \(α\) 
one is allowed to use also the values of
\(φ(α)\),
which may depend not only on
\(α\) 
but also on what restrictions have been placed on
\(α\),
and on what properties of
\(φ(α)\)
the creating subject may have inferred from these.
\end{quotation}

In terms of Brouwer’s proof,
the objection amounts to this.
In the construction of
\(n(f)\)
from
\(f(p)\),
and hence from
\(f\),
the Creating Subject may use, 
besides the values in
\(f\),
the values in 
\(D(γ,α_f)\),
and the latter
may depend
on restrictions that it
has imposed on
\(f\);
whereas application of the Fan Theorem would require that
\(n(f)\)
can be calculated from
only values in 
\(f(p)\),
and hence from only values
in
\(f\).

To my mind,
Myhill is right about this.
Under the assumption of 
\(x > y \leftrightarrow x \gro y\),
the inference from 
\(\forall f(D(γ,α_f) > 0)\)
to
\(\forall f(D(γ,α_f) \gro 0)\),
that is,
to the existence of a method that transforms
any 
\(f\)
into a proof of
\(D(γ,α_f) \gro 0\),
entails
no condition on the information
about \(f\)
that this method may depend on
to arrive at that proof;
but it is from that proof that,
via 
\(n(e)\)
and
\(n(p)\),
\(n(f)\)
is constructed.
On the other hand,
use of the Fan Theorem
does impose
the condition Myhill states,
because it is required by
WC-N,
on which the proof of that theorem depends.

Brouwer and Myhill both appeal to the impossibility of 
deciding the question of rationality of a real number from
a freely chosen initial segment of a given length without restrictions;
but Brouwer
in his reductio argument
applies the Fan Theorem
to steer toward it,
whereas,
as Myhill observes,
it should have kept him from applying the Fan Theorem in the first place.%
\footnote{Also Heyting’s version of Brouwer’s proof~\citep[p.117-118]{Heyting1956}
would be affected by Myhill’s criticism,
as Myhill indicates~\citep[p.175]{Myhill1968}.
Note that in the third edition of Heyting’s book
this matter is not brought up~\citep[p.121–122]{Heyting1972},
but Heyting 
in his note 4 to the reprint of Brouwer’s paper
in the
\textit{Collected Works}~\citep[p.603]{Brouwer1975}
in effect concurs with Myhill’s objection,
without a reference, although Myhill’s papers are included in the bibliography,
and without suggesting a way to save the theorem.
The remainder of Heyting’s note there contains a related remark from Brouwer, 
part of which will be discussed,
in a different context, 
in subsection~\ref{L073} below.}
(This situation may arise, 
more generally, 
whenever an application of a theorem overlooks
a condition on the possibility to do so.)

In his letter to Brouwer,
Myhill goes on to propose a repair.
Suppose that
for every
\(f \in J\)
the proposition
\(f\in\numrat\)
could be decided
before 
\(n(e)\)
choices
have been made in
it:
\begin{quote}
Now in my formalism this is immediately contradictory,
because it implies that the species of all real numbers
in
\([0,1]\)
would be split up into the rational[s] 
and the irrationals,
q.e.d.
\end{quote}
The letter unfortunately does not present that formalism itself.
For the contradiction, 
Brouwer’s Negative Continuity Theorem
(see Appendix~\ref{L057})
would suffice,
but it is likelier that Myhill was thinking of
Brouwer’s theorem
that fully defined functions on 
\([0,1]\) 
are uniformly continuous~\citep{Brouwer1927B}.
That theorem is proved 
from the Bar Theorem,
by way of deriving the Fan Theorem
and then applying the latter in a way
that respects the condition Myhill refers to in his letter.
In 1963 Kreisel had proved the Bar Theorem from three 
continuity axioms~\citep[p.IV-20]{Kreisel1963};
Myhill preferred that proof to Brouwer’s own two%
\footnote{See footnote~\ref{L142} above.},
and included 
these axioms in his published formalism
of 1968~\citep[p.166–167]{Myhill1968}.

\subsection{An argument from 1954}\label{L069}

The occasion for the following counterexample,
from 
‘Points and spaces’~\citep[p.4]{Brouwer1954A},
is a discussion of the
Principle of the Excluded Middle.
In constructive mathematics
any open problem  is a weak counterexample
to PEM,
but this argument shows that,
as long as there are open problems at all,
in whatever domain of mathematics,
a garden-variety instance of quantified PEM
is not valid either:

\begin{wce}[{{\citep[p.4]{Brouwer1954A}}}]\label{L104}
There is no hope of showing that
\begin{equation}
\forall x\in\numreal(x \in \numrat \vee x \not\in \numrat)
\end{equation}
\end{wce}

Brouwer’s theorem from 1927 
that all full functions on the unit continuum are uniformly continuous
entails the stronger result that that statement is contradictory,
but uses the Bar Theorem as a lemma.
Veldman has shown that such machinery 
is not needed for the weaker theorem
that all full functions on the unit continuum are continuous~\citep{Veldman1982}.
Brouwer did in his paper also establish,
without the Bar Theorem,
the again weaker Negative Continuity Theorem:
Every full function on the unit continuum is negatively continuous~\citep[p.62]{Brouwer1927B}. 
That
would still yield this weak counterexample.
Some details of the Negative Continuity Theorem
are presented in Appendix~\ref{L057},
so as to make some remarks on
the resemblance between Brouwer’s proof of it
and his 1954 plausibility argument for Weak counterexample~\eqref{L104};
the former may well have inspired the latter,
which
runs as follows.

\begin{pargmt}\label{L068}
Let
\(A\)
be a proposition
that is at present untestable,
and
\(γ\)
a
drift
with
rational counting numbers
\(\seq{c_n(γ)}_n\)
and an irrational
kernel
\(c(γ)\).
Then truth of 
\(A\)
and rationality of
the conditional checking number
\(C(γ,A)\)
are equivalent:
the choices in 
\(C(γ,A)\)
become a fixed rational
as soon as
a construction for
\(A\)
has been found,
and only then.
So as long as 
the proposition
\(A\)
is not testable,
the proposition
\(C(γ,A) \in \numrat\)
is not testable,
and,
as decidability implies testability,
the proposition
\(C(γ,A) \in \numrat \vee C(γ,A) \not\in \numrat\)
is not provable.
\end{pargmt}
For the discussion of \bks{+}{}
it will be not the counterexample itself that is important,
but that the observation
that
\(A \leftrightarrow C(γ,A) \in \numrat\)
is also made in Brouwer’s own text
(see subsection~\ref{L072} for the relevance of this).

\section{BKS and CS}

It is characteristic of the philosophical instinct
of Kripke,
Kreisel,
and Myhill
that they decided to take Brouwer’s arguments at face value and analyse them.
In the remark that is the motto for the present paper,
Kreisel says that
‘Brouwer’s views may be wrong or crazy
(e.g.~self-contradictory),
but one will never find out without looking
at their more dubious aspects’~\citep[p.159]{Kreisel1967b}.
I read this with the emphasis on ‘looking’.
The attitude Kreisel is countering here is that
of Heyting and Kleene,
described on p.\pageref{L132} above.
This is also seen in Kreisel’s
review  of Kleene and Vesley’s
\textit{Foundations of Intuitionistic Mathematics}:
‘Chapter \textsc{iv} is handicapped by the author’s obvious wish to avoid the dubious notions used in Brouwer’s refutation of 
\(\ast\) [= 
\(\forall α(\neg\neg\exists x(αx=0) \rightarrow \exists x(αx=0))\)]
instead of studying them’~\citep[p.261]{Kreisel1966}.
The next two subsections discuss the respective analyses
of Brouwer’s Creating Subject arguments 
by
Kripke and Kreisel.

\subsection{BKS}\label{L085}
The schema Kripke introduced around 1965 
to reconstruct Brouwer’s Creating Subject arguments
is
\begin{equation}
\exists α\left[%
\begin{gathered}
\forall n(α(n)=0  \vee α(n)=1) \\
\wedge \\
\forall n(α(n)=0) \leftrightarrow \neg A\\
\wedge\\
\exists n(α(n)=1) \rightarrow A
\end{gathered}
\right]
\tag{\bks{-}{}}
\end{equation}
There is also a strong version,
\bks{+}{}
(see below).%
\footnote{Note that the notation 
\bks{}{},
\bks{+}{},
\bks{-}{}
does not correspond to 
KS, KS+, KS-
as used by Dragálin~\citep[p.132]{Dragalin1988}:
his 
KS
is 
\bks{-}{},
his KS+
is
\bks{+}{},
his KS-
is only the bi-implication in \bks{-}{}.
I use 
\bks{}{}
as a general term.}
The traditional tag for the schema
has 
‘KS’ instead of
‘BKS’;
but,
as will be explained in section~\ref{L032},
the latter,
for 
‘Brouwer-Kripke Schema’
is more appropriate.

Thus formulated,
a witness for
\(\exists n(α(n)=1)\)
need not be unique.
But it is means no loss of generality\label{L045}
to assume that it is:
if 
\(α\)
possibly contains more than one 
\(1\),
then 
\(α^\ast\)
contains at most one,
while preserving the properties 
guaranteed by
\bks{-}{}:
\begin{equation}
α^\ast(n) = 
\begin{dcases*}
α(n) & \text{if  
\(\forall m < n(α(m)=0 )\)}\\
0 & \text{otherwise}
\end{dcases*}
\end{equation}

Since the assumption of uniqueness seems not often needed,
I will not treat it as the default,
and instead be explicit when invoking it.
If 
\bks{}{}
is understood as formulated by 
Kripke in his  Amsterdam lecture of 2016 
(see
subsection~\ref{L031} below),
then it will not be unique,
and the same holds for 
Myhill’s formulation~\citep[p.151]{Myhill1970}.

Kripke did not publish the schema,
but 
it circulated widely~\citep[p.336]{Myhill1968b}.
In particular, 
Kripke stated 
it
in a letter to Kreisel~\citep[footnote 8]{Kripke2018},
who went on to introduce
the alternative
\cs{}{},
which will be discussed in subsection~\ref{L125}.
The first references to BKS in print are
made in 1967, 
by 
Kreisel~\citep[p.174]{Kreisel1967b}
and
Myhill~\citep[p. 295]{Myhill1967}.

As a special case of
\bks{-}{},
one obtains for species 
\(X\)~\citep[p.128]{Kreisel1970a}: 
\begin{equation}
\forall x\exists α\left[%
\begin{gathered}
\forall n(α(n)=0  \vee α(n)=1) \\
\wedge \\
\forall n(α(n)=0) \leftrightarrow x \not\in X\\
\wedge\\
\exists n(α(n)=1) \rightarrow x \in X
\end{gathered}
\right]
\notag
\end{equation}
The quantifier order can be reversed,
by coding the sequences obtained for the various values of
\(x\)
into one.
This can be done because
the  infinitely many infinite sequences
\(\seq{α(x)}_k\)
may be constructed step by step
in a zigzag manner:
\mathlist{α(0)(0),
α(1)(0),
α(0)(1),
 α(2)(0),
\dots}.
Using the pairing function
\(j\)
to code 
‘the 
\(k\)-th value of the 
\(x\)-th sequence’
into one number,
define the sequence
\(β\)
by
\(β(j(x,k)) = α(x)(k)\).
This yields
\begin{equation}
\exists β\forall x\left[%
\begin{gathered}
\forall n(β(n)=0  \vee β(n)=1) \\
\wedge \\
\forall n(β(j(x,n))=0) \leftrightarrow x \not\in X\\
\wedge\\
\exists n(β(j(x,n))=1) \rightarrow x \in X
\end{gathered}
\right]
\tag{\bks{-}{S}}
\end{equation}
As the way
\(β\)
is constructed from the
sequences
\(α(x)\)
does obviously not depend on the exact values that the latter may take,
one concludes
to the correctness of a more general
choice principle
[e.g.~\citealp[p.14]{Kleene.Vesley1965}; 
\begin{pri}\label{L141}
\begin{equation}
\forall n\exists αA(n,α) \rightarrow \exists β\forall n A(n,{(β)}_n)
\tag{AC-NF}
\end{equation}
with
\({(β)}_n \coloneqq λm.β(j(n,m))\).
\end{pri}

For inhabited species
\(X\),
one derives from
\bks{-}{S}
\begin{equation}\label{L087}
\exists f\forall x\left[%
\begin{gathered}
\neg\exists n(f(n)=x) \rightarrow x \not\in X\\
\wedge\\
\exists n(f(n)=x) \rightarrow x \in X
\end{gathered}
\right]
\tag{\bks{-}{SF}}
\end{equation}
by
putting
\(f(j(x,k))=a\)
if
\(β(j(x,k))=0\),
and
\(f(j(x,k))=x\)
if
\(β(j(x,k))=1\).
\bks{-}{SF}
was first published 
(as the general weak version of
\bks{}{})
by Kreisel~\citep[p.128]{Kreisel1970a},
who does not derive it from
\bks{-}{S}
but
remarks that
it 
follows from 
\cs{}{};
see Theorem~\ref{L086}
and its proof
below.

As an example of the use of \bks{-}{},
here is an alternative to Plausibility argument~\ref{L005} for
Weak counterexample~\ref{L013}.
\begin{pargmt}\label{L020}
(In the style of~\citep[p.245]{Hull1969})
Let 
\(A\)
be a proposition 
that is at present not
testable.
Applying \bks{-}{} to
\(A\)
and to
\(\neg A\)
gives
\begin{align}
&  \exists α[\forall n(α(n)=0) \leftrightarrow \neg A)
\wedge
\exists n(α(n)=1) \rightarrow A]\label{L016}\\
& \exists β[\forall n(β(n)=0) \leftrightarrow \neg\neg A)
\wedge
\exists n(β(n)=1) \rightarrow \neg A]\label{L017}
\end{align}
Define the real numbers
\(r\)
\begin{equation}
r(n) = 
\begin{dcases*}
0 & \text{if  \(\forall m \leq n(α(m)=0)\)}\\
2^{-k}\phantom{-} & \text{if  \(k \leq n \wedge \forall m < k(α(m)=0) \wedge α(k)=1\)}
\end{dcases*}
\end{equation}
and
\(s\)
\begin{equation}
s(n) = 
\begin{dcases*}
0 & \text{if  
\(\forall m \leq n(β(m)=0)\)}\\
-2^{-k} & \text{if   \(k \leq n \wedge \forall m < k(β(m)=0) \wedge β(k)=1\)}
\end{dcases*}
\end{equation}
Finally,
define the real number
\(t\) 
by
\begin{equation}
t(n) = r(n) +s(n)
\end{equation}
By~\eqref{L016}
and~\eqref{L017},
\begin{equation}
\exists n(α(n)=1) \wedge \exists n(β(n)=1) \rightarrow A \wedge \neg A
\end{equation}
hence
\begin{align}
& \exists n(α(n)=1)  \rightarrow  \forall n(β(n)=0)  \rightarrow t > 0\label{L014}\\
& \exists n(β(n)=1)  \rightarrow  \forall n(α(n)=0)  \rightarrow t < 0\label{L015}
\end{align}
and,
by contraposing both and using
\(t=0 \rightarrow \neg(t > 0) \wedge  \neg(t < 0)\),
\begin{equation}
t=0 \rightarrow \forall n(α(n)=0) \wedge \forall n(β(n)=0)
\end{equation}
In turn,
by~\eqref{L016} %
and~\eqref{L017}, %
\begin{equation}
\forall n(α(n)=0) \wedge \forall n(β(n)=0) \rightarrow \neg A \wedge \neg\neg A
\end{equation}
so
\begin{equation}
t \neq 0
\end{equation}
By 
def.~\(<\),~\eqref{L015} and then~\eqref{L017},
\begin{align}
& t < 0   \rightarrow   \neg(t > 0 )   \rightarrow   \neg\forall n(β(n)=0)    \rightarrow \neg\neg\neg A \rightarrow \neg A\label{L022}\\
\intertext{and similarly, from~\eqref{L014} and then~\eqref{L016},}
& t > 0   \rightarrow   \neg(t < 0 )   \rightarrow   \neg\forall n(α(n)=0)  \rightarrow \neg\neg A
\end{align}
which together yield
\begin{equation}
t < 0  \vee  t > 0 \rightarrow 
\neg A \vee \neg\neg A
\end{equation}
Hence,
under the hypothesis that 
\(A\) 
cannot be tested,
\(t < 0 \vee t > 0\) cannot be proved,
and,
by implication, 
neither
can
\(t \klo 0 \vee t \gro 0\).
\end{pargmt}

A reconstruction of Argument~\ref{L021},
yielding a weak counterexample only to
\(t \klo 0 \vee t \gro 0\),
would proceed in the same way,
but,
whereas in~\eqref{L022} we had to make the step
\( t < 0   \rightarrow   \neg(t > 0 )\),
we can now,
at the corresponding place,
reason that
\(t \klo 0   \rightarrow   \exists n(α(n) = 1) \);
and in the case of
\(\gro\)
this yields a stronger conclusion.
So
instead we have
\begin{align}
& t \klo 0   \rightarrow   \exists n(β(n) = 1) \rightarrow \neg A\\
\intertext{and}
& t \gro 0   \rightarrow   \exists n(α(n) = 1) \rightarrow  A
\end{align}
which together yield
\begin{equation}
t \klo 0   \vee  t \gro 0 \rightarrow 
A \vee \neg A \rightarrow \neg A \vee \neg \neg A
\end{equation}
Hence,
under the hypothesis that 
\(A\) 
cannot be tested
(or not be decided),
\(t \gro 0 \vee t \klo 0\) cannot be proved.

Reconstructions
using 
\bks{-}{}
of other Creating Subject arguments
as given by Heyting in his
\textit{Intuitionism}~\citep{Heyting1956}
(and thus including,
of Brouwer’s arguments discussed above,
the one from 1949 but not that from 1954)
were carried out by Hull in 1969~\citep{Hull1969}.
Hull also presents an alternative argument for
Heyting’s claim that 
the virtual order of the continuum cannot be expected to be a pseudo-order
(subsection~\ref{L018} above).
While Hull appeals to
\bks{-}{} and continuity~\citep[p.244]{Hull1969},
he indicates 
that 
Heyting’s argument
can be reconstructed using only this
variant of 
\bks{-}{}
with an additional clause~\citep[p.246]{Hull1969}: 
\begin{equation}
\exists α\left[%
\begin{gathered}
\forall n(α(n)=0  \vee α(n)=1) \\
\wedge \\
\forall x(α(n)=0) \leftrightarrow \neg A\\
\wedge\\
\exists x(α(n)=1) \rightarrow A\\
\wedge\\
\forall x\exists y(t(y) > x) \rightarrow (\forall x(α(t(x))=0) \rightarrow A \vee \neg A)
\end{gathered}
\right]
\tag{\bks{-}{W}}
\end{equation}
where 
\(t\)
is an arbitrary term with infinite range
(the formalism had no variables for functions).
The new clause is an expression of the idea that,
since the Creating Subject does not know beforehand
where,
if anywhere,
a
\(1\)
will occur in
\(α\),
knowing that 
\(A\)
will not be implied by the presence of a
\(1\)
in 
\(α\)
at
a position of the form
\(t(x)\)
is only possible if
\(A\)
has been decided
already.
The informal explanation 
Hull provides is that
‘the only way we could know something will not be proved
on a Wednesday is that either it has already been proved
on some other day,
or that its negation has already been proved’~\citep[p.246]{Hull1969};
this also explains my choice of the tag.%
\footnote{Hull’s advisor Myhill
generalised
\bks{-}{W}
by using a decidable subspecies of choice sequences
instead of just a term containing
\(x\),
and called it
the ‘Never-on-Sunday schema’~\citep[p.158]{Myhill1970}.}
Now,
even if some proposition has already been proved on some other day,
the Creating Subject surely has the freedom to prove it again on Wednesday;
what is meant here clearly is that
the first time the Subject proves 
\(A\)
is not a Wednesday.
For the additional clause to reflect that idea accurately,
it must be assumed
that 
\(α\)
takes the value
\(1\)
at most once.
By the remark
on 
p.\pageref{L045} above, 
this assumption can be made safely.

Hull does not go on to give 
a reconstruction of
Heyting’s argument
with
\bks{-}{W},
but one is as follows.
\begin{pargmt}[for Weak counterexample~\ref{L046}]
Let 
\(A\)
be a proposition 
that is not yet testable.
Applying \bks{-}{W} to
\(A\)
and
\(t(y)=2y\),
together with the fact that
\(\forall x\exists y(2y > x)\)
is true in the domain of the natural numbers,
gives
\begin{equation}\label{L036}
\exists α\left[%
\begin{gathered}
\forall n(α(n)=0  \vee α(n)=1) \\
\wedge \\
\forall n(α(n)=0) \leftrightarrow \neg A\\
\wedge\\
\exists n(α(n)=1) \rightarrow A\\
\wedge\\
\forall n(α(2n)=0) \rightarrow A \vee \neg A
\end{gathered}
\right]
\end{equation}
and,
similarly,
\begin{equation}\label{L039}
\exists β\left[%
\begin{gathered}
\forall n(β(n)=0  \vee β(n)=1) \\
\wedge \\
\forall n(β(n)=0) \leftrightarrow \neg A\\
\wedge\\
\exists n(β(n)=1) \rightarrow A\\
\wedge\\
\forall n(β(2n+1)=0) \rightarrow A \vee \neg A
\end{gathered}
\right]
\end{equation}

As mentioned (p.\pageref{L045}),
we may assume that
\(α\)
and
\(β\)
take the value
\(1\)
at most once.

Define the real numbers
\(r\)
and
\(s\)
by
\begin{align}
& r(n) = 
\begin{dcases*}
2^{-n} & \text{if  
\(\forall m \leq n(α(m)=0 \wedge β(m)=0)\)}\\
2^{-k}\phantom{-} & \text{if  
\(k \leq n \wedge (α(k)=1 \vee β(k)=1)\)}
\end{dcases*}\\
 \intertext{and}
& s(n) = 
\begin{dcases*}
2^{-n}  & \text{if  
\(\forall m \leq n(α(m)=0 \wedge β(m)=0)\)}\\
2^{-k} & \text{if 
\(k \leq n \wedge Even(k) \wedge (α(k)=1 \vee β(k)=1)\)}\\ 
2^{-n} & \text{\parbox[c]{20em}{if  
\(\exists k \leq n (Odd(k) \wedge (α(k)=1\vee β(k)=1))\)}}
\end{dcases*}
\end{align}

Then
\(0 < r\) because 
\(\neg\forall m(α(m)=0 \wedge β(m)=0)\).

If
\(s \apart 0\),
then by the second clause
\(s=r\),
so 
\(r \apart 0 \rightarrow s=r\)
and therefore,
by
\(s < r   \rightarrow s \neq r\)
and~\eqref{L040},
\(s < r   \rightarrow s=0\).
By definition,
\(0 < s   \rightarrow s \neq 0\).
Combining the two gives
\begin{equation}\label{L035}
s < r \vee   0 < s  \rightarrow s=0 \vee s \neq 0 
\end{equation}

If
\(s=0\),
then
\begin{equation}\label{L037}
\neg \exists k (Even(k) \wedge α(k)=1 \vee β(k)=1)
\end{equation}
hence
\begin{equation}
\neg \exists k (Even(k) \wedge α(k)=1)
\end{equation}
so
\begin{equation}
\forall n(α(2n)=0)
\end{equation}
By the last clause in~\eqref{L036},
\begin{equation}
A \vee \neg A
\end{equation}
hence
\begin{equation}
\neg A \vee \neg\neg A
\end{equation}

If
\(s \neq 0\),
then
\begin{equation}\label{L038}
\neg\neg \exists k (Even(k) \wedge α(k)=1 \vee β(k)=1)
\end{equation}
and therefore
\begin{equation}
\neg \exists k (Odd(k) \wedge α(k)=1 \vee β(k)=1)
\end{equation}
by the assumption that
\(α\)
and
\(β\)
take the value
\(1\)
at most once.
So
\begin{equation}
\neg \exists k (Odd(k)  \wedge β(k)=1)
\end{equation}
and hence
\begin{equation}
\forall n(β(2n+1)=0)
\end{equation}
By the last clause in~\eqref{L039},
\begin{equation}
A \vee \neg A
\end{equation}
hence
\begin{equation}
\neg A \vee \neg\neg A
\end{equation}

Together,
the arguments for the cases 
\(s=0\)
and
\(s\neq 0\)
yield
\begin{equation}
s=0 \vee s \neq 0  
\rightarrow 
\neg A \vee \neg\neg A 
\end{equation}
which shows that,
as long as 
\(A\)
cannot be tested,
\(s=0 \vee s \neq 0\)
cannot be proved;
by~\eqref{L035},
this means that,
as long as 
\(A\)
cannot be tested,
\(0 < s \vee s < r \)
cannot be proved.
On the other hand,
\(0 < r\)
can be proved;
hence
\(0 < r \rightarrow 0 < s \vee s < r\)
cannot be proved,
and this means that
axiom~\eqref{L027}
does not hold for
\(<\).
\end{pargmt}

The same thought that  justifies
\bks{-}{W}
also  justifies a version
of
\bks{}{}
that I here call
\bks{-}{R}
and which was
proposed and defended by 
Erdélyi-Szabó
and
Scowcroft~\citep[p.1027–1028]{Erdelyi-Szabo2000}.
If,
upon proving
\(A\),
the Creating Subject chooses a random positive number,
then nothing can be said about that number
until 
\(A\)
has been proved:
\begin{equation}
\exists α\left[%
\begin{gathered}
\exists n(α(n) > 0) \rightarrow A\\
\wedge\\
\neg\exists n(α(n) > 0) \rightarrow \neg A\\
\wedge\\
\forall k>0(\neg\exists n(α(n) = k) \rightarrow A \vee \neg A)\\
\wedge\\
\forall k>0\forall n((α(n) = k) \rightarrow \forall m \geq n(α(m) = k))
\end{gathered}
\right]
\tag{\bks{-}{R}}
\end{equation}

There is a strengthening of 
\bks{-}{} 
that is not only
of mathematical but also of philosophical significance~[\citealp[p.295]{Myhill1967};~\citealp[p.96]{Troelstra1969}]:
\begin{equation}
\exists α\left[%
\begin{gathered}
\forall n(α(n)=0 \vee α(n)=1) \\
\wedge \\
\exists n(α(n)=1) \leftrightarrow A
\end{gathered}
\right]
\tag{\bks{+}{}} 
\end{equation}
In subsection~\ref{L031},
Kripke’s reason for accepting 
\bks{-}{}
but not
\bks{+}{}
will be discussed.%
\footnote{Note that it is not the case that in reconstructions of 
Creating Subject arguments in terms of \bks{+}{} and \bks{-}{},
for strong counterexamples 
the strong schema is used,
and for weak counterexamples the weak.}

As in the case of 
\bks{-}{},
for species 
\(X\)
there is~\citep[p.104]{Troelstra1969}: 
\begin{equation}
\exists β\forall x\left[%
\begin{gathered}
\forall n(β(n)=0  \vee β(n)=1) \\
\wedge \\
\forall n(β(j(x,n))=0) \leftrightarrow x \not\in X\\
\wedge\\
\exists n(β(j(x,n))=1) \leftrightarrow x \in X
\end{gathered}
\right]
\tag{\bks{+}{S}}
\end{equation}

Classically, 
\bks{+}{} 
is just a weak comprehension principle~\citep[p.74]{Dalen1982b}.
In the  setting of
classically defined models for formal intuitionistic theories,
a consistency proof of 
\bks{+}{}
was given by Scott~\citep{Scott1970b};
his
conjecture that 
\bks{-}{}
is weaker than 
\bks{+}{}
was proved,
again for a classically defined model,
by
Krol’~\citep{Krol1978}.

\subsection{CS}\label{L125}
The axiom schemata that make up
(what became known as)
the
‘Theory of the Creating Subject’ 
(CS)
are 
due to Kreisel~\citep{Kreisel1967b},
who had received a letter from Kripke
stating
\bks{-}{}~\citep[footnote 8]{Kripke2018}.%
\footnote{Kreisel showed his schemata
in a letter to Gödel of July 6, 1965~\citep[item  011182]{Godel.Papers}.} 
In 1969 Anne Troelstra, 
in his early and influential treatment of Kreisel’s schemata,
added a schema
(see below)
and spoke of
‘Brouwer’s theory of the creative subject’~\citep[ch.16]{Troelstra1969}
where 
Kreisel had used the term
‘thinking subject’.
Since Brouwer had indeed used the term
‘creative subject’,
or rather ‘creating subject’%
\footnote{See p.\pageref{L115}
above.},
but had not presented an explicit theory
such as Kreisel proposed,
overall I think it is best to speak of
‘the Kreisel-Troelstra Theory of the Creating Subject’%
\label{L138}%
.

The basic notion in Kreisel’s original schemata
was
\(Σ \vdash_n A\),
where
\(Σ\)
is a variable ranging over Creating Subjects 
and which means that
‘the (thinking) subject
\(Σ\)
has evidence for asserting
\(A\)
at stage
\(m\)’~\citep[p.159]{Kreisel1967b}.
Troelstra 
assumes the existence of only one Creating Subject,
as do most later presentations.
I will do the same here,
and
say something about the reasons why in
subsection~\ref{L105}.

Troelstra  writes
\(\vdash_n A\)
with the meaning
‘the creative subject has evidence for
\(A\)
at stage
\(m\)’~\citep[p.95]{Troelstra1969},
but he is sensitive~\citep[p.105–106]{Troelstra1969}
to a possible ambiguity
in the phrase
‘to have evidence for’:
must the Creating Subject be aware of this
at stage
\(m\),
or is evidence closed under,
say, 
trivial
(but unrealised)
deducibility?
It seems to me that,
to the extent that the aim is a reconstruction of Brouwer’s ideas,
a construction can only 
contribute to making
some proposition evident
once the Creating Subject has made the connection between the two in an appropriate act,
so that evidence is not closed under unrealised deductions (trivial or not).
Thus,
to use Troelstra’s phrase,
I take the Creating Subject to establish one conclusion at a time.%
\footnote{For Troelstra~\citep{Troelstra1969}
this question comes up because of its relation to the paradox in the Theory of the Creating Subject that he there devises.
See subsection~\ref{L077} below.}

Instead of the turnstile of Kreisel and early Troelstra,
too reminiscent of formal derivability relations,
I prefer to use the propositional operator 
\(\csop_{n}\)
as in
Troelstra and Van Dalen’s~\citep{Troelstra.Dalen1988},
and understand
\(\csop_{n}A\)
as
‘By stage 
\(n\) the Creating Subject has made
\(A\) evident’.

Like BKS,
CS comes in a weak and a strong version.
In the original presentation~\citep{Kreisel1967b}
there is only the
weak version
\cs{-}{}.%
\begin{equation}
	\forall n(\csop_{n}A \vee \neg
	\csop_{n}A)
	\tag{\cs{-}{1}}
\end{equation}
That is,
for any stage,
it is decidable for the Creating Subject whether 
by that stage it has
made
\(A\)
evident.
\begin{equation}
	\forall n\forall m(\csop_{n}A \rightarrow
	\csop_{n+m}A) \tag{\cs{-}{2}}
\end{equation}
The Creating Subject never forgets what it has made evident.
(This axiom was added by Troelstra~\citep[p.95]{Troelstra1969}.)
\begin{equation}
(A
\rightarrow
\neg\neg\exists n\csop_{n}A)
\wedge
(\exists n\csop_{n}A
\rightarrow
A)
\tag{\cs{-}{3}}
\end{equation}
If
\(A\)
is true,
then it is not possible that the Creating Subject will never make it evident;
and the Creating Subject makes no mistakes.
Kreisel named
\cs{-}{3}
‘the principle of Christian Charity’.

In the stronger version~\citep[p.96]{Troelstra1969},
\cs{+}{1}~= \cs{-}{1} 
and 
\cs{+}{2}~= \cs{-}{2},
but
\cs{-}{3} is strengthened to
\begin{equation}
\exists n\csop_{n}A
\leftrightarrow
A  
\tag{\cs{+}{3}}
\end{equation}
A proposition 
\(A\) 
is true if and only if the Creating Subject
has made
\(A\)
evident by some stage.

Splitting \cs{-}{3}:
\begin{gather}
A
\rightarrow
\neg\neg\exists n\csop_{n}A
\tag{\cs{-}{3a}}\\
\exists n\csop_{n}A
\rightarrow
A
\tag{\cs{-}{3b}}
\end{gather}
and \cs{+}{3}:
\begin{gather}
A
\rightarrow
\exists n\csop_{n}A
\tag{\cs{+}{3a}}\\
\exists n\csop_{n}A
\rightarrow
A
\tag{\cs{+}{3b}}
\end{gather}

In  terms of \cs{-}{},
Plausibility argument~\ref{L005}
for Weak counterexample~\ref{L013}
could be recast as follows: 
\begin{pargmt}\label{L049}
Let 
\(A\) 
be a proposition 
that is, 
at present,
not
testable.
The 
choice
sequence 
\(α\) 
is defined as follows:
\begin{equation} %
α(n) =
\begin{dcases}
0 & \text{if  } \forall m \leq n(\neg \csop_{m}A \wedge \neg \csop_{m} \neg A)\\
2^{-k} & \text{if  } k \leq n \wedge \forall m < k(α(m)=0) \wedge \csop_{k}  A\\
-2^{-k} & \text{if  } k \leq n \wedge \forall m < k(α(m)=0) \wedge \csop_{k}  \neg A
\end{dcases}
\end{equation}
\end{pargmt}
With \cs{-}{2} and \cs{-}{3b} 
one shows that for no
\(n\)
both
\(\csop_{n}  A\)
and
\(\csop_{n}  \neg A\),
then
with \cs{-}{1}
that 
\(α\)
converges
and therefore determines a real number
\(r\),
and with
the contraposition of \cs{-}{3a}
that
\(r \neq 0\).
Then,
analogously to the 
reasonings~\eqref{L009} and~\eqref{L010},
\begin{equation}\label{L047}
\begin{aligned}
& r  < 0 & \text{(assumption)}\\
& \neg(r  > 0) & \text{(def.~\(<\))}\\
& \neg\exists n\csop_{n}A  & \text{(def.~\(r\))}\\
& \neg A & \text{(\cs{-}{3a})}\\
& \text{\(A\) has been tested} & \\
\end{aligned}
\end{equation}
and,
symmetrically,
\begin{equation}\label{L048}
\begin{aligned}
& r  > 0 & \text{(assumption)}\\
& \neg(r  < 0) & \text{(def.~\(>\))}\\
& \neg\exists n\csop_{n}\neg A & \text{(def.~\(r\))}\\
& \neg\neg A & \text{(\cs{-}{3a})}\\
& \text{\(A\) has been tested} & \\
\end{aligned}
\end{equation}

An immediate consequence of
\cs{+}{} 
that will be used further on
is:
\begin{thm}[{{\citep[p.103]{Troelstra1969}}}]\label{L086}
Let the
species 
\(X \subseteq \numnat\)
be inhabited.
Then it is 
enumerable by the 
Creating Subject.
\end{thm}
\begin{prf}
Let
\(j \colon \numnat^{\numnat} \rightarrow \numnat\)
be the Cantorian pairing function
\begin{equation}
j(n,k) = \frac{{(n+k)}^2 + n + 3k}{2}
\end{equation}
As
\(X\)
is inhabited,
the Creating Subject has at some stage 
\(k\)
proved that
\(a \in X\) 
for some
\(a\).
Without loss of generality,
it may be assumed that
\(k=0\)
and
\(a=0\),
i.e.,
\( \csop_{0} 0 \in X\).
Set 
\(f(0) = a = 0\)
and
define,
more generally,
\begin{equation}
f(j(n,m)) = \left\{
\begin{array}{ll}
f(0) & \text{if not }  \csop_{m} n \in X\\
n & \text{if \makebox[\widthof{not }]{}} \csop_{m} n \in X
\end{array}
\right.
\end{equation}
Then,
by
\cs{+}{3},
\( n \in X \leftrightarrow \exists m \csop_{m} n \in X\), 
and, 
by the definition of
\(f\),
\(\exists m\csop_{m} n \in X
\leftrightarrow \exists m(f(j(n,m)) = n)\);
hence,
\begin{equation}
n \in X \leftrightarrow \exists m(f(j(n,m)) = n)
\end{equation}
\end{prf}
For the Creating Subject
this function 
\(f\) 
is
computable,
as,
by
\cs{+}{1}, 
for any given 
\(m\)
and
\(n\),
\(\csop_{m} n \in X\)
is decidable.

Kreisel had become interested in
the Creating Subject
after his work on the Theory of Constructions~\citep{Kreisel1962,Kreisel1965};
the Theory of the Creating Subject
was meant
as an enrichment of that  theory~\citep[p.180]{Kreisel1967b}.
The particular use to which Kreisel wanted to put his Theory of the Creating Subject 
is related to the question of completeness of Heyting’s predicate logic.
Gödel had shown,
and in 1962 Kreisel had published~\citep{Kreisel1962b}, 
a theorem to the effect that
a completeness proof of Heyting’s predicate logic
(relative to a Tarski-style notion of validity but,
as was shown thereafter,
also to Beth and Kripke models)
would entail the validity of MP
in the form
\[
\forall α\neg\neg\exists x φ(α,x) \rightarrow \forall α\exists x φ(α,x)
\]
where φ is a primitive recursive predicate,
and α ranges over choice sequences
but not lawless ones.%
\footnote{The way Gödel
treated
negation in his argument
–~negation as absence of models~–
is decidedly un-Brouwerian,
as he will have been fully aware of.
Veldman
developed a notion of model
in which negation is handled differently,
which allowed him to
give an intuitionistically correct completeness proof
while blocking
Gödel’s argument.
But Veldman did not think that his semantics
(nor similar ones)
shed much light on the notion 
‘intuitionistically true sentence’~\citep[p.159]{Veldman1976}.}
So a rejection of this principle entails the incompleteness
of Heyting’s predicate logic.
As noted in subsection~\ref{L066},
Brouwer’s Creating Subject argument
for the conclusion 
\[
\neg\forall  x\in\numreal(x \neq 0 \rightarrow x \apart 0)
\]
entails a strong counterexample to
a strong version of MP,
and Kreisel was interested to see
if an argument like Gödel’s could
be used
also to derive
incompleteness of Heyting’s predicate logic
from that strong counterexample.
He did not succeed without invoking Church’s Thesis or something of
a similar nature~\citep[p.182]{Kreisel1967b},
which Brouwerian intuitionists do not accept
(see subsection~\ref{L034} below).

Be that as it may,
in preparing the final version of ‘Informal rigour’,
Kreisel discovered that he could
obtain Brouwer’s strong counterexample
\textit{without}
using the Theory of the Creating Subject~\citep[p.180–182]{Kreisel1967b}
(Vesley in 1968 devised a Schema for the purpose;
see subsection~\ref{L092} below).
It seems he then lost interest in the topic;
when in 1982 Van Dalen showed that the Theory of the Creating Subject
is conservative over Heyting Arithmetic~\citep{Dalen1982},%
\footnote{But Krivtsov~\citep{Krivtsov1999} 
sees the need for an alternative proof, which he supplies (for a more limited result).}
Kreisel,
in his 
\textit{\german{Zentralblatt}}
review~\citep{Kreisel1984},
took that
 to be 
‘further evidence for the mathematical sterility\label{L096} of CS’.
By the list in subsection~\ref{L034} below,
this, 
seems not borne out by the
(mostly later)
facts.
But I agree when he adds that considerations of 
\cs{}{} 
remain of interest 
‘as a philosophical object lesson’.

\subsection{Comparing BKS and CS}

A conspicuous difference between
\bks{}{}
and
\cs{}{} 
is that
the former
 is extensional,
whereas the latter is,
with its reference to 
the Creating Subject 
and the stages of its activity over time,
highly intensional.
On the other hand,
there seems to be no justification
of
\bks{}{}
except its derivation from
\cs{}{}.

These points were highlighted by
Myhill~\citep[p.295–296]{Myhill1967}
and 
Kreisel~\citep[p.128]{Kreisel1970a}.
Myhill
has a preference for
\bks{-}{}
over
\cs{-}{}
because
‘we think the loss of extensionality too high a price to pay in technical facility’~\citep[p.175]{Myhill1968}.%
\footnote{I am not sure why Myhill did not include 
\bks{+}{}
in his system,
the occurrence of which in Brouwer’s work
he had been the one to find (see subsection~\ref{L072});
perhaps because Brouwer made no further use of it.
Note that in later applications of \bks{}{},
it is the strong version that is used;
see  subsection~\ref{L091}.}
To illustrate this,
he points out 
that
the mistake 
in 
Brouwer’s argument from 1949
(which he refers to in Heyting’s version~\citep[p.118]{Heyting1956})
would have been avoided if
\bks{}{}
had been used
instead of
\cs{}{}
(see subsections~\ref{L066} and~\ref{L067} above).
Thus,
in his formalisation of intuitionistic analysis,
Myhill
chooses to adopt \bks{-}{} 
instead of
\cs{-}{}. %
But the reason is pragmatic;
he held that 
\cs{-}{}  
is correct,
and that it provides a 
‘deeper analysis’ 
than \bks{-}{}~\citep[p.296]{Myhill1967}.
Van Dalen~\citep[p.19]{Dalen1978a}
showed that
the addition of
\cs{+}{}
to a
theory of intuitionistic analysis
that contains
\bks{+}{}
is conservative.

\subsection{Applications beyond Brouwerian counterexamples}\label{L091}

In the discussion after Kreisel’s presentation of ‘Informal rigour’
in 1965,
Heyting had warned that 
the method of the Creating Subject
‘is not central in intuitionistic mathematics.
It can only be applied to show that certain propositions
of which nobody believed that they could be true,
are actually false’~\citep[p.173]{Kreisel1967b};
and we saw above 
(p.\pageref{L096})
that
eventually
also
Kreisel 
came to consider 
\cs{}{}
mathematically sterile.
Various theorems from
\bks{}{}
(and hence also from \cs{}{})
show
that
\cs{}{}
and
\bks{}{}
do have 
applications beyond Brouwerian counterexamples:
\begin{enumerate}

\item (Van Dalen) \bks{+}{} can be used to construct 
a topological model for the theory of species of natural numbers~\citep{Dalen1974}).

\item  (Van Dalen) \bks{+}{} can be used to construct an intuitionistic analogue
of the powerset of 
\(\numnat\)~\citep{Dalen1977}.%
\footnote{That paper was inspired by 
Troelstra~\citep[p.104]{Troelstra1969} and 
Kreisel~\citep[p.128]{Kreisel1970a}.}

\item (Burgess) 
If the
‘basic system’
of
principles common to intuitionistic and classical analysis
as defined by Kleene and Vesley~\citep[p.8]{Kleene.Vesley1965}
is extended with 
\bks{+}{},
a theorem can be obtained that is
classically equivalent to Souslin’s Theorem~\citep{Burgess1980}.%
\footnote{Gielen,
Veldman,
and De Swart~\citep[p.134]{Gielen.Swart.Veldman1981}
do not accept this proof because they argue that 
\bks{+}{}
can be justified only for 
‘determinate’
propositions
\(A\);
see subsection~\ref{L090} below.}

\item (Van Dalen) Using \bks{+}{},
one shows that
every negatively defined
dense subset 
\(X \subset \numreal\)
is unsplittable,
i.e.,
if 
\(X\)
is the disjoint sum of
\(Y\)
and
\(Z\),
either 
\(Y=X\)
or
\(Z=X\)~\citep{Dalen1999a}. 
(See also subsection~\ref{L092}
below.)

\item (Lubarsky, Richman, and Schuster) 
In the presence of countable choice,
\bks{+}{} 
is equivalent
to
‘Every open subspace of a separable space is separable’
and also to 
 ‘Every open subset of a separable metric space is a countable union of open balls’~\citep{Lubarsky.Richman.Schuster2012}.

\item (Kachapova) There is an intuitionistic formal theory SLP 
for higher types and lawless sequences,
together with a Beth model for it,
that includes
\cs{+}{}%
\footnote{The review in the
\textit{\german{Zentralblatt}}
opines that
‘The purist will be disappointed since SLP proves some principles, 
like weak continuity or the Kripke’s schema which are debatable in a predicative setting; 
also, 
the formal Church’s thesis is inconsistent with SLP’~\citep{Benini}.
By the main theme of the present paper,
from the Brouwerian point of view,
this is
all as it should be.}
and which is 
equiconsistent with TI, 
a subtheory of classical typed set theory that is stronger than
classical second-order arithmetic~\citep{Kachapova2015}.%
\footnote{This makes a considerable step towards fulfilling
(also) 
Brouwer’s prediction that Hilbert’s program 
for a constructive consistency proof of classical mathematics
would succeed,
although not based, 
as he thought, 
on the consistency of PEM,
but on his ideas about the Creating Subject.
Brouwer thus was not,
before Gödel,
sceptical about the possibility of 
such a consistency proof;
but he doubted its value
for the foundations of mathematics~[\citealp[p.252n4]{Brouwer1925E};~\citealp[p.377]{Brouwer1928A2};~\citealp[p.164]{Brouwer1929A}].}
\end{enumerate}
Here one should also mention 
Van Rootselaar’s
suggestion to use 
\cs{}{}
also for defining various notions of choice sequence
by characterising the Creating Subject’s knowledge of them~\citep[p.196]{Rootselaar1970},
the analysis 
by Friedrich and Luckhardt
of
the role of 
\bks{+}{}
in establishing certain uniformity principles~\citep{Friedrich.Luckhardt1980},
and
Schuster and Zappe’s
use
of versions of
\bks{+}{}
to classify countability statements~\citep{Schuster.Zappe2008}.

\section{A Brouwerian justification of CS and BKS}\label{L078}

\subsection{A remark on Brouwerian logic}\label{L103}

As my primary interest in this paper is in questions of justification and use
of Creating Subject arguments in a Brouwerian framework,
as opposed to other forms of constructivism,
and
\cs{}{}
and
\bks{}{}
were both intended to clarify those arguments,
I will insist here on construing
the logic  
in a Brouwerian way.
Brouwer’s conception of logic
can be found
in his dissertation from 1907~\citep[ch.3]{Brouwer1907}
and in 
‘The unreliability of the logical principles’
from 1908~\citep{Brouwer1908C, Atten.Sundholm2017}.
It differs considerably from
that embodied in the now ubiquitous
combination of
Natural Deduction and what has become known as 
‘the BHK interpretation’~\citep{Sundholm.Atten2008,Atten2009e}.
For the present purpose,
the main points are:
\begin{enumerate}
\item\label{L102}
Intuitionistic mathematics consists primarily 
in the act of effecting mental constructions 
of a certain kind;
its objects, 
relations, 
and proofs
exist only in so far as they have been constructed in
these acts.
Neither these acts nor the resulting constructions are of a linguistic nature, 
but when series of construction acts and their results are described in a language,
the descriptions may come to exhibit linguistic patterns. 
Intuitionistic logic is the mathematical study of these patterns, 
and in particular of those that characterise correct inferences.%
\footnote{As Heyting put it,
‘every logical theorem
\elide
is but a mathematical theorem of extreme generality’~\citep[p.6]{Heyting1956}.
Note that intuitionistic logic,
thus conceived,
is not the logic
‘underlying’
intuitionistic mathematics,
as is sometimes said;
quite the opposite.}

\item\label{L100}
Sundholm and I suggest the following characterisation
of intuitionistically correct inference~\citep[p.26]{Atten.Sundholm2017}:
‘A correct inference is one where the construction required by its conclusion can be
found from hypothetical actual constructions for its premisses. 
That is, 
we assume that
constructions for the premisses 
\textit{have been} 
effected. 
The hypotheses here are epistemic
ones, 
in that the premisses are 
\textit{known}. 
Thus, 
they differ from assumptions of the usual
natural deduction kind, 
which merely assume that propositions are true. 
For Brouwer’s
conception of truth, however, 
only these epistemic assumptions play a role, 
since for him
to assume that a proposition is true is to assume that one has a demonstration of it, 
that is,
that one knows that it is true.
\elide
We will write
\(A \rightarrow B\)
for
“A 
(hypothetical)
actual construction for 
\(A\)
can be continued into a construction for
\(B\)”.’

\item The valid logical principles then are 
rules under which all of mathematics 
(however it develops)
is closed.
\end{enumerate}
Sundholm has claimed that 
\cs{+}{3b}
‘is simply wrong’ 
if implication is understood
as a relation between a non-epistemic assumption and a consequent~\citep[p.19]{Sundholm2014}. 
I am inclined to agree with this, 
and refer to pages 18–20 of his
paper for further discussion of this point.%
\footnote{In conversation, Sundholm emphasised that that discussion does not
target the Creating Subject arguments in Brouwer’s own setting.} 
In the early discussions of 
\cs{}{}
and
\bks{}{},
I have found no one whose explanation of these principles involves 
implications with non-epistemic
assumptions;%
\footnote{In particular Myhill 
is explicit about his epistemic conception~[\citealp[p.295]{Myhill1967}~\citealp[p.326]{Myhill1968b}].}
it seems that it was not realised at the time that this meant that
the treatment of assumptions 
(and hence of implication) 
in Natural Deduction
renders the latter inappropriate for a faithful modelling of 
reasoning in Brouwerian mathematics.

\subsection{BKS is derivable from CS}\label{L080}

\begin{thm}[{{[\citealp[p.295-296]{Myhill1967}; \citealp[p.191-192]{Rootselaar1970}; \citealp[p.96]{Troelstra1969}]}}]
\cs{-}{}  implies \bks{-}{},
and
\cs{+}{}  implies \bks{+}{}.
\end{thm}

\begin{prf}
For both cases,
define
\begin{equation}
α(n) = 
\left\{
\begin{array}{ll}
0 & \text{if }  \neg\csop_{n} A\\
1 & \text{if \makebox[\widthof{\(\neg\)}]{}} \csop_{n} A
\end{array}
\right.
\end{equation}
By \cs{-}{1}
we have
\(\forall n(α(n)=0  \vee α(n)=1)\).

Assume
\(\forall n(α(n)=0)\).
Then,
by the definition of
\(α\),
\(\forall n \neg\csop_{n} A\), 
hence
\(\neg\exists n \csop_{n} A\),
and so,
by the contraposition of
\cs{-}{3a}
and
\(\neg\neg\neg A \leftrightarrow \neg A\),
\(\neg A\).
Conversely,
assume 
\(\neg A\).
Then,
by the contraposition of
\cs{-}{3b},
\(\neg\exists n \csop_{n} A\),
hence
\(\forall n \neg\csop_{n} A\),
and so,
by the definition of
\(α\),
\(\forall n(α(n)=0)\).

Assume
\(\exists n(α(n)=1)\).
Then,
by
the definition of
\(α\),
\(\exists n \csop_{n} A\),
and hence,
by
\cs{-}{3b},
\(A\).

To obtain moreover
\bks{+},
assume
\(A\).
Then,
by
\cs{+}{3a},
\(\exists n \csop_{n} A\),
and therefore,
by the definition of
\(α\),
\(\exists n(α(n)=1)\).
\end{prf}

\subsection{A justification of CS}\label{L105}

As will now be explained,
the conception of the Creating Subject that suits
Brouwer’s foundational thought
is that of an
ideal subject,
who can,
by itself, 
do whatever mathematics can in principle be done,
and whose activity
is structured as an
ω-sequence.
It remains a schematic subject in the sense
that it does not have a particular history,
but is the correlate of possible histories:
in each possible history,
the subject as thought of in that history
is an instantiation of the schematic Creating Subject.

The idea that the Creating Subject 
carries out its constructive activities in 
an ω-sequence of stages
is not made explicit by Brouwer,
but is implied 
whenever he states the view that the individual mathematical objects
that can be built up from the basic intuition of two-ity
are either finite or
or potentially infinite constructions.%
\footnote{Early examples of statements to that effect are found in
the dissertation~\citep[p.142–143]{Brouwer1907} (1907)
and in
‘\german{Die möglichen Mächtigkeiten}’~\citep{Brouwer1908A}
(1908);
for a late one,
see the quotation from
‘Guidelines of intuitionistic mathematics’
(1947)
on p.\ref{L126}
below.)}
At any given moment,
only finitely many constructions will have been carried out,
with an open horizon for further ones.
To see its mathematical activity as an ω-sequence of stages,
the Creating Subject first looks back at its earlier activity
(reflection),
and projects an initial segment of an
ω-sequence
(e.g., that of the natural numbers)
onto its various earlier acts.
The Creating Subject also sees that it can do this again and again in the future.
It is clear that actual human beings
carry out only finitely many mathematical acts
in their finite lives;
but in the study of Brouwerian constructibility the choice
to allow for potentially infinite time is
as reasonable
as it is in the study of computability.
Indeed,
ituitionism is a theory about an idealised mathematician in the same sense
as 
Turing’s theory of computable numbers is a theory about an idealised
(human or artificial)
computer,
and
Chomsky’s theory about grammar is a theory about an idealised speaker.
(For further discussion of these parallels,
see~\citealp[esp.~p.234–235]{Gil1983};~\citealp[p.148-151]{George1988};~\citealp[ch.6]{Atten2004};~\citealp{Atten-forthcomingB}.)

The idealisation involved in such theoretical accounts
also includes the idea
that we never make mistakes.
As Troelstra aptly put it at the beginning of the chapter on Creating Subject arguments in
\textit{Principles of Intuitionism},
\begin{quote}
The central idea is that of an idealized mathematician 
(consistent with the subjectivistic viewpoint of intuitionism, 
we may think of ourselves; 
or even better, 
to obtain the required idealization, 
we may think of ourselves as we should like 
to be).~\citep[p.95]{Troelstra1969}
\end{quote}
Van Dantzig~\citep{Dantzig1949}
turned the fact that such idealisations are required
to make Brouwer’s Creating Subject arguments work into an
objection to them;
but this is misguided in view of the nature of the theoretical model
of our mathematical activity that intuitionism presents.%
\footnote{Markov 
in effect shared Van Dantzig’s objection that
the assumption that the Creating Subject
never proves a proposition
\(A\)
does not imply that
\(A\)
is false;
perhaps the Creating Subject
just lost interest in the problem!
(Markov makes this remark in his Russian translation~\citep[p.192]{Heyting1965}
of Heyting’s
\textit{Intuitionism}~\citep{Heyting1956};
the remark is translated into French by Margenstern~\citep[p.290]{Margenstern1995}.)
However,
allowing for such a possibility is already to restrict
the Creating Subject~\citep[section 3]{Atten-forthcomingB}.}

Another idea common to these analyses of 
constructibility,
computability,
and grammaticality
is that whatever falls under these notions 
can in principle be fully mastered by a single subject.
The Creating Subject is essentially singular.%
\footnote{Compare on this point~\citep[p.109–110n15]{Cellucci1974}.}
This may seem to be contradicted by
Brouwer’s remark that one and the same untestable proposition can give rise
to different Creating Subject sequences~\citep[p.204n2]{Brouwer1954F},
or,
as he put it more vividly in the 1934 Geneva lectures:
\begin{quote}
If I
would give the definition of 
[the Creating Subject sequence]
\(s\) 
to one hundred different persons, 
who are all going
to work in a different room, 
it is possible that one of these one hundred persons at
one time will choose an interval not covered by an interval chosen by one of the
others.~[\citealp[lecture 2]{Brouwer1934}, trl.~\citealp[p.45]{Niekus2010}]
\end{quote}
But
I read this not as an acknowledgement that the notion of Creating Subject would be fundamentally
plural,
in such a way that mathematical truth would not depend on a single one of them,
but rather as an acknowledgement of the schematic character of the notion,
and that different real or imagined instantiations may proceed differently 
from one another.%
\footnote{This may also be behind Brouwer’s objection~\citep[p.11]{Brouwer1930B}
to Fraenkel’s remark that
all who agree that a given mathematical question is meaningful will give the same
answer to it,
but Brouwer would perhaps say that already out of a distrust of language.}
This again corresponds to the theoretical role of
the Universal Turing Machine.
The latter is defined in general terms,
and then in thought experiments actual human beings theorise about
different runs of it with different programs and different inputs,
so that each of these runs could be said to present a possible
history of the Machine’s computations,
without any of them being the actual history,
because,
being schematic,
the Machine has no history.

With this conception of the Creating Subject
in place,
the schemata for
\cs{+}{}
are justified as follows,
in keeping with
the conception of logic
described in subsection~\ref{L103}.
 
\cs{+}{1} is correct because,
if 
\(n\) 
is in the past,
the Creating Subject inspects its perfect memory~\citep[section 3]{Brouwer1933A2}; 
and if 
\(n\) 
is in the future,
the Creating Subject  postpones its decision for the finite number of stages required
and then sees again.

\cs{+}{2}
is justified by the Creating Subject’s perfect memory.

\cs{+}{3b} 
can be read as an assertion of the correctness of the Creating Subject’s thought.%
\footnote{In a draft version of the
‘\dutch{stellingen}’
for his thesis defense,
Brouwer pointed out that
in mathematics there should be no hypothesis
‘I reason correctly’:
‘First of all,
reasoning is an
\textit{act},
in which the self is not objectivated;
secondly,
these words have meaning,
and,
a fortiori,
meaning as a foundation for something else,
only on the basis of 
already existing mathematical systems,
and therefore of already existing logic;
and,
thirdly,
in particular the word 
“correctly”
means nothing but 
“mathematically correct”,
and therefore presupposes mathematics and logic.’~\citep[p.147, trl.~MvA]{Brouwer2001}.
Creating Subject arguments clearly are an example of what Brouwer
(in 1907)
called
‘second-order mathematics’~\citep[p.119n.]{Brouwer1907},
which allows one 
to describe,
and 
to reason about,
mathematical acts mathematically.
In such an objectivation,
there certainly is a place for descriptions of aspects of that activity such as
\cs{+}{3b};
Brouwer’s point,
applied here,
is just that accepting
\cs{+}{3b}
is in no way a condition for engaging in
mathematical activity as such.
Another example
(one that Brouwer remarks on in his dissertation)
is that of mathematical induction as a construction act
(as opposed to an axiom (schema))~\citep{Atten-forthcomingC}.}
Martino~\citep[section 5.5 ]{Martino1985}
argues that for
\cs{+}{3b}
it does not matter whether the proof that provides evidence for
\(A\)
comes from within the Creating Subject or is given to it by someone else;
and that there is therefore nothing solipsistic about it.
However,
it seems to me that the Creating Subject cannot simply accept a proof from outside,
and must go through the steps to reconstruct the evidence.
In the end,
then,
its own ability to make
\(A\)
evident is what really matters.
The Creating Subject must take a proof that is communicated to it
not as itself evidence, 
but as a set of instructions to obtain that evidence~\citep[p.169]{Brouwer1907}.

\cs{+}{3a} 
does not mean,
as it may seem at first sight,
‘If 
\(A\)
is true,
whether the Creating Subject knows it or not,
a time will come when the Creating Subject proves it’.
On the Brouwerian conception of logic,
it rather gives expression to
the fact
that
whenever 
the Creating Subject has effected
a construction for 
\(A\),
and thereby experienced the truth of
\(A\),
it has done so at a particular moment in time,
and that moment can be made explicit~[\citealp[p.295]{Myhill1967};~\citealp[section 5.5]{Martino1985};~\citealp[p.242]{Dummett2000}].

In a discussion of
\bks{+}{}
that in effect concerned
\cs{+}{3a},
Kleene and Moschovakis suggested to Myhill that 
‘this presupposes a linear ordering of all possible proofs,
so that the apparent reference to time is really a reference
to Gödel-numbers of 
(informal) 
proofs’~\citep[p.295]{Myhill1967}.%
\footnote{Webb~\citep[p.211]{Webb1980}
writes of
\cs{}{}
that
‘ it implies that all possible proofs can be arranged in
an 
\(ω\)-sequence, 
which acutely conflicts with the impredicative nature of
intuitionistic implication.’
I argue in~\citep{Atten-E} that intuitionistic implication
is predicative.}
I agree with Myhill that this idea does not accord
with Brouwer’s speaking of
‘experiencing the truth’,
and would explain that by the fact that for Brouwer
the number is extracted from that experience
but not from its content.
Moreover,
a Brouwerian understanding of
the species of 
‘all possible proofs’
would itself depend on time,
as it must be understood as an essentially growing species,
for which no single construction method can be given:
\begin{quote}
As further examples of denumerably unfinished sets
we mention:
the totality of definable points on the continuum,
and a fortiori the totality of all possible mathematical systems.~[\citealp[p.148–149]{Brouwer1907};
trl.~\citealp[p.82]{Brouwer1975}]
\end{quote}
and,
in a notebook around 1907,
\begin{quote}
The totality of mathematical theorems is, 
among other things, also a set that is denumerable, but never finished.~\citep[Notebook VIII, p.44, trl.~MvA]{Brouwer19041907}%
\footnote{‘\dutch{Het aantal wiskundige stellingen is o.a.~ook een \german{Menge}, die aftelbaar is, maar nooit af.}’}
\end{quote}
The idea also occurs in Kreisel:
\begin{quote}
First of all,
very little of the 
‘thinking subject’ 
is used in the derivation 
[in a reconstructed Creating Subject argument], 
Instead of writing 
\(Σ \vdash_n A\), 
I could write
\(Σ_n \vdash A\)
and read it as: 
the 
\(n\)-th proof establishes 
\(A\). 
In other words,
the essential point would not be
the individual subject, 
but the idea of proofs arranged in an 
\(ω\)-order.
\elide
Also,
the sequence 
\(Σ_n\)
is not itself considered to be given by a rule.~\citep[p.179]{Kreisel1967b}
\end{quote}
The question is whether 
the existence of the sequence 
\(Σ_n\) 
could be justified
constructively
if the
\(ω\)-order
it depends on
is not the one that is induced by
the fact that
the Creating Subject’s
activities
unfold over time,
so that,
genetically,
the individual subject is essential after all.

A further analysis of
the reference to time in
\cs{+}{3a}
than has been advanced so far,
and which
will be important also for the discussion of Kripke’s own
objection to \bks{+}{} 
in subsection~\ref{L031}, 
depends on two distinctions.

The first is that between tokens and types.
The other 
is
a distinction between three 
meanings of the term ‘construction’,
drawn attention to by Sundholm~\citep[p.164]{Sundholm1983}:
\begin{enumerate}
\item process of construction (as it unfolds in time),
\item object obtained as the result of such a process,
\item  construction-process as object
(the objectification of a process of construction).%
\footnote{\label{proof-trace}
There is a
fourth sense~[\citealp{Sundholm1993};~\citealp[p.68]{Sundholm.Atten2008}]: 
‘by abstracting of the objectified act with respect
to subject and time, 
[one obtains a construction in the sense of] 
a prescription or blueprint for construction acts’.
I will not be concerned with 
the abstraction from a (particular) subject,
because my discussion concerns only the one Creating Subject;
whereas abstraction from (particular) time is dealt with
by applying the type-token distinction.}
\end{enumerate}

When Brouwer in his dissertation writes that
‘strictly speaking the construction of intuitive mathematics in itself is an 
\textit{act} 
and not a 
\textit{science}~\citealp[p.99n]{Brouwer1907}; 
trl.~\citealp[p.61n1]{Brouwer1975}, modified]%
\footnote{I here translate ‘\dutch{daad}’ by ‘act’ instead of ‘action’.}
he is thinking of constructions in the first sense;
and it is clear that constructions in the other two senses
presuppose for their existence a construction process in the first sense.
Constructions in the first sense are ontologically prior to the others.
The objectification of a process happens in an act of reflection;
this possibility to reflect on our acts will turn out to be crucial
to a Brouwerian argument for \bks{}{}.

This can be connected to the type-token distinction as follows.
At the most concrete level,
construction processes occuring at different times are
for that reason different processes.
But we may come to see 
that processes that are different in this sense
have various things in common,
and we may therefore see them
as instantiations or tokens of the same type of construction process.
The same can be done for constructions in the other two senses,
constructed objects
and
the objectified processes.
For example,
this allows us to observe that an act in which we construct the number 
\(2\)
and
an act in which we construct the number 
\(3\)
have in common that the objects constructed in them are of the same type,
that of natural number.
In the extreme case,
we may even come to
identify  processes with one another,
and identify the objects constructed in them.
This is the 
sense in which we can say,
for example,
that when constructing the number 
\(2\) 
time and again,
each time we carry out the same construction process
in which we construct the same object.

The notion of proof is related to that of a construction
in a straightforward way.
A clear statement of this relation was made by Heyting:
\begin{quote}\label{L079}
If mathematics consists of mental constructions,
then every mathematical theorem is
the expression of a result of a successful construction.
The proof of the theorem
consists in this construction itself,
and the steps of the proof are the same as the
steps of the mathematical construction.~\citep[p.107]{Heyting1958A}
\end{quote}
So for proofs
the same threefold distinction can be made
as was introduced for constructions,
and we may consider each case as a type or as a token,
depending on our purpose.%
\footnote{This is closely related to Martin-Löf’s notion of a ‘proof-trace’~\citep{Sundholm1993}.}

The question what counts as a proof of a proposition
\(A\) 
must be a question about the proof type;
but whenever a proof of 
\(A\)  
is given to us,
what is given to us is primarily a proof token,
as intuitionistically types only exist as abstractions from tokens. 
The parameter 
\(t\) 
is part of the proof token,
as
the specific time at which a concrete process occurs is
constitutive of its identity;
but it is not part
of the proof type,
as proofs of the same type may be constructed at different times.

Logical principles are formulated at the level of proof types,
so as to allow these principles to express the relevant kind of generality:
such a principle is of the form
‘Whenever I have a proof of the premises,
I can obtain from it a proof of the conclusion’.
But
observe that 
whenever I have a proof of the premises,
this is given to me first of all as a proof token,
and only abstractively as a type.
This is because, 
intuitionistically,
types do not exist independently from their tokens.

This difference creates room for the following argument.
When claiming that an implication 
\(A \rightarrow B\)
holds,
what is claimed
is that a certain relation holds between the type 
‘proof of 
\(A\)’
and the type
‘proof of 
\(B\)’,
namely,
that there is a construction method to
convert 
\textit{any} 
token proof of 
\(A\)
into a token proof of 
\(B\).
So whenever the Creating Subject wishes to apply that method
to a token proof of the antecedent
(modus ponens),
the parameter 
\(t\) 
is always available among the data.
It is 
\textit{not} 
claimed 
that 
\(t\)  
is among the specifications of a proof type;
nor is that required to make good on the claim
\(A \rightarrow B\).
In other words,
the idea is 
\emph{not}
that a value for 
\(t\)  
can somehow be extracted from the propositional content of 
\(A\).

Perhaps it is objected that,
even if we acknowledge this distinction between construction types
and tokens,
and that givenness of a token proof of the antecedent of an implication includes givenness of 
\(t\),
this 
does not suffice to consider 
\(t\)  
part of the
\emph{mathematical}
data we have on which to base a construction for the consequent.

First note that this objection,
if correct,
would apply to
\bks{-}{}
as well.
That schema too requires
that the time parameter 
\(t\) 
be a
\textit{mathematical} 
datum.
Otherwise the sequence that witnesses the existential quantifier in 
\bks{-}{},
a witness the constructibility of which Kripke in effect justifies in terms of 
\(t\),
would not exist as an object of pure  mathematics.
(This again illustrates the theme from Kreisel and Myhill that
\bks{}{}
is an extensional principle but is justified by consideration of intensional aspects.) 
%

A reply to the objection itself must begin with the observation that
for Brouwer the Creating Subject is an idealised intuitionistic mathematician,
and this includes 
a property of consciousness that Husserl calls
‘inner time awareness’.
Brouwer
does not mention Husserl
in his dissertation,
nor elsewhere,
but he does
confirm there that inner time awareness is what he has in mind:
\begin{quote}
Of course we mean here 
\textit{intuitive time}
which must be clearly distinguished from
\textit{scientific time}.
Very much a posteriori,
only by means of experience 
it becomes clear that the latter
can 
suitably be introduced
as a one-dimensional coordinate equipped with a one-parameter group
for cataloguing phenomena.~[\citealp[p.99n]{Brouwer1907};
trl.~\citealp[p.61n2]{Brouwer1975}].
\end{quote}

To see why 
the time parameter in this sense
is 
\textit{mathematically} 
relevant for the Creating Subject,
recall Brouwer’s statement in his dissertation that mathematics is first of all an 
\textit{act}.
In a much later passage,
from 1947,
he said the same thing more specifically:
\begin{quote}
Intuitionistic mathematics is a mental construction,
essentially independent of language.
It comes into being
by self-unfolding of the basic intuition of mathematics,
which consists in the abstraction of two-ity.~[\citealp[p.339]{Brouwer1947},
trl.~\citealp[p.477]{Brouwer1975}]
\end{quote}
By 
‘self-unfolding’
is meant 
that in our mathematical acts we first construct certain basic objects,
the nature of which Brouwer specifies but is not relevant now,
and then to those 
apply the same mathematical acts
to construct further objects.
This activity thus has an iterative structure,
and
induces a linear order on the constructions
that the Creating Subject has effected.

Extracting the time parameter 
\(t\)
associated to a token construction that makes a proposition
\(A\)
true
is just making the position of the token in that induced ordering explicit.
To do this,
the Creating Subject needs to reflect on its mathematical activity so far and
apply mathematics to it.
This is an instance of what Brouwer,
in his dissertation,
describes as
‘viewing mathematical activity mathematically’,
a form of self-reflexivity 
which he calls 
‘second-order mathematics’
and which is itself of a  mathematical nature~[\citealp[p.98n]{Brouwer1907};
trl.~\citealp[p.61n1]{Brouwer1975}].
As the linear ordering of the Creating Subject’s constructions
is not only the order in which it
becomes aware of the objects,
but indeed the order in which they come into being,
this order is a mathematical fact.
It is this consideration that
justifies the strong Brouwer-Kripke Schema
in a Brouwerian setting.
Indeed,
more generally,
the schemata
\cs{}{}
are pieces of
general self-knowledge
of the Creating Subject,
obtained by applying mathematics
to mathematical activity as given in reflection.%
\footnote{Metschl~\citep{Metschl2000}
analyses the role of epistemic obligations
in the formation of the Creating Subject’s knowledge
(in general; not with an eye on Creating Subject arguments in particular).}

There are different ways of  relating 
the
\(ω\)-ordering of the inner time parameter 
\(t\)
to the
\(ω\)-ordering
of the elements of a sequence required for a witness
of the existential quantifier in \bks{}{}.
The
first element of the witnessing sequence
may be correlated to
the beginning of all of the Creating Subject’s
activity,
or to the moment at which
the Creating Subject has
indeed begun working towards a decision
of
\(A\).
However,
either of these sequences can be mapped to the other in an order-preserving way.
Similarly,
one master sequence tracking
all of the Creating Subject’s other mathematical activity
would be sufficient,
given that from it any
sequence
pertaining only to the activity related to a given
proposition
\(A\)
may be extracted as needed.

It was remarked above that
the number 
\(2\) 
could be constructed at many different times.
Would it then in principle be possible, 
in intuitionistic mathematics, 
to use different numbers 
\(2\)  
and the times that they occur?%
\footnote{I thank Saul Kripke for raising this question, in conversation.}
Indeed, 
by the same reasoning as 
applied to proofs in the present account of
\bks{+}{},
one could use the time parameter in the Creating Subject’s activity to
assign a different number
to each token construction of the number
\(2\).
But,
as emphasised in the quotation from Heyting
(p.\pageref{L079} above),
any successful construction process may be considered to be a proof 
(whether subsequently expressed in a theorem or not);
to construct a token of the number
\(2\)
then is,
in that sense,
to construct
a token of the proof (type)
that the number 
\(2\),
considered as a construction type,
exists.
It therefore seems to me
that even though
the notions of proof and object are different,
assigning numbers also to objects
would not 
yield mathematical possibilities 
beyond those given by 
\bks{+}{}.

\section{BKS in Brouwer}\label{L075}
\bks{}{}
is found
in Brouwer’s published work:
on several occasions,
he reasons from the antecedent of an instance
of 
\bks{-}{}
to its consequent,
and he once
establishes an equivalence 
for untested
\(A\)
from which
\bks{+}{}
for such propositions
follows
immediately
and hence,
it can be argued,
for other propositions as well.
These occurrences are discussed below.
Brouwer does not make any fanfare about these particular inferences,
and in particular he does not pause to formulate the general principles. 
On any interpretation of Brouwer
according to which he in effect accepted 
\cs{}{}
on the grounds given in subsection~\ref{L105}, 
this silence is not surprising, 
because Brouwer will then have carried out his reasoning 
in the corresponding terms  directly.
After all,
the derivations of
weak and strong \bks{}{}
from
the respective versions of
\cs{}{}
are very short, 
and Brouwer never quite bothered to make
the latter explicit either.

When\label{L032} Kripke isolated
\bks{-}{},
he was not aware of occurrences 
of 
\bks{-}{}
and
\bks{+}{}
in Brouwer.%
\footnote{Personal communication.}
Clearly,
the schema was picked up on in the literature
because of Kripke’s rediscovery,
not because of Brouwer’s earlier and implicit use.
So the name
‘Brouwer-Kripke Schema’%
\label{L139}%
,
which seems to have been introduced in print by
De Swart in 1977~\citep[p.578]{Swart1977},%
\footnote{He had used it in 1976 in his Nijmegen dissertation~\citep[p.34-35]{Swart1976b}.
Gielen, Veldman and De Swart in 1981
speak of
‘the Brouwer-Kripke Axiom’~\citep[p.122,126]{Gielen.Swart.Veldman1981}.
In a manuscript 
‘The trustworthiness of intuitionistic principles’ (1983),
Gielen~\citep{Gielen1983}
uses ‘Brouwer’s Scheme’
for a version of 
\bks{}{}
with a restriction on the content
\(A\)
which is discussed in 
subsection~\ref{L090} below.
He may have come to do this because
the Nijmegen intuitionists
considered
that restriction 
to be a necessary condition 
for the schema’s validity,
and to be motivated by their interpretation of Brouwer,
whereas it
is absent from Kripke’s
formulation.
(De Swart’s dissertation~\citep[p.35]{Swart1976b}
presents a result on 
\bks{}{}
from Gielen’s ‘\dutch{Verzamelingen}’.
I have not seen that unpublished manuscript
from 1976 at the latest,
and do not know what term he used there.)}
is appropriate.

For the special case of lawless sequences,
the connection between the strong counterexample to
\(\forall x\in \numreal(x \neq 0 \rightarrow x  \apart 0)\)
and MP
had been made by Kripke in 1965,
when he proved that MP implies Brouwer’s theorem~\citep[p.103–104]{Kripke1965}.
In that paper 
\bks{-}{}is not yet stated,
but comes near the surface;
especially when on p.104 it is suggested to replace
lawless sequences by sequences based on the solving of problems.
For 
lawless sequences, 
the arguments are simpler and
\bks{-}{} is not needed.
Kripke’s reference for the Creating Subject arguments there is
Heyting’s
\textit{Intuitionism}~\citep{Heyting1956},
which
includes Brouwer’s argument from 1949,
but not that of 1954 
(subsections~\ref{L066} and~\ref{L069} above);
and there is no paper by Brouwer among the references.
Kripke
remarks that
‘I think it probable that such treatments in FC [= Kreisel’s formal theory of lawless sequences~\citep{Kreisel1958b}]
will extend
to all the counterexamples to classical theorems which Brouwer gives by this method;
but I have not made a survey of the literature’~\citep[p.104]{Kripke1965}.

\subsection{\texorpdfstring{\bks{-}{}}{BKS-}}\label{L095}

As observed by
Gielen,
De Swart,
and Veldman~\citep[p.128]{Gielen.Swart.Veldman1981},
in the proof from 1949 discussed
as Proof~\ref{L074} above,
Brouwer 
reasons along the lines of \bks{-}{}.
This occurs where,
from the hypothesis that
for a direct checking number
\(D(γ,p_f) \gro 0\)
holds,
a proof is constructed 
for a previously defined real number
\(f\)
that 
\(f \in\numrat \vee f \not\in \numrat\).
This reasoning 
clearly shows the pattern
that makes it amenable 
to reconstruction with
\bks{-}{}:
from the 
(hypothetical)
existence of a value distinct from 
\(0\) 
in a certain infinite sequence,
the truth of a priorly given proposition can be concluded to.%
\footnote{Note that Myhill’s criticism and repair of the proof
(see subsection~\ref{L067})
concerns only the part that
comes after this.}
In this argument,
no
witnessing sequence
for 
\bks{-}{}
is actually constructed;
rather,
the hypothesis
(towards a contradiction)
that
\(>\)
implies
\(\gro\)
amounts to a hypothesis that such a sequence
has been constructed.
(To accept the implication of
\(\gro\)
by
\(>\)
would 
in fact
come down to accepting
Markov’s Principle for real numbers given by arbitrary kinds of 
converging choice sequences
(see page~\pageref{L006} above).)

Brouwer  also argues in this manner
in
a weak counterexample in
‘Points and spaces’~\citep{Brouwer1954A}. 
Some definitions are needed first.

\begin{dfn}[{{\citep[p.8-9]{Brouwer1954A}, simplified}}]
A 
\textit{spread direction}
is a tree 
over the natural numbers
such that each node
\(p\)
allows either all natural numbers
as its immediate descendants,
or all natural numbers 
\(\leq m_p\)
for some
\(m_p\).

A
\textit{spread}
is a species of infinite paths through
a spread direction.

A subspecies of a spread direction
is
\textit{thin}
if no node in it is a descendant of any of the other nodes.

A subspecies of a spread direction
such that no infinite path through the spread direction
can fail to intersect it is called a
\textit{crude block}.

A decidable,
thin subspecies of a spread direction
that is a block is called a
\textit{proper block}
or simply a
\textit{block}.%
\footnote{The terminology in the 
\textit{Cambridge Lectures},
held from 1946 to 1951,
was different~\citep[p.21–22]{Brouwer1981A}.
The crude block of
‘Points and spaces’
was there called a
\textit{barrage};
what was there called
a
\textit{crude block}
is the positive counterpart of a barrage,
i.e.~a  subspecies of a spread direction
such that every infinite path through the spread direction
intersects it.
The notions of
\textit{proper block}
or simply
\textit{block},
defined as
that of
thin and decidable
crude block,
likewise were positive.}
\end{dfn}

\begin{dfn}[{{[e.g.~\citealp[p.10]{Brouwer1954A}; \citealp[p.44]{Brouwer1981A}]}}]\label{L097}
Let
\(R\)
be
a finite or infinite sequence
of species
\(N_{ν}\),
such that 
the
\(N_{ν}\)
are all disjoint and each completely ordered
by a respective
\(<_{N_{ν}}\).
Then the
\textit{ordinal sum}
of the 
\(N_{ν}\)
is their union species
\(M\)
equipped with a complete order
\(<_M\)
such that,
for 
\(x \in N_{ν1}\)
and
\(y \in N_{ν2}\),
\begin{equation}
x <_M y \equiv N_{ν1} <_R N_{ν2} \vee (N_{ν1}=N_{ν2}=N_k \wedge x <_{N_k} y)
\end{equation}
The  
\textit{well-ordered species} 
are defined inductively:
\begin{enumerate}
\item A species containing exactly
\(1\)
element is a
\textit{basic species}.
Basic species are well-ordered species.

\item The ordinal sum of an infinite sequence of previously acquired
disjoint well-ordered species is again a well-ordered species.

\item The ordinal sum of a non-empty finite sequence of disjoint previously acquired
well-ordered species is again a well-ordered species.
\end{enumerate}
\end{dfn}

\begin{wce}[{{\citep[p.12]{Brouwer1954A}}}]
There is no hope of showing that
\begin{equation}\label{L}
K \text{ is a block } \rightarrow K \text{ is a well-ordered block} 
\end{equation}
for arbitrary spread directions 
\(K\).
\end{wce}

\begin{pargmt}
Let
\(K\)
be a spread direction,
and
\(X_n\)
the species of nodes of
\(K\)
at depth
\(n\).
An
\(n\)-union
is a union of species
\(X_n\);
it is
constructed step by step
from a choice sequence
\(α\)
in which choice
\(α(n)\)
determines whether
the
\(n\)-th species
\(X_n\)
will be included in the
\(n\)-union
or not.

Let 
\(A\) 
be a proposition 
that is at present not
testable.
The Creating Subject constructs 
the 
\(n\)-union
\(n_A\)
by making its sequence of choices
\(α\)
as follows:
\begin{itemize}
\item As long as,
by the choice of 
\(α(n)\),
the Creating Subject
has obtained evidence neither of 
\(A\)
nor of 
\(\neg  A\),
\(α(n)\) is chosen to be negative.

\item If between the choice of 
\(α(r-1)\) and
\(α(r)\),
the Creating Subject
has obtained evidence
either of 
\(A\)
or of
\(\neg A\),
\(α(r)\) 
is chosen to be positive.

\item For all
\(n > r\),
\(α(n)\)
is
chosen to be
negative again.
\end{itemize}
Then 
the 
\(n\)-union
\(n_A\)
is a block of 
\(K\)
because
\(n_A\)
cannot be empty:
its being empty would imply
\(\neg(A \vee \neg A)\),
which is contradictory.
Note that,
by definition,
the elements of
\(n_A\)
cannot reside at different depths.
On the other hand,
for it to be a well-ordered block,
it would have to be known 
what the 
nodes in
\(n_A\)
are,
and hence what their depth is,
and this can only be known if
\(A \vee \neg A\)
is known.
So as long as 
\(A\)
has not been decided,
\(n_A\)
is a block, 
but not a well-ordered block.
\end{pargmt}

The argument
depends on the fact
that the sequence
\(α\)
satisfies
the properties
\begin{equation}
\begin{gathered}
\forall n(α(n)=0  \vee α(n)=1) \\
\wedge \\
\forall n(α(n)=0) \leftrightarrow \neg (A \vee \neg A)\\
\wedge\\
\exists n(α(n)=1) \rightarrow A \vee \neg A
\end{gathered}
\end{equation}
and these are precisely the properties guaranteed by
\bks{-}{}.
MP,
if it were valid for sequences like
\(α\),
would allow one to accept 
\(n_A\)
as a well-ordered block.%
\footnote{The referee pointed out that there is a much simpler argument:
\begin{quote}
An example as sought for is already given by
\(B \coloneqq \{s \mid α(\length(s)) \neq 0\}\),
where 
\(α\)
satisfies:
\(\neg\neg\exists n[α(n) \neq 0]
\wedge
\forall m\forall n
[
(α(m) \neq 0
\wedge
α(n) \neq 0)
\rightarrow
m=n
]\),
and one is unable to prove:
\(\exists n[α(n) \neq 0]\).
\end{quote}
The referee adds that the notion
of a block is unlikely to be fruitful in constructive mathematicsl.
I mention these points for their intrinsic interest;
for the discussion of Brouwer’s willingness to reason along the lines
of
\bks{-}{},
both are moot.}

A third place
where Brouwer  
proceeds thus
is in
‘Intuitionistic differentiability’,
also from 1954~\citep{Brouwer1954E}.
Brouwer had studied intuitionistic differentiation only in 1923~\citep{Brouwer1923A}; 
this counterexample may have been motivated by
Van Rootselaar’s then recent work~\citep{Rootselaar1952,Rootselaar1954}.%
\footnote{The latter reference is to 
Van Rootselaar’s
dissertation,
supervised by Heyting
and defended in 1954.
In his
\textit{Intuitionism}
Heyting
does not treat the topic, 
and refers~\citep[p.96]{Heyting1956}  
to
Van Rootselaar~\citep[ch.5]{Rootselaar1954}.}

\begin{dfn}[{{\citep{Brouwer1954E}}}]
An interval 
\([r,s] \subset \numreal\)
is 
\textit{substantial}
if
\(r \apart s\).
Let 
\(f \colon \numreal\rightarrow\numreal\)
be a total function,
\(r\)
an arbitrary
\(x\in\numreal\),
\(s=i_{s1},i_{s2},\dots\)
an infinite sequence of substantial intervals
containing
\(r\)
that converges positively to
\(r\).
Associated to each
\(i_{sv} = [i_{sv1},i_{sv2}]\)
is the 
\textit{difference quotient}
\[
d_{sv} = \frac{f(i_{sv2})-f(i_{sv1})}{i_{sv2}-i_{sv1}}
\]
\(f\) is
\textit{strongly differentiable}
in
\(x=r\),
and has
the
\textit{strong differential quotient 
\(c\in\numreal\) in 
\(x=r\)},
if
\begin{equation}\label{L061}
\forall n\exists m\forall s\forall v(\absval{i_{sv2}-i_{sv1}} < 2^{-m} \rightarrow \absval{d_{sv}-c} < 2^{-n})
\end{equation}
\(f\) is
\textit{weakly differentiable}
in
\(x=r\),
and has
the
\textit{weak differential quotient 
\(c\in\numreal\) in 
\(x=r\)},
if
\begin{equation}
\neg\exists s\exists m\forall v(\absval{d_{sv}-c} > 2^{-m})
\end{equation}
\end{dfn}

\begin{wce}[\citep{Brouwer1954E}]\label{L129}
There is no hope of showing that
\begin{equation}
\text{\(f\) is weakly differentiable in 
\(x=r\)} \rightarrow \text{\(f\) is strongly differentiable 
in 
\(x=r\)}
\end{equation}
\end{wce}

\begin{pargmt}
Define the family of functions
\(ω_{\nu \in \numnat^{+}} \colon \numreal \rightarrow [0,\frac{1}{4}]\)
(Figure~\ref{L063})
by
\begin{equation}
ω_\nu(x) =
\begin{dcases*}
0 & \text{if  
\(x \leq -2^{-\nu + 1}\)}\\
\sqrt{-3 \cdot 2^{-\nu}x - x^2 - 2^{-2\nu + 1}}  & \text{if  
\(-2^{-\nu + 1} \leq x \leq -2^{-\nu}\)}\\
0 & \text{if  
\(-2^{-\nu} \leq x \leq 2^{-\nu}\)}\\
\sqrt{3 \cdot 2^{-\nu}x - x^2 - 2^{-2\nu + 1}}  & \text{if  
\(2^{-\nu} \leq x \leq 2^{-\nu+1}\)}\\
0  & \text{if  
\(x \geq 2^{-\nu+1}\)}
\end{dcases*}
\end{equation}

\pgfplotsset{width=1.1\textwidth,height=5cm}

\pgfplotsset{every axis/.append style={%
line width=0.3pt,
tick style={line width=0.3pt,black}}}

\begin{figure}
\centering
\begin{tikzpicture}
\begin{axis}[
axis x line=bottom,
axis y line=center,
black,
mark=none,
samples=100, 
tips=never,
xlabel=
\(\),ylabel=
\(\),
xmin=-1.05,
xmax=1.05,
xtick={-1,-0.667,-0.5,-0.333,-0.25,-0.1667,-0.125,-0.083,-0.0625,0,0.0625,0.083,0.125,0.1667,0.25,0.333,0.5,0.667,1},
xticklabels={%
	\(-1\),
	\(-x_1\), 
	\(-\frac{1}{2}\),
	\(-x_2\), 
	\(\), 
	\(-x_3\), 
	\(\), 
	\(\),
	\(\),
	\(0\),
	\(\frac{1}{16}\),
	\(\), 
	\(\frac{1}{8}\), 
	\(\), 
	\(\frac{1}{4}\),
	\(x_2\),
	\(\frac{1}{2}\),
	\(x_1\), 
	\(1\)},
ytick={0.03125,0.0625,0.125,0.25}, 
yticklabels={%
	\(\),
	\(\frac{1}{16}\),
	\(\frac{1}{8}\),
	\(\frac{1}{4}\)}
]
\addplot [
domain=0.5:1,
line width=1pt,
] {sqrt(3*2^(-1)*x-x^2-2^(-1))}; 

\addplot [
domain=-1:-0.5,
line width=1pt,
] {sqrt(-3*2^(-1)*x-x^2-2^(-1))}; 

\addplot [
domain=0.25:0.5,
line width=1pt,
] {sqrt(3*2^(-2)*x-x^2-2^(-3))}; 

\addplot [
domain=-0.5:-0.25,
line width=1pt,
] {sqrt(-3*2^(-2)*x-x^2-2^(-3))}; 

\addplot [
domain=0.125:0.25,
line width=1pt,
] {sqrt(3*2^(-3)*x-x^2-2^(-5))}; 

\addplot [
domain=-0.25:-0.125,
line width=1pt,
] {sqrt(-3*2^(-3)*x-x^2-2^(-5))}; 

\addplot [
domain=0.0625:0.125,
line width=1pt,
] {sqrt(3*2^(-4)*x-x^2-2^(-7))}; 

\addplot [
domain=-0.125:-0.0625,
line width=1pt,
] {sqrt(-3*2^(-4)*x-x^2-2^(-7))}; 

 \addplot [
 domain=-1.3:-1,samples=10,
line width=1pt,
 ] {0*x}; %

 \addplot [
 domain=1:1.3,samples=10,
line width=1pt,
 ] {0}; %

 \addplot [
 domain=-0.0625:0.0625,samples=10,
line width=1pt,
 ] {0*x}; %

\addplot [
domain=0:1,samples=10,
] {0.25*sqrt(2)*x}; %

\addplot [
domain=-1:0,samples=10,
] {0.25*sqrt(2)*-x}; %

\node [align = center] at (0.75,0.05) {\({ν=1}\)};
\node [align = center] at (-0.75,0.05) {\({ν=1}\)};
\node [align = center] at (0.375,0.05) {\({ν=2}\)};
\node [align = center] at (-0.375,0.05) {\({ν=2}\)};

\node at (0.6,0.3) {\({{\Delta y}/{\Delta x}=\frac{1}{4}\surd{2}}\)};
\end{axis}
\end{tikzpicture}%
\caption{\(\sum\limits_{ν} ω_{ν} \) for the first four values of 
\(ν\).\label{L063}}
\end{figure}

Let 
\(A\) 
be a proposition 
that is at present not
decidable.
The Creating Subject constructs a  
choice
sequence 
\(ζ\)
of 
total real functions 
\(ζ(n) \colon \numreal \rightarrow [0,\frac{1}{4}]\):
\begin{itemize}
\item As long as,
by the choice of 
\(ζ(n)\),
the Creating Subject
has obtained evidence neither of 
\(A\)
nor of 
\(\neg  A\),
\(ζ(n)\) is chosen to be 
\(ω_{n}\).

\item If between the choice of 
\(ζ(m-1)\) and
\(ζ(m)\),
the Creating Subject
has obtained evidence
of 
\(A \vee \neg A\),
\(ζ(n)\) 
for all
\(n \geq m\)
is
chosen to be
the constant function
\(\lambda x.0\).
\end{itemize}

Define
the function
\(Z \colon \numreal \rightarrow [0,\frac{1}{4}]\)
by
\begin{equation}
Z(x) = \sum_{ν=1}^{\infty} ζ(ν)(x)
\end{equation}
From the definition of the 
\(ζ(n)\)
we have
\begin{equation}\label{L060}
\neg\exists n(ζ(n)=\lambda x.0) \rightarrow \neg (A \vee \neg A)
\end{equation}
which together with
\(\neg\neg (A \vee \neg A)\)
gives
\begin{equation}\label{L058}
\neg\neg \exists n(ζ(n)=\lambda x.0)
\end{equation}
So the function 
\(Z\)
cannot fail to diverge from 
\(\sum_{ν}ω_{ν}\).

In the remainder we only consider
\(r=0\).

We have
\begin{equation}\label{L059}
\exists n(ζ(n) = \lambda x.0) \rightarrow \neg\exists s\exists m\forall ν(\absval{d_{sν}} > 2^{-m})
\end{equation}
because if 
\(\exists n(ζ(n) = \lambda x.0)\)
then for any 
\(s\),
once
\(ν\)
has become large enough
and
\(i_{ν}\)
small enough,
all
\(d_{sν}\)
vanish.
From~\eqref{L058}
and~\eqref{L059}
we have
\begin{equation}
\neg\neg\neg \exists s\exists m\forall ν(\absval{d_{sν}} > 2^{-m})
\end{equation}
which is equivalent to
\begin{equation}
\neg \exists s\exists m\forall ν(\absval{d_{sν}} > 2^{-m})
\end{equation}
So
\(Z\)
is weakly differentiable
in
\(x=0\)
and there has the weak differential quotient
\(0\).

Now suppose
that
\(Z\)
is also strongly 
differentiable
in
\(x=0\)
and there has the strong differential quotient
\(0\).

From the definition of the functions
\(ω_{ν}\)
it can be derived that there are two lines
from the origin that
are tangent to all of them
and bound all of them from above.
These lines have slopes
\(\surd{2}/4\)
and
\(-\surd{2}/4\)
(Figure~\ref{L063}),
and pass through each
\(ω_{ν}\)
at the points with the respective
\(x\)
coordinates
\begin{equation}
x_{ν}=\frac{2^{(-2ν+1)/2}}{3}  
\end{equation}
and
\(-x_{ν}\).
For any
\(m \in\numnat\),
we can find
a
\(ν\)
such that
all positive values
of
\(ω_{ν}\)
lie
in the interval
\((-2^{-m},2^{-m})\).
For such a
\(ν\),
the difference quotients
relative to the function
\(\sum_{ν}ω_{ν}\)
of the intervals
\([-x_{ν},0]\)
and
\([0,x_{ν}]\)
attain
the absolute value
\(\frac{1}{4}\sqrt{2}\);
and these intervals
can be used in the construction of 
sequences
\(s\).
In contrast,
for 
\(c=0\),~\eqref{L061}
gives
\begin{equation}\label{L062}
\forall n\exists m\forall s\forall ν(\absval{i_{sv2}-i_{sv1}} < 2^{-m} \rightarrow \absval{d_{sv}} < 2^{-n})
\end{equation}
so for each
\(s\)
the difference quotients 
relative to 
\(Z\)
of any of its intervals
will
for all
\(n\)
eventually
become
\( < 2^{-n}\)
in absolute value,
and therefore  
\(\klo \frac{1}{4}\sqrt{2}\).
It follows that on such intervals
none of the
values
\(Z(x)\)
will have been
generated by an
\(ω_{ν}\)
anymore.
It must then be the case that
\begin{equation}\label{L064}
\exists n(ζ(n)=\lambda x.0)
\end{equation}
and hence,
by the definition of the
\(ζ(n)\),
\begin{equation}
A \vee \neg A
\end{equation}
By hypothesis,
\(A\)
is not yet decidable;
therefore,
\(Z\)
is not strongly differentiable.
\end{pargmt}

In this argument,
the argument
from the assumption that
\(Z\)
is strongly differentiable in
\(x=0\)
amounts to an argument that,
once 
\(A\)
is given,
there exists a sequence
\(ζ^\ast\)
defined by
\begin{equation}
ζ^\ast(n) = 
\begin{dcases*}
0  & \text{if  
\(\forall k \leq n(ζ(k) \neq \lambda x.0)\)}\\
1 &  \text{otherwise}
\end{dcases*}
\end{equation}
which has,
in particular, 
the property that
\(\exists n(ζ^\ast(n)=1) \rightarrow A \vee \neg A\);
hence it is,
in effect,
a use of an instance of
\bks{-}{}.%
\footnote{The referee pointed out that there is a considerably simpler
plausibility argument for Weak counterexample~\ref{L129}:
\begin{quote}
Let
\(f\)
be the function from
\(\numreal\)
to 
\(\numreal\)
such that
\(\forall x \in \numreal 
[f(x) \apart 0
\leftrightarrow
\exists n[α(n) \neq 0
\wedge
\exists a \in (0,1)
[x= a/(n+1) + (1-a)/n
\wedge
f(x)=\inf(a,1-a)/n
]
]
]\).
If one assumes:
\(\neg\neg\exists n[α(n) \neq 0]
\wedge
\forall m\forall n
[
(α(m) \neq 0
\wedge
α(n) \neq 0)
\rightarrow
m=n
]\),
and is unable to prove:
\(\exists n[α(n) \neq 0]\),
then 
\(f\)
is weakly differentiable at
\(0\)
(with outcome \(0\))
but one is unable to prove:
\(f\)
is strongly differentiable at
\(0\).
\end{quote}
The referee adds that the example shows that the notion
of weak differentiability is not very useful.
I mention these points for their intrinsic interest.}

A (mis)application of
MP
would
have led one to conclude 
from~\eqref{L058}
to~\eqref{L064};
see on this point also
subsection~\ref{L122}
below.

\subsection{\texorpdfstring{\bks{+}{}}{BKS+}}\label{L072}
That Brouwer had asserted
\bks{+}{}
himself
was seen by
Myhill~\citep[p.295]{Myhill1967},
with reference to the point in
‘Points and spaces’
where 
Brouwer
indeed 
can be said to demonstrate
\bks{+}{}
when
he
proves the equivalence
\(A \leftrightarrow C(γ,A) \in \numrat\)
(see Plausibility argument~\ref{L068} in subsection~\ref{L069} above).
From the equivalence
an
explicit construction for a sequence 
\(α\)
witnessing 
\bks{+}{}
is obtained immediately by
setting
\begin{equation}
α(n) =  
\begin{cases}
1 & \text{if } \exists k \leq n(C(γ,A)(k)=C(γ,A)(k+1))\\
0 & \text{otherwise}
\end{cases}
\end{equation}
A hypothesis in Brouwer’s argument is that 
\(A\) 
has been recognised neither as tested nor as testable.
Strictly speaking,
Brouwer therefore demonstrates
\bks{+}{} only for  
\(A\)
satisfying that condition;
but there is nothing in the construction 
of the checking number
\(C(γ,A)\)
itself
that hinges on this.
For Brouwer’s purpose in the argument at hand,
a weak counterexample to the quantified form of the Principle of the Excluded Middle,
there naturally is no use for greater generality.
There is,
however, 
in his argument from 1949,
discussed in subsection~\ref{L066},
because various real numbers may well have already been proved to be rational
before the construction of the drift
\(γ\)
has 
even begun.
Indeed,
as was noticed in Definition~\ref{L070} in subsection~\ref{L071},
Brouwer’s notion of a checking number of a drift through 
\(A\)
imposes no condition on the status of
\(A\),
and explicitly allows for the case that
\(A\)
has been decided before the construction of the sequence begins.
As a general principle,
then,
\bks{+}{}
will be proved in the same way.
There are,
however,
no later counterexamples in which Brouwer
exploits the equivalence again,
and which could therefore be reconstructed
as uses of 
\bks{+}{}.
(Brouwer stopped publishing altogether the next year.)
Brouwer’s construction dependent on a proposition
\(A\)
of a number that is rational if and only if
\(A\)
is true
will have been Brouwer’s closest approximation to  truth values.%
\footnote{Joachim Lambek recounts~\citep[p.62]{Lambek1994}
how when he met Brouwer in 1953
(at the conference where Brouwer presented ‘Points and spaces’ as a lecture series),
the latter expressed doubts whether Wittgenstein had made any
contributions to logic. 
When Lambek suggested that he had,
because he had come up with truth tables,
Brouwer asked:
‘What are truth tables?’}

\section{Discussion of objections to \cs{}{} and \bks{}{}}

From the beginning,
several objections have been raised 
to both the weak and strong versions of
\cs{}{}
and 
\bks{}{}.
Here I present a survey,
including only arguments 
that
either were made within a broadly Brouwerian framework
or 
would
(prima facie)
be translatable into one.
For example,
it may be observed 
that it is highly unlikely that a meaning explanation
of the propositional operator 
\(\csop_{n}\)
will be found 
in Martin-Löf’s Constructive Type Theory:
in that theory there are no choice sequences as standard objects,%
\footnote{As Martin-Löf 
shows~\citep{Martin-Lof1990}, 
a simple kind of choice sequence can be defined as a nonstandard type.}
mathematical truth is not tensed,
and first-person knowledge plays a role, 
but not in the content of mathematics.
But to advance an observation of this difference as an 
objection would be question-begging.
The real discussion between the opposing views of Brouwer and
Martin-Löf must be about
more general matters than any specific principle in mathematics or logic.
The most interesting and potentially most effective objections
to 
\cs{}{}
and
\bks{}{}
are those that 
aim to show that
accepting them goes
against the philosophical intentions,
accepted evidence,
or standards
of Brouwer’s program itself.%
\footnote{This reflects a view on philosophical controversy 
developed by Henry Johnstone~\citep{Johnstone1959}.} 
(The objection of conflict with Church’s Thesis
in the next subsection
is,
in this respect, 
best seen as a boundary case.)
I will attempt to show that
none of the objections discussed below succeeds.
However,
several objections point to delicate matters of mathematical or philosophical interest.
The grouping of the objections is loosely thematic.

\subsection{\texorpdfstring{\bks{-}{} contradicts Church’s Thesis}{BKS- contradicts Church’s Thesis}}\label{L034}

The following result
is
due to Kripke:
\begin{thm}\label{L084}
\bks{-}{} entails that there exists an effective but non-recursive function.
\end{thm}

Kripke did not publish his argument,
but it was presented by Kreisel in 1970.
It uses the formulation for  species 
\bks{-}{SF} (p.\pageref{L087} above),
which for convenience is repeated here:
\begin{equation}
\exists f\forall x\left[%
\begin{gathered}
\neg\exists n(f(n)=x) \rightarrow x \not\in X\\
\wedge\\
\exists n(f(n)=x) \rightarrow x \in X
\end{gathered}
\right]
\tag{\bks{-}{SF}}
\end{equation}

\begin{prf}[{{\citep[p.145n10]{Kreisel1970a}}}]
CT is formulated as
\begin{equation}\label{L088}
\forall f \exists e \forall n \exists p [T_1(e,n,p) \wedge U(p)=f(n)]
\end{equation}
where
\(T\)
is 
Kleene’s
\(T\)-predicate.

Combining 
\bks{-}{SF} 
and CT
yields
\begin{equation}\label{L089}
\exists e \forall n \left[%
\begin{gathered}
\neg\exists m\exists p (T_1(e,m,p) \wedge U(p)=n) \rightarrow n \not\in X\\
\wedge\\
\exists m\exists p (T_1(e,m,p) \wedge U(p)=n)\rightarrow n \in X
\end{gathered}
\right]
\end{equation}

Consider the species
\(X_0 = \{n \mid \forall y \neg T_1(n,n,y)\}\).
Kleene has shown
that 
it
is not recursively enumerable,
which may be paraphrased by saying that
if a recursive function is total,
it does not have 
\(X_0\) as its range:
\begin{equation}\label{L098}
\forall m\exists p T_1(e,m,p) \rightarrow 
\neg\forall n
\left[%
\begin{gathered}
\neg\exists m\exists p (T_1(e,m,p) \wedge U(p)=n) \rightarrow n \not\in X_0\\
\wedge\\
\exists m\exists p (T_1(e,m,p) \wedge U(p)=n)\rightarrow n \in X_0
\end{gathered}
\right]
\end{equation}
Kreisel points out that Kleene’s theorem has a proof in Heyting Arithmetic. 
Instantiating~\eqref{L089} 
with
\(X=X_0\) 
yields an
\(e\)
that
satisfies the antecedent of~\eqref{L098},
so that instance 
and~\eqref{L098}
are contradictory,
hence
the conjunction of
CT and \bks{-}{SF} 
is,
as Heyting Arithmetic is not contested.%
\footnote{Brouwerians may find the inclusion 
in its logic of Ex Falso Sequitur Quodlibet objectionable~\citep{Atten2009e}, 
but that is not relevant here.}
\end{prf}

In Kreisel’s presentation,
the proof applies
\bks{-}{SF}
to
\(X=X_0\) 
only indirectly,
after combining it with CT;
it thus leaves implicit that
the
\(f\)
that a direct application of
\bks{-}{SF}
to
\(X=X_0\) 
would give
is a counterexample to CT\@.
The existence of a counterexample
is more prominent in 
Myhill’s version.
It uses
%
%
\begin{lem}[{{\citep[p.296–297]{Myhill1967}}}]\label{L081}
Let
\(A(n,x,y)\)
be an arbitrary predicate
with only the free variables
\(n,x,y\),
all ranging over 
\(\numnat\).
Then
\bks{-}{}
implies that
\(\exists f \neg\exists n\forall x\forall y(f(x)=y \leftrightarrow A(n,x,y))\).
\end{lem}

Myhill does not go on to give a proof of 
this lemma,
but one is given by Dragálin;
I here modify it slightly
to highlight the fact that
it
in effect
consists 
in an application
of
\bks{-}{S}
to the species
\(X_A=\{n \mid \forall k A(n, j(n,k), 0)\} \).

\begin{prf}
This is an adaptation of the proof given by Dragálin~\citep[p.134–135]{Dragalin1988}.

Applying
\bks{-}{S}
to the species
\(X_A=\{n \mid \forall k A(n, j(n,k), 0)\} \)
yields a sequence
\(β\)
such that
\begin{equation}\label{L082}
\forall x
(\forall y(β(j(x,y))=0) \leftrightarrow x \not\in X_A)
\end{equation}
Set 
\(f=β\).
Then
\(\neg\exists n\forall x\forall y(f(x)=y \leftrightarrow A(n,x,y))\).
For assume,
towards a contradiction,
that
for some
\(z\),
\begin{equation}
\forall x\forall y(f(x)=y  =y \leftrightarrow A(z, x, y))
\end{equation}
This implies,
given that
\(x=j(v,k)\)
for unique
\(v\)
and
\(k\),
that
\begin{equation}\label{L083}
\forall v(\forall k(f(j(v,k)) = 0) 
\leftrightarrow 
\forall k  A(z,j(v,k),0))
\end{equation}
Instantiating~\eqref{L082}
with
\(y=z\)
and~\eqref{L083}
with
\(v=z\),
it follows that
\begin{equation}
z \in X_A \leftrightarrow z \not\in X_A
\end{equation}
\end{prf}

\begin{prf}[of Theorem~\ref{L084}]
Set
\(A(n,x,y) \coloneqq \exists w (T_1(n, x, w) \wedge U(w)=y)\)
and
apply Lemma~\ref{L081}.
\end{prf}

If 
\cs{+}{} 
is used instead
of one of the variants of
\bks{-}{},
Theorem~\ref{L084} 
is proved as follows:
\begin{prf}[of Theorem~\ref{L084}]\label{L127}
This proof has been adapted from that given by Van Dalen~\citep[p.40n3]{Dalen1978}.%
\footnote{In the exposition I gave of this proof
in~\citep{Atten2008a},
the appeal to Post’s Theorem,
while to my mind correct –~it is derivable from
MP for  primitive recursive predicates~\citep[p.205]{Troelstra.Dalen1988},
which is a highly plausible principle~\citep[section 2]{Atten-forthcomingB}~–
is wholly superfluous.}

Let 
\(K\) 
be a species
that is
recursively enumerable, 
but not recursive.
Define
\begin{equation}
f(n,m) = \left\{
\begin{array}{ll}
0 & \text{if }  \neg\csop_{m} n \not\in K\\
1 & \text{if \makebox[\widthof{\(\neg\)}]{}} \csop_{m} n \not\in K
\end{array}
\right.
\end{equation}
Then
\begin{equation} 
n \not\in K \leftrightarrow \exists m f(n,m)=1
\end{equation}
From left to right,
this follows from 
\cs{+}{3a};
from right to left,
this follows from 
\cs{+}{3b}.
For the Creating Subject
this function 
\(f\) 
is
computable,
as,
by
\cs{+}{1}, 
for any given 
\(m\)
and
\(n\),
\(\csop_{m} n \not\in K\)
is decidable.
Assume that 
\(f\) 
is moreover recursive.
Then the species 
\(S=\{n \mid \exists m f(n,m)=1\}\) 
is recursively enumerable,
but as
\(S=K^c\),
this contradicts the hypothesis.%
\footnote{According to Dirk van Dalen
(personal communication),
this argument was considered common
knowledge at the Summer Conference on Intuitionism and Proof Theory,
SUNY at Buffalo, 1968.
It is,
however,
not found in either of the two publications that came out of that meeting,
Troelstra’s 
\textit{Principles of Intuitionism}~\citep{Troelstra1969}
and the proceedings edited by
Kino,
Myhill,
and Vesley~\citep{Kino.Myhill.Vesley1970}.}
\end{prf}

Although Theorem~\ref{L084}  can be,
and has been,
taken as casting doubt on
\bks{-}{},
the real target of
any such doubts is of course the notion of 
constructive non-recursive sequence itself;
in its simplest form,
it is the doubt that lawless sequences are
individual mathematical entities.
In subsection~\ref{L090},
I argue that Brouwerians have no reason
to harbour that doubt.
Unsurprisingly,
\bks{-}{}
is not the only
principle in the theory of choice sequences
that contradicts CT:
\begin{thm}[{{\citep[p.211]{Troelstra.Dalen1988}}}]\label{L108}
WC-N contradicts CT\@.
\end{thm}

\begin{prf}
Apply WC-N 
to
CT as formulated in~\eqref{L088};
then for every function
\(f\),
computed by a recursive function with index
\(e\),
there is an
\(m\)
such that
any function 
\(g\)
that agrees with
\(f\)
on the arguments
\(1, \dots, m\)
is likewise
computed by the recursive function with index
\(e\).
But then 
agreement of two functions on
\(1, \dots, m\)
implies
agreement 
everywhere,
which is absurd.
\end{prf}


\subsection{\texorpdfstring{\bks{+}{} and MP together imply DNE and PEM}{BKS+ and MP together imply DNE and PEM}}\label{L122}

In 1980,
Joan Moschovakis
found a connection between
\bks{-}{} and MP:
\begin{thm}[{{\citep[p.250-251]{Moschovakis1981}}}]
Let
\(A\)
be any proposition.
Then
\bks{-}{} and MP imply Double Negation Elimination,
\(\neg\neg A \rightarrow A\).
\end{thm}

\begin{prf}
Assume
\(\neg\neg A\).
Applying
\bks{-}{} 
gives
\begin{equation}
\exists α\left[%
\begin{gathered}
\forall n(α(n)=0  \vee α(n)=1) \\
\wedge \\
\forall n(α(n)=0) \leftrightarrow \neg A\\
\wedge\\
\exists n(α(n)=1) \rightarrow A
\end{gathered}
\right]
\end{equation}
We then have
\(\neg\neg A \rightarrow \neg\forall n(α(n)=0)\)
and therefore
\(\neg\forall n(α(n)=0)\)
and
\(\neg\neg\exists n(α(n)=1)\).
By
MP,
now
\(\exists x(α(n)=1)\)
and
hence
\(A\).
So
\(\neg\neg A \rightarrow A\).
\end{prf}

Troelstra and Van Dalen
give the similar direct argument
for the conclusion that
\bks{+}{} and MP imply PEM,
shown by applying 
(in effect)
\bks{-}{}
to
\(A \vee \neg A\)~\citep[p.237]{Troelstra.Dalen1988}.%
\footnote{Compare also Myhill’s refutation of MP 
from 1963~\citep{Myhill1963},
which does not use 
\bks{}{}
(which had not been isolated yet)
but whose setting is that of Kreisel’s Theory of Constructions.}

Depending on one’s views,
these arguments either refute
\bks{-}{}
and
\bks{+}{},
or MP\@.
However,
as Troelstra and Van Dalen 
remark,
their refutation of MP from \bks{+}{}
goes through
‘in an axiomatic setting’~\citep[p.237]{Troelstra.Dalen1988},
and the same can be said of
the above theorem.
What would be needed for contentual
arguments are notions of sequence for which
both MP and \bks{-}{},
respectively \bks{+}{}
hold.
But MP is only plausible
(to my mind, highly so)
for recursive sequences,
whereas the Brouwerian justification of 
either version of 
\bks{}{}
depends on the sequence of acts of the Creating Subject,
which,
as recalled in the previous subsection, 
is certainly not recursive.

\subsection{\cs{}{} leads to a paradox}\label{L077}

In 1969,
Troelstra constructed the following paradox in
\cs{}{}.
Assume that at each stage 
\(m\),
the Creating Subject obtains evidence for
one proposition,
\(A^{(m)}\):
\begin{quote} 
[S]ince it is natural
to assume that we know when a conclusion
has the form
‘%
\(a\) is a lawlike sequence’
we have:\\
\begin{tabular}{ll}
\hspace{1em} & 
\(A^{(m)}\)  is a conclusion of the form ‘%
\(a\) is a lawlike sequence’ or\\
\hspace{1em} & 
\(A^{(m)}\)  is a conclusion of another kind.
\end{tabular}

Then it is possible for us to enumerate 
the
\(A^{(m)}\) 
of the form
‘%
\(a\) is lawlike’;
let
\(A^{(bx)}\)
be the
\(x\)\textsuperscript{th}
conclusion of this form,
stating
‘%
\(a_x\) is a lawlike sequence’.
Then
\begin{equation*}
\bigwedge x\bigvee a(A^{(bx)} \equiv \text{\(a\) is a lawlike sequence})
\end{equation*}
and so we conclude to the existence of a
\(b^\prime\)
such that
\begin{equation*}
b^\prime(x,y) \equiv a_x(y),
\end{equation*}
Intuitively
\(c = \lambda x.b^\prime(x,x) + 1\)
is a lawlike sequence,
but then we ought to be able to indicate a
\(z \in N\)
such that
\begin{equation*}
A^{(bz)} \equiv \text{\(c\) is a lawlike sequence}
\end{equation*}
which implies:
\begin{equation*}
\bigwedge x(b^\prime(x,x) + 1 = b^\prime(z,x) )
\end{equation*}
which is contradictory.~\citep[p.105–106]{Troelstra1969}
\end{quote}
However,
the steps towards a contradiction
depends on Markov’s Principle.
Associate to each sequence
\(a_n\)
a sequence
\(a^*_n\)
such that
\(a^*_n(m)\)
is 
\(0\)
if by stage 
\(m\)
the construction of
\(a_n\)
has not yet begun,
and
\(1\)
if it has.
Then
it can be argued
from the essential freedom of the Creating Subject,
that
\begin{equation}
\forall n\neg\neg\exists k(a^*_n(k)=1)
\end{equation}
because the Creating Subject is free to work
towards constructing 
its
\(n\)-th 
lawlike sequence,
but that 
\begin{equation}
\forall n\exists k(a^*_n(k)=1)
\end{equation}
is not true,
because the same freedom allows the Creating Subject,
for given
\(n\),
not to begin 
constructing its
\(n\)-th 
lawlike sequence
by whatever stage
\(k\)
the lawlike sequence
\(a^*_n\)
indicates.
MP is not correct for the sequences
\(a_n\);
but this means that 
constructively
the enumeration
\(b\)
does not exist.
For the details of this argument,
the reader is referred to~\citep{Atten-forthcomingB},
where it is also argued that
this dependence on
Markov’s Principle
survives in both the solution proposed by 
Troelstra~\citep[p.106–107]{Troelstra1969}
and that proposed by Niekus~\citep[section 3]{Niekus1987},
neither of which can therefore be accepted.
Briefly sketched,
Troelstra’s solution consists in a stratification of propositions and the objects they are about,
so that the enumeration 
\(b\)
will be of a higher level than,
and hence not occur among,
the
\(a_n\).
Niekus’ solution consists in a limitation of reference to stages to those in the future,
and insisting that the conclusion drawn at each stage be new;
the combined effect is that
when the sequence
\(c\)
is defined,
the sequences it depends on
will only be constructed in the future,
so that
\(c\)
is itself not among them.
I will have occasion to discuss Niekus’ solution further in
subsection~\ref{L065},
because it involves questions about how to interpret Brouwer that would have to be discussed also
if Troelstra’s Paradox had not existed.

Troelstra
calls the sequence 
\(c\) 
‘lawlike’.%
\footnote{Similarly,
Van Dalen~\citep[p.265]{Fraenkel.etal1973}
and 
Troelstra~\citep[p.126]{Troelstra1981b}
call a sequence whose existence is guaranteed by 
\bks{}{}
‘lawlike’.}
In view of the non-recursiveness of the
Creating Subject’s activity
(see above),
this may be surprising.
However,
earlier in that chapter
Troelstra explained his terminology:
\begin{quotation}
If we have a definite prescription
involving the actions of the creative subject
(by means of a relation like
\(\vdash_n A\))
for determining the values of a sequence,
we speak of an
empirical
sequence.

Our idea of lawlike sequence does not exclude empirical sequences,
at least not as long as we are willing to consider reference to our own
course of activity by means of
\(\vdash_n\)
as ‘definite’.

It is clear,
however,
that e.g.~primitive recursive functions are lawlike in a
stricter,
more objective sense;
their values are independent of 
\textit{future 
decisions
about the order in which we want to make deductions}.~\citep[p.96–97, emphasis mine]{Troelstra1969}
\end{quotation}
So Troelstra is clear that
Creating Subject sequences
(‘empirical sequences’)
do depend on free choices of the Creating Subject,
and that they are not lawlike in a sense 
that would include their being
independent of the Creating Subject’s free choices.
Thus,
he introduces a further term to reflect this distinction:
\begin{quote}
If a sequence
ξ
is defined by a complete description from sequences
\(χ_1, χ_2, \dots\),
without reference to the creative subject,
we shall call
ξ
\textit{mathematical}
or
\textit{absolutely lawlike}
in
\(χ_1, χ_2, \dots\).~\citep[p.97, emphasis mine]{Troelstra1969}
\end{quote}
Troelstra’s motivation for the wide conception of lawlikeness
would seem to be that also in
Creating Subject sequences,
once the Creating Subject is about to choose 
the 
\(n\)-th value,
it cannot do so freely,
as that value is fully determined by something else.
It is just that
this something else itself 
has come about in an exercise of the essential freedom of the Creating Subject:
the Creating Subject’s activity up to the moment of choosing
the 
\(n\)-th value.%
\footnote{This is also remarked on by Niekus~\citep[p.441]{Niekus1987}.}

To avoid the connotation of predetermination often associated with
lawlikeness,
a connotation that Troelstra as we saw does not wish to evoke in all cases,
in the later presentation of Troelstra’s paradox
in 
\textit{Constructivism in Mathematics},
‘lawlike’
is replaced by
‘fixed by a recipe’~\citep[p.845]{Troelstra.Dalen1988}.
We will briefly return to this in subsection~\ref{L065}. 

It is sometimes suggested that the Creating Subject’s activity
is determined by a law,
but just one that it does not know itself.
This suggestion,
however,
is not consistent with the Brouwerian conception of
mathematical existence:
if some object exists mathematically,
this is only because the Creating Subject has defined and constructed it;
all mathematically relevant truth is brought about
in that activity,
and therefore known to the Creating Subject~\citep[section 4]{Atten-forthcomingB}.

\subsection{Proofs are not different at different times}\label{L031}
At the Brouwer memorial symposium in Amsterdam 
in December 2016,
Kripke stated his
motivation for \bks{-}{}
as follows:%
\footnote{This corresponds to the end of his contribution to the present volume,~\citep{Kripke2018}.}
\begin{quotation}
Intuitionistically, to prove a conditional 
\(A \rightarrow B\), 
one
must have a technique so that from any proof of A
one can get a proof of B.

Now, the idea of the schema, based on Brouwer’s
own arguments about the creative subject, is that one
imagines a sequence in time that is 
\(0\)  
as long as 
\(A\) 
has not been proved but is 
\(1\)  
as soon as it has been
proved.

Then one claims that a proof that the sequence is
always 
\(0\) 
amounts to a proof that 
\(A\) 
can never be
proved, i.e.~that it is absurd.

If the sequence gets the value 
\(1\), then 
\(A\) 
has been
proved.~\citep[slides 35 and 36]{Kripke2016}
\end{quotation}


Note here the appeal to the intuitionistic notion of truth:
since one reasons about one ideal subject,
who can carry out whatever mathematical construction can be carried out at all,
then if it is known that 
\textit{it}
will never carry out a construction for a certain proposition 
\(A\)
this can only be because it is not 
\emph{possible} 
to do so.

Kripke then went on
to state his objection to \bks{+}{}:
\begin{quotation}
What would justify the strong form?

Well, as Myhill says,%
\footnote{[Note MvA:~\citep[p.295]{Myhill1967}.]} 
if 
\(A\) 
is ever proved, 
it is
proved at some definite time, 
and then one gets a
value of 
\(α\) 
that is not equal to 
\(0\) 
simply by looking
at your watch (or calendar, as the case may be)

However:

Is  the  idea  of  the  intuitionistic  conditional
(remember that from a proof of 
\(A\), 
one can get a
proof of 
\(B\)) 
such that intuitionistic proofs, which
indeed must take place at some definite time, have
the time of their occurrence as part of the proof, so
that the same proof would be different if it
happened to occur at a different time?

I think not.

In that case, the strong form does not appear to be
justified, but the weak form is.~\citep[slides 36 and 37]{Kripke2016}
\end{quotation}


When Kripke above speaks of 
‘proofs, which
indeed must take place at some definite time’,
I would argue this must be understood as referring primarily
to proofs as token construction processes,
and,
founded on that,
to token constructions in the other two senses;
for it is the tokens that are bound to a definite time.
And when Kripke asks whether the idea is that
‘the same proof would be different if it
happened to occur at a different time’,
I would argue that
that phrase is ambiguous
but that one of the disambiguations
is indeed is correct:
the sameness 
is to be understood as sameness of type,
but the difference as difference  between
tokens.

The conclusion I draw from
Kripke’s objection to 
\bks{+}{}
 is  that he conceives of 
constructive mathematics,
truth and proof in a different way than Brouwer did. 
But,
symmetrically, 
I don’t think that the fact
that for Brouwer 
the strong form is defensible is, 
by itself, 
an argument against Kripke’s conceptions.

\subsection{CS depends on a function that is not unitype}


The type-token distinction has been appealed to before in a similar situation,
by Timothy Williamson in a discussion not of Kripke’s Schema but of Fitch’s Paradox~\citep{Williamson1988}.
Fitch’s Paradox itself,
which from the premise that all truths are knowable
derives the conclusion that all truths are known,
will not concern us here.%
\footnote{Congenial discussions are those by
Sundholm~\citep[p.20–21]{Sundholm2014}
and
Klev~\citep{Klev2016}.}
But in his discussion,
Williamson comes to consider the intuitionistic meaning of 
the principle that if some proposition
\(A\)
is proved, 
then it is proved
at a definite time
\(t\);
Williamson writes this schematically as
\(A \rightarrow \known{A}\),
and he
notes that
there are two different forms of the consequent to consider: 
either 
\(t\) 
is given explicitly as a parameter, 
or 
there is an existential quantifier over moments in time.
The latter corresponds to 
\cs{+}{3a}.

As in my defense of
\cs{+}{3a},
Williamson holds that the distinction between
types and tokens is pertinent here,
but he aims to use it to opposite effect:
like Kripke, 
he argues for the conclusion
that the principle does not hold,
but for a different reason.

First note the particularity of Williamson’s position that 
the notion of type recognised turns out 
to be so specific that types
cannot be divided into subtypes;
the only way the objects of a type can differ
is numerically.
Applied to proofs,
this means that
for Williamson 
two proof tokens are of the same type
if and only if these proof tokens have 
exactly the same
intermediate and final conclusions~\citep[p.431n16]{Williamson1988}.

His argument against
\(A \rightarrow \known{A}\)
then comes down to this.

Let 
\(a_1\)
and
\(a_2\) 
be token proofs of proposition 
\(A\)
effected at times
\(t_1\)
and
\(t_2\),
respectively.
They are both proofs of
\(A\),
so they are of the same type.
Assume that there is a function
\(f\)
that proves the implication
\(A \rightarrow \known{A}\).
Williamson requires that this function
be
‘unitype’,
that is,
if its arguments are of the same type, 
then so should its values be.

First consider the case where the conclusions of
the two token proofs
\(f(a_1)\)
and
\(f(a_2)\)
are given in the first form,
with an explicit parameter
\(t\).
Then
function
\(f\)
takes the token proof 
\(a_1\) 
to a token proof of the proposition that 
\(A\) 
is proved at time 
\(t_1\)
and the token proof 
\(a_2\) 
to a token proof of the proposition that
\(A\) 
is proved at time 
\(t_2\).
These token proofs have different conclusions,
so
\(f\)
does not fulfill the requirement of being unitype.

In the other case, 
the conclusions of
the two token proofs
\(f(a_1)\)
and
\(f(a_2)\)
have the form of an existentially quantified statement,
and 
are identical.
Williamson points out that in that case
these token proofs 
still differ in another respect,
because 
their intermediate conclusions  
to which existential generalisation is applied
must be different.
Then 
\(f\)
is still not unitype.%
\footnote{Williamson goes on to argue that if
\(A\)
has been decided,
a unitype function may still be obtained by mapping every
proof token of
\(A\)
to
a proof that 
\(A\) was proved at time
\(t\)
for some fixed
\(t\)
(for example,
the 
\(t\) at which
\(A\)
was proved for the first time).
But he observes that this will not work for as yet
undecided
\(A\),
and in any case,
the idea to map
a proof token of 
\(A\)
obtained at 
\(t_1\)
to the moment
\(t_0\)
associated with a different proof token
goes wholly against the Brouwerian descriptivist conception:
there is no intrinsic relation between the act in which
the former proof token was created and 
\(t_0\).}

But the notion of type Williamson employs here  is far too restrictive,
and in his discussion I do not find a motivation for it.
Compare for example
the Curry-Howard-De Bruijn isomorphism,
where a proposition is considered to be the type of all its proofs,
however different the latter may be in their inner workings.
Similarly,
Brouwerian species
–~species is the word Brouwer would have used for type~–
may collect construction objects,
so in particular proofs, 
based on whatever similarity we can see between them
in reflection upon the underlying construction acts,
in second-order mathematics~\citep[section 3]{Atten-E}.
In such a setting,
Williamson’s argument would not succeed.
The reason why Williamson actually wants it to succeed is
that he takes the principle
\(A \rightarrow \known{A}\)
to contradict the
idea
that there may well be truths that are not ever known.
He regards the latter a truism,
and perhaps in some or many domains it is.
But in pure mathematics as Brouwer conceives of it,
as an activity in which the mathematical facts
are 
brought about,
‘truth is only in reality,
i.e.~in the present and past experiences of consciousness’~\citep[p.1243]{Brouwer1949C},
and therefore to assume that
\(A\)
is true
is to assume that 
it is known that
\(A\)
is true.

\subsection{\texorpdfstring{\bks{-}{} contradicts \(\forall α\exists β\)-continuity}{BKS- contradicts AαEβ-continuity}}

A generalisation of
C-N
is 
\(\forall α\exists β\)-continuity~[\citealp[IV-21]{Kreisel1963};~\citealp[p.135]{Kreisel1965};~\citealp[p.73]{Kleene.Vesley1965}].
\begin{pri}[\(\forall α\exists β\)-continuity]
\begin{equation}
\forall α\exists β A(α,β)
\rightarrow
\exists Φ\forall α A(α,Φ(α))
\end{equation}
\end{pri}
where
\(Φ\)
is a continuous functional.

No formulation of it is found in Brouwer,
and,
as
Myhill observed~\citep[p.174]{Kreisel1967b},
Brouwer in 1949~\citep{Brouwer1949A}
in effect constructed
a counterexample  to it
–~because,
in Proof~\ref{L074} above,
the functional that
to each real number
\(f\)
associates 
the direct checking number
\(D(γ,α_f)\)
is
discontinuous.
Myhill also established
\begin{thm}[{{[\citealp[p.293–296]{Myhill1967}; \citealp[p.173–174]{Kreisel1967b}]}}]
\(\forall α\exists β\)-continuity 
is not consistent with
\bks{-}{} 
\end{thm}

\begin{prf}[{{\citep[p.173–174]{Kreisel1967b}}}]
Let 
\(α\)
range over real numbers.
By applying 
\bks{-}{}
to the proposition
\(α \in \numrat\) 
we have
\begin{equation}
\forall α\exists β\left[%
\begin{gathered}
\forall n(β(n)=0  \vee β(n)=1) \\
\wedge \\
\forall n(β(n)=0) \leftrightarrow α \not\in \numrat\\
\wedge\\
\exists n(β(n)=1) \rightarrow α \in \numrat\\
\end{gathered}
\right]
\end{equation}
Now if
\(β=Φ(α)\)
for a continuous functional
\(Φ\),
we would have
\begin{equation}
\forall α\exists m\forall γ(\bar{α}m=\bar{γ}m \rightarrow ((\exists n(Φ(α)(n)=1) \rightarrow γ \in \numrat))
\end{equation}
then the continuity of 
\(Φ\)
contradicts the fact that any initial segment of an
\(α \in \numrat\)
can be continued into an irrational number.
\end{prf}

The system FIM of Kleene and Vesley 
includes
\(\forall α\exists β\)-continuity
(there named
‘Brouwer’s Principle for functions’~\citep[p.72]{Kleene.Vesley1965},
which is
perhaps unfortunate 
because,
as Kleene and Vesley are aware, 
Brouwer never formulated it
and moreover accepted discontinuous functionals).
Dragálin calls the discovery of the resulting inconsistency of \bks{-}{}
and FIM
‘the most dramatic moment in the history of Kripke’s scheme’~\citep[p.135]{Dragalin1988}.
Yet
it is clear that 
FIM,
in the form in which it was published,
was not designed to be
suitable for Brouwer’s Creating Subject arguments:
\begin{quote}
Another point in the formalization in FIM criticized by 
Myhill (1967) 
is that 
the free choice sequences are extensional.%
\footnote{[Note MvA:~\citep[p.286]{Myhill1967}.]}
Certainly, 
extensional free choice sequences are intuitionistically acceptable; 
for these, 
one restricts the freedom of the choices only by 
the choice law adopted in advance. 
Since in fact the intuitionistic theory of the continuum can 
be developed using only extensional choice sequences, 
it seems more interesting to do so. 
The complication of nonextensional free choice sequences 
(where at each choice one picks 
both a function value and a new choice law within the preceding one) 
can be left until a 
need arises, 
as perhaps for the formalization of Brouwer’s ‘historical’ arguments.~\citep[p.138n2]{Kleene1968}
\end{quote}
The inadequacy of FIM if it were meant
as a formalisation of Brouwerian foundations
is exemplified by 
the facts that
in FIM, 
MP is formally undecidable~\citep[p.131]{Kleene.Vesley1965},
and that if FIM proves the existence of an individual choice sequence
with a certain property then it proves the existence of a recursive 
sequence with that property.~\citep[p.101]{Kleene1969}%
\footnote{A phenomenological analysis~\citep{Atten2007}
leads to the conclusion that
also non-lawlike choice sequences are individual mathematical objects}
At the same time,
as Vesley remarks in his very useful retrospective~\citep[p.324]{Vesley1980},
FIM could be seen as supporting Brouwer’s
claim,
in his
‘Second act of intuitionism’,
that choice sequences bring something new:
in FIM it is shown that
the Bar Theorem does not hold classically~\citep[p.87—88]{Kleene.Vesley1965}.%
\footnote{Vesley there~\citep[p.326]{Vesley1980}
also rightly emphasises the fact that Myhill and Kreisel
had taken the initiative to see 
\cs{}{} and \bks{}{}
not only as means to formalise Brouwerian counterexamples,
but as potentially useful in positive intuitionistic developments.
On the latter,
see subsection~\ref{L091}.}
And FIM as published may of course be extended
with principles that bring it closer to Brouwerian analysis;
Vesley’s Schema is an example
(see subsection~\ref{L092} below).

The conflict between
\(\forall α\exists β\)-continuity
and
\bks{-}{}
can be
mitigated in two different ways.
Troelstra has remarked that
there may be interesting subuniverses of choice sequences,
such as that described by his own axiomatic theory GC,
that are not closed under \bks{-}{},
in which case it does not serve to refute
\(\forall α\exists β\)-continuity~[\citealp[sections 6 and 9]{Troelstra1968};~\citealp[p.137]{Troelstra1977};~\citealp[p.235]{Troelstra1983}].
Secondly,
the weaker variant
\(\forall α\exists! β\)-continuity,
which is
\(\forall α\exists β\)-continuity
except that the antecedent
requires
uniqueness of
\(β\),
is often sufficient in analysis,
is not inconsistent with
KS~[\citealp[p.152]{Myhill1970}; \citealp{Krol1978}],
and 
follows from
C-N~\citep[p.89]{Kleene.Vesley1965}.

\subsection{\texorpdfstring{\bks{+}{} must be restricted to determinate properties}{BKS+ must be restricted to determinate properties}}\label{L090}

In reaction 
to the conflict of
\bks{-}{} 
(and hence
\bks{+}{})
with 
\(\forall α\exists β\)-continuity,
Johan de  Iongh proposed
(but,
characteristically, 
did not publish)
a restriction
on the propositions that
\bks{+}{}
may be applied to.
His students
Gielen, 
Veldman, 
and De Swart
describe it as follows:
\begin{quotation}
Application of the axiom should be restricted to propositions which are 
determinate in the sense that they do not depend on infinite objects 
which still are under 
construction, 
and are created more or less freely, 
not merely being developed from 
their previously given definition. 
I may need some more thought, 
but not more 
information, 
in order to know if a determinate proposition is true. 

\elide

The restriction proposed by J.J.~de Iongh, 
seems rather natural: 
as long as information 
about a proposition 
\(P\) 
has not yet been completed, 
I cannot really start 
to think about its truth.~\citep[p.126-127]{Gielen.Swart.Veldman1981}
\end{quotation}
Application of 
\bks{-}{}
to a proposition of the form
\(α \in \numrat\) 
is then obviously no longer permitted
for 
\(α\)
that have not been specified by a law.

In the view of De Iongh and some of those inspired by him,
a construction process that is not governed by 
a finite, 
full definition 
does not yield a construction object,
and remains only a construction project
(the term is De Iongh’s~\citep[p.204]{Swart1992}).
The term
‘construction project’
itself is considered to be primitive
and it is acknowledged that construction projects
may involve making more or less free choices~[\citealp{Gielen1983};~\citealp[p.204]{Swart1992}].
It then follows
that lawless sequences are no
construction objects,
but remain partially defined construction projects
(the Creating Subject knows the choices made so far).
But one wonders 
if an identity criterion for partially defined construction projects 
that is allowed to depend on the moment in time at which the project is begun
would not quickly lead 
to an identity criterion for lawless sequences~\citep{Atten2007}.

De Swart~\citep[p.208]{Swart1992}
explains
‘α is lawless’
by the negative characterisation
‘i.e.~there is no finite law that determines α’,
and then rejects that notion
because there seems to be no finite definition of the general notion
‘finite law’
that it presupposes.%
\footnote{He presents a second argument: 
propositional functions should 
be extensional,
but 
‘α is lawlike’ 
yields a true assertion
for the sequence defined as 
\(λx.0\)
and false for a sequence
that ‘accidentally’ only contains zeros
(the quotation marks here are De Swart’s).
However,
if an infinite sequence contains only zeros,
this can only be because the Subject had imposed a law
on it with that consequence;
to suppose otherwise would be to introduce a realist
conception.
Similarly,
I do not understand De Swart’s claim on the same page
that
‘one might imagine lawless sequences α and β which
are intensionally different, 
but extensionally the same’.
(See on this point also Martino’s comment~\citep[p.394]{Martino1985}
on Troelstra’s description of his
‘abstr’ 
operator as a thought experiment.)}

This stands in contrast to the approach of
Troelstra and Van Dalen,
in which lawless sequences are
indeed construction objects,
which as such are to be considered individuals,
and which can be quantified over.
More generally,
on their approach the explanation of
\(\forall α\),
where
\(α\)
ranges over sequences of a certain type,
requires an analysis of what exactly is given to
the Creating Subject when it gives itself a sequence of that type.
It does not require a construction method that would generate
all sequences of that type.

It seems to me that the approach of Troelstra and Van Dalen captures 
Brouwer’s
descriptions such as the following more accurately:
\begin{quote}\label{L126}
Intuitionistic mathematics is a mental construction,
essentially independent of language.
It comes into being by self-unfolding of the basic intuition
of mathematics,
which consists in the abstraction of two-ity.
This self-unfolding allows us in the first instance to
survey in one act not only a finite sequence of mathematical
systems,
but also an infinitely proceeding sequence,
defined by a law,
of mathematical systems previously defined by induction.
But in the second instance it allows us as well to create a
sequence of mathematical systems which infinitely proceeds
in complete freedom or is subject to restrictions which
may be varied in the course of the progress of the
sequence.~[\citealp[p.339]{Brouwer1947};
trl.~\citealp[p.477]{Brouwer1975}]
\end{quote}
and
\begin{quotation}
The first act of intuitionism completely separates mathematics
from mathematical language. \elide
And the basic operation of mathematical construction is the
mental creation of the two-ity of two mathematical systems
previously acquired,
and the consideration of this two-ity as a new mathematical
system.

It is introspectively realized how this basic operation,
continually displaying unaltered retention by memory,
successively generates each natural number,
the infinitely proceeding sequence of the natural numbers,
arbitrary finite sequences and infinitely proceeding
sequences of mathematical systems previously acquired,
finally a continually extending stock of mathematical
systems corresponding to `separable' systems of classical
mathematics.

The second act of intuitionism recognizes the possibility of
generating new mathematical entities:
First,
in the form of infinitely proceeding sequences whose terms
are chosen more or less freely from mathematical entities
previously acquired;
in such a way that the freedom existing perhaps at the first
choice may be irrevocably subjected,
again and again,
to progressive restrictions at subsequent choices,
while all these restricting interventions,
as well as the choices themselves,
may,
at any stage,
be made to depend on possible future mathematical
experiences of the creating subject \elide.~\citep[p.2]{Brouwer1954A}
\end{quotation}
There is no suggestion in these quotations that the infinite sequences described 
are not on a par
and should be divided into those constructions that are proper objects
and those that remain construction projects.
The sequences described in the second act are generated as mathematical entities
just as much as those of the first act.%
\footnote{In a lecture from 1951,
Brouwer even says that the second act is a special case of the first~\citep[p.93n]{Brouwer1981A}.
An analysis of choice sequences as individual,
mathematical objects
can be carried out phenomenologically~\citep{Atten2007}.}
But then so are,
in particular,
the sequences in which every term is chosen freely;
this is a positive characterisation of lawless sequences.%
\footnote{De Swart’s negative characterisation 
seems to be acceptable as an entailment of 
this positive characterisation
on the ground that,
whether the extension of the concept of finite law
is finitely definable or not,
it is certainly part of the meaning of ‘law’ that
it is incompatible with freedom of each choice.}

To a description of the 
‘second act of intuitionism’ 
in another paper than the one just quoted from,
but of the same period,
Brouwer added the following footnote:
\begin{quote}
In former publications I have sometimes admitted restrictions of freedom
with regard also to future restrictions of freedom.
However this admission is not justified by close introspection
and moreover would endanger the simplicity and rigor
of future developments.~\citep[p.142]{Brouwer1952B}
\end{quote}
This has sometimes been taken to be a rejection of
lawless sequences
(which may be seen as governed by the second-order restriction that there will be no
first-order restrictions).
However,
the formulation of the second act there,
just as the one quoted here in the text,
clearly leaves the freedom to make every choice freely.
One may phrase this as a second-order restriction,
but it is not necessary to bring in that concept.%
\footnote{Also Martino~\citep[p.396]{Martino1985}
argues in favour of accepting lawless sequences.}
Brouwer’s footnote I take to doubt rather  the viability of a general theory of 
higher-order
restrictions.
Whether that doubt is justified is a question that I will not go into here.%
\footnote{As Kreisel observed
(speaking of lawless sequences),
the acknowledgement of
a higher-order restriction
rather simplifies further developments~\citep[p.180]{Kreisel1967b};
a point that Brouwer in writing that footnote may have missed~\citep[p.131–132]{Troelstra1977}.
It can be argued~\citep[p.41–42]{Atten2007}
that the so-called 
‘elimination theorem’
for lawless sequences
of Kreisel~\citep{Kreisel1968a}
and Troelstra~\citep[ch.3]{Troelstra1977},
does not entail their elimination from the intuitionistic ontology,
an entailment that indeed was denied by Kreisel himself~\citep[p.225–226]{Kreisel1968a}.}

Gielen, De Swart, and Veldman
reflect on the reason why De Iongh’s
suggestion was not widely adopted:
\begin{quote}
One may wonder why the restriction 
[i.e., on \bks{-}{}] 
we discussed has not found 
general acceptance among those who work on intuitionism. 
One reason for this is perhaps that J.J.~de longh did not advertise his views, 
strongly. 
Besides, 
Brouwer himself on at least one occasion did not follow this path [9]. 
He proved a theorem which comes down to
\[
\neg\forall α[\neg\neg \exists n[α(n)=0] \rightarrow \exists n[α(n)=0] ]
\]
by using his principle without any restrictions.~\citep[p.128]{Gielen.Swart.Veldman1981}
\end{quote}
(Their reference ‘[9]’ is to
Brouwer’s 1949 paper~\citep{Brouwer1949A},
discussed in subsection~\ref{L066}
above.)
The distinction between a construction proper
and a construction project
was well known to Brouwer.
It is essential to his notion of denumerably unfinished sets:
\begin{quote}
[H]ere we call a set
denumerably unfinished
if it has the following properties:
we can never construct in a well-defined way
more than a denumerable subset of it,
but when we have constructed such a subset,
we can immediately deduce from it,
following some previously defined mathematical process,
new elements which are counted to the original set.
But from a strictly mathematical point of view this set does not exist as a whole,
nor does its power exist;
however we can introduce these words here as an expression
for a known intention.~[\citealp[p.148]{Brouwer1907};
trl.~\citealp[p.82]{Brouwer1975}]
\end{quote}
But in the quotations from 1947 and 1954 above we do not see Brouwer say,
analogously,
that sequences that are not completely defined
do from a strictly mathematical point of view not exist as objects,
but that terms for them are introduced as expressions for
a known intention
(namely,
to begin and continue a construction project of a certain kind).
This explains the fact noted in the latter half of
Gielen, De Swart, and Veldman’s reflection.

Still,
the distinction at the basis of De Iongh’s view
between construction processes that are
governed by a full definition of the object under construction
and those that,
as a matter of principle,
cannot be thus governed,
is a principled one of mathematical relevance,
and it is important to realise that,
if a proposed axiom turns out not to hold in general,
it may still hold for one of these two subclasses.

\subsection{\texorpdfstring{A palatable substitute for \bks{-}{}}{A palatable substitute for BKS-}}\label{L092}

Vesley~\citep{Vesley1970}
proposed a schema that 
is implied by 
\bks{-}{}
but does not imply it~\citep[p.199]{Vesley1970}, 
yet allows alternative arguments for Brouwerian counterexamples
without appealing to the Creating Subject.
It is also,
unlike
\bks{-}{},
consistent with FIM
and with
\(\forall α\exists β\)-continuity.
Vesley considers it,
in the title of his paper,
‘a palatable substitute for Kripke’s Schema’.
Moschovakis showed that 
in
FIM extended with Vesley’s Schema
it is consistent to hold that all sequences
are not not generally recursive~\citep{Moschovakis1971}.

The idea behind
Vesley’s Schema is the assertion that
every continuous function whose domain
is a negatively defined dense subset of the continuum
can be extended to a continuous function on the
full continuum
(which classically is the case even without the condition of being negatively defined).
This is immediately related to Brouwer’s counterexamples
because a set such as
\(\{x \in \numreal \mid x > 0\}\)
is an example of such a negatively defined
subset.
The formal version of VS
does not translate this idea literally,
so as to avoid problems that may arise from the notion
of partially defined function,
but distills:
\begin{equation}
\text{Dense}(\neg A(α), v) \wedge \forall α(\neg A(α) \rightarrow \exists βB(α,β)) 
\rightarrow \forall α\exists β(\neg A(α) \rightarrow B(α,β)) \tag{VS}
\end{equation}
where 
\(v\)
is the universal spread,
and
density of the species
\(\{α \mid \neg A(α)\}\)
in
\(v\)
means that
for any initial segment
one can find a
\(β \in v\)
that extends it
and for which
\(\neg A(β)\)
holds.

Vesley shows that
VS is a fragment of
\bks{-}{} (his Theorems~1 and 2),
and in his Theorem~3
demonstrates various Brouwerian strong counterexamples
from it;%
\footnote{VS has been investigated further by
Scowcroft~\citep{Scowcroft1989}.}
in fact its consequence
\begin{equation}\label{L094}
\text{Dense}(\neg A(α), v) \wedge \forall α(\neg A(α) \rightarrow \exists b B(α,b)) 
\rightarrow \forall α\exists b(\neg A(α) \rightarrow B(α,b)) 
\end{equation}
suffices.
Thus,
Vesley establishes in an informative and formally precise way
that the counterexamples in question do not require
the full strength of 
\bks{-}{}~\citep[p.198]{Vesley1970}.

Similarly,
Van Dalen’s
result~\citep{Dalen1999a}
that
\bks{+}{}
entails that
every negative,
dense subset 
\(X \subset \numreal\)
is unsplittable 
(listed in subsection~\ref{L091}
above)
follows from VS;%
\footnote{I thank Joan Moschovakis for pointing this out.}
but the proof from 
\bks{-}{}
(implicit in Vesley, 
explicit in Van Dalen)
is 
epistemologically
to be preferred for the reason given in the previous paragraph.%
\footnote{Loeb~\citep{Loeb2009} 
has shown that Van Dalen’s theorem is even equivalent
to a weakening of 
(the assertion motivating)
VS,
namely to
‘Every function
\(\numreal \backslash\{0\} \rightarrow \numnat\) 
is sequentially nondiscontinuously extensible to a
function 
\(\numreal \rightarrow \numnat\)’.
In another recent formal investigation of VS, 
Vafeiadou~\citep[ch.2, sections 7 and 8]{Vafeiadou2012}
investigates
VS with 
\(\exists! βB(α,β)\)
instead of
\(\exists βB(α,β)\).}

Epistemologically,
on the other hand,
it not at all obvious that,
once one accepts a Brouwerian setting,
VS can be justified 
on grounds that
are clearer 
than those on which
\bks{-}{}
is justified.
In his 
‘\german{Autorreferat}’
in
\textit{\german{Zentralblatt}}~\citep{Vesley1971},
Vesley points out that his paper
gives no intuitive motivation
for VS; 
none is given in that review either,
nor,
it seems,
elsewhere.
Perhaps such a motivation is found one day;
but in the meantime,
one seems to have no foundational alternative to accepting
VS
because of its derivability from
\bks{-}{}.%
\footnote{Vesley~\citep[p.203]{Vesley1970}
sees some evidence for VS in 
the penultimate paragraph of
Brouwer’s 1949 paper on the non-equivalence
of the positive and negative order relations~\citep{Brouwer1949A};
but this is 
an application of
\bks{-}{}
(see subsection~\ref{L095}).
(I thank Joan Moschovakis for discussion
of Vesley’s remark.)
Vesley's reference to p.123 is to the version of Brouwer’s paper in the
\textit{KNAW Proceedings},
which corresponds to
page 38 of the version in
\textit{Indagationes Mathematicae} 
listed in his references.}

Overall,
it is clear that there are various contexts in which
VS
may function as a
‘palatable substitute’
for
\bks{}{},
but that 
the Brouwerian
foundational perspective
is not one of them.
Myhill
concluded that
‘Vesley's work seems 
intended as an extremely illuminating technical contribution rather than as 
a historically accurate rendering of Brouwer's intentions’~\citep[p.176]{Myhill1968}.

\subsection{\texorpdfstring{Brouwer does not want to accept \bks{-}{}}{Brouwer does not want to accept BKS-}}\label{L073}

In  four papers~\citep{Niekus1987,
Niekus2005,
Niekus2010,
Niekus2017},
Joop Niekus argues that Brouwer avoided 
the reasoning step that in effect appeals to
\cs{-}{3b},
\(\exists n\csop_{n}A \rightarrow A\),
and this is taken to indicate that Brouwer had reason to accept neither
\cs{-}{3b} nor,
by implication,
\bks{-}{}.%

To motivate his claim that Brouwer avoids the use of 
\cs{-}{3b},
Niekus reconstructs the argument
from 1948,
discussed above in subsection~\ref{L007}.
In his reconstruction,
he defines
\(<\)
as follows:
\begin{quote}
All handling of real numbers is done via their generating sequences.
For example,
for the real numbers
\(a\)
and
\(b\),
generated by 
\(a_n\)
and
\(b_n\),
\(a < b\)
holds if
\(\exists n \exists k \forall m ((b_{n+m}-a_{n+m}) > 2^{-k})\)
holds.~\citep[p.41]{Niekus2010}%
\footnote{Also~[\citealp[p.432]{Niekus1987};~\citealp[p.5]{Niekus2017}].}
\end{quote}
The reconstruction proceeds as in Plausibility argument~\ref{L049}
on
p.\pageref{L049}
above.
Niekus notes 
that Brouwer does not reason thus:
\begin{equation}\label{L050}
\begin{aligned}
& α > 0 & \text{(assumption)}\\
& \exists n\csop_{n}A & \text{(def.~α)}\\
& A   &  \text{(\cs{-}{3b})}
\end{aligned}
\end{equation}
His diagnosis is that
\begin{quote}
The use of (4) [=\cs{-}{3b}]
would simplify his argument,
and he would not need to resort to an untested proposition,
but could have used an undecided one.
It seems to us he does not want to use 
(4).~\citep[p.36]{Niekus2010}%
\footnote{Similarly,~\citep[p.434]{Niekus1987}.)}
\end{quote}
The simplification,
while modest,
is there,
but it is not so clear that
an untested proposition is
something to be
‘resorted to’:
the notion of testability is
not more complicated
than that of decidability,
as
testability of 
\(A \)
is 
decidability of
\(\neg A \).
And 
if one’s faith that there will always be undecided
propositions is based on the faith that there
will always be completely open problems,\label{L130}
then 
it is as reasonable to expect
that there will always be
untested propositions,
because 
a problem is not
completely open
if for the corresponding proposition
\(A \)
one already possesses
a proof
of 
\(\neg\neg A \).
Two arguments 
where Brouwer could have
used an undecided proposition,
but invokes
an untested one, 
can be found in 
‘Points and spaces’~\citep[p.533]{Brouwer1954A}
and in the
\textit{Cambridge Lectures}~\citep[p.50]{Brouwer1981A}.

Be that as it may,
in fact Brouwer in his 1948 paper had a good reason for using an untested proposition:
it allows him to draw 
not only the conclusion that
\(\neq\)
is an essentially negative property,
but also that
\(<\)
is,
which, 
unlike the former, 
could not have been arrived at
using an undecided proposition
(see the remarks
to Plausibility arguments~\ref{L005} 
and~\ref{L021} 
in subsection~\ref{L007}).

Niekus does not comment on that,
but it leads us to
the main problem with Niekus’ account:
his definition of
\(<\),
and his interpretation of that sign in Brouwer’s text,
correspond
not to that of Brouwer’s
\(<\),
but 
to that of
Brouwer’s
\(\klo\).
Yet,
Brouwer in this paper
uses both
\(<\)
and
\(\klo\),
and distinguishes them at the end
as the 
‘the simple negative property 
\(ρ > 0\)’
and
‘the constructive property 
\(ρ \gro 0\)’;%
\footnote{The  texts by Brouwer that Niekus discusses in his papers  
are 
the second Vienna lecture~\citep{Brouwer1930A},
the second Geneva lecture~\citep{Brouwer1934},
and 
‘Essentially negative properties’
from 1948~\citep{Brouwer1948A}.
In each of these, Brouwer distinguishes between the two orderings.
The same contrast is referred to in the opening line of the companion paper
on the non-equivalence of the two order relations~\citep{Brouwer1949A}.
That paper is referred to by Niekus~\citep{Niekus2017},
but in a different matter,
which I discuss in subsection~\ref{L065}.}
Brouwer evidently trusts that the reader will look up the definitions
in his earlier papers.
As a consequence,
when Niekus observes that
Brouwer does not make the steps in~\eqref{L050},
the explanation is not,
as he claims it is~[\citealp[p.434]{Niekus1987};~\citealp[p.230]{Niekus2005};~\citealp[p.36]{Niekus2010};~\citealp[p.9, p.13]{Niekus2017}],
that Brouwer did not want to appeal to
\cs{-}{3b}.
It is,
at this precise point,
even irrelevant
whether 
\cs{-}{3b}
is acceptable;
as Brouwer is reasoning in terms of
his
\(>\),
not
\(\gro\),
the assumption
\(α > 0\)
does not allow him
to derive,
towards a subsequent application of
\cs{-}{3b},
its antecedent first.
That derivation would require accepting the very form of MP
that Brouwer’s 1948 paper 
provides a weak counterexample to
(Weak counterexample~\ref{L028} above).%
\footnote{Niekus does not comment on that part of Brouwer’s paper.}
Correspondingly,
the explanation of Brouwer’s use of an untested proposition
lies not
in
an unwillingness to use 
\cs{-}{3b},
but in the wish to establish a property
of the negatively defined
\(<\).
Niekus’ observation on Brouwer’s reasoning 
therefore does not lead to a correct argument that
Brouwer did not want to accept
\cs{-}{3b}
(or
\bks{-}{}).

The misconstrual of Brouwer’s argument leads Niekus to introduce
an alternative propositional operator to that in
the Theory of the Creating Subject.
In the latter
the proposition
\(\csop_{n}A\)
is given the meaning
‘By stage 
\(n\) 
the creating subject has made 
\(A\) evident’,
where
\(n\)
ranges over all stages,
past,
present,
and future.
Niekus offers
\(G_{n}\),
which is used to reason about future stages only:
\begin{quotation}
The
\(G\)
is used in temporal logic to express
‘it is going to be the case that’,
and we shall use it similarly.

We imagine our future to be covered by a discrete sequence
of 
\(ω\)
stages,
starting with the present stage as stage
\(0\),
and we define for a mathematical assertion
\(φ\)
\begin{equation*}
G_nφ
\end{equation*}
as:
at the 
\(n\)-th stage 
from now we shall have a proof of 
\(φ\).
The introduction of this term enables is to refine the notion of proof.

In intuitionism,
stating
\(φ\)
means stating the possession of a proof of
\(φ\).
We now demand of such a proof that it can be carried out
\textit{here and now},
i.e.~all information for the proof is available at the present stage.
If future information is involved we use
\(G_nφ\).~\citep[p.37]{Niekus2010}%
\footnote{Correspondingly,~[\citealp[p.434]{Niekus1987};~\citealp[p.226]{Niekus2005};~\citealp[p.8–9]{Niekus2017}].}
\end{quotation}
(For the use of ‘we’ in the explanation of 
\(G_nφ\),
see subsection~\ref{L065} below.)
In one paper~\citep[p.226]{Niekus2005}
the choice is made to let
the values for
\(n\)
at
\(0\),
so that
\(A \leftrightarrow G_0A\).%
\footnote{To make the parallel to the use of
\cs{}{}
as close as possible,
one might want to say that,
in proving that the sequence
α
converges,~\eqref{L054} 
is used for the special case
\(n=0\).}
For this operator,
the analogues to
\cs{+}{2}
and
\cs{+}{3a}
\begin{equation}
\forall n\forall m({G_n}A \rightarrow G_{n+m}A)
\end{equation}
and
\begin{equation}\label{L055}
A \rightarrow \exists n G_nA
\end{equation}
are valid.
But we also have
\begin{thm}[{{\citep[p.37–38]{Niekus2010}}}]
The analogues to 
\cs{+}{3b},
\begin{equation}\label{L053}
\exists n G_nA \rightarrow A
\end{equation}
and to 
\cs{+}{1},
\begin{equation}\label{L054}
\forall n(G_nA \vee \neg G_nA)
\end{equation}
are not valid.
\end{thm}
\begin{prf}
Let
\(A\)
be an undecided proposition.
Define the sequence
\(α\)
by
\begin{equation} 
α(n) =
\begin{dcases}
0 & \text{if  } \neg G_0(A \vee \neg A)\\
1 & \text{otherwise.}
\end{dcases}
\end{equation}
So α will remain 
\(0\) 
until we have decided 
\(A\),
and then it becomes constant 
with all remaining values being
\(1\).
Let
\(k>0\)
be an arbitrary natural number,
and set
\(B \coloneqq α(k) = 0 \vee α(k)k = 1\).
Then we do not now have
\(B\),
for that would mean that we already know now
whether by the 
\(k\)-th choice from now
we shall have a proof of
\(A \vee \neg A\),
which is impossible.
But obviously
\(G_{k+1}(α(k)=0 \vee α(k) = 1)\),
so
\(G_{k+1}B\).
Hence~\eqref{L053}
is not valid.
Now
set
\(A \coloneqq  α(k) = 1\).
If~\eqref{L054}
were valid,
then
we would
already know now
whether by the 
\(k\)-th choice from now
we shall have a proof of 
\(A \vee \neg A\),
which is impossible.
\end{prf}

\begin{crl}[{{\citep[p.38]{Niekus2010}}}]
\bks{-}{} is not derivable from the schemata for 
\(G\).
\end{crl}

Niekus observes that Brouwer’s reasoning in 
Plausibility argument~\ref{L005}
can be construed in terms of
\(G\),
analogously to
Plausibility argument~\ref{L049}
on
p.\pageref{L049},
and that this
requires
(the contraposition of)~\eqref{L055},
but neither~\eqref{L053}
nor~\eqref{L054}.
But
he goes further and claims~\citep[p.38-39]{Niekus2010}
that not only 
is a reconstruction in terms of
\(G_n A\)
instead of
\(\csop_{n}A\)
possible,
this is what Brouwer had in mind.
In evidence he cites a
handwritten note of Brouwer’s,
of which he
later
says that it provides
‘an even more decisive argument
against 3 [=\cs{-}{3b}]’~\citep[p.13]{Niekus2017}
than Brouwer’s not appealing to 
\cs{-}{3b}
in the 1948 paper
(which,
as we saw above,
actually is not an argument against \cs{-}{3b}).
This note,
in the English translation in the
\textit{Collected Works},
runs as follows:
\begin{quotation}
Further distinctions in connection with the excluded middle.

\(\overline{a}\)
will mean: 
\(a\) 
is non-contradictory.

\(\underline{a}\)
will mean: 
\(a\) 
is contradictory.

\(b \text{ \textit{implies} }a\) 
will mean: 
\textit{from now on} 
I have an algorithm which enables me to derive
\(a\) 
from 
\(b\).

The principle of testability can assert:

either: 
\textit{from now on} 
either 
\(\overline{a}\) 
or 
\(\underline{a}\) 
holds, notation: 
\(|a\).

or: from a certain moment in [the] future on either 
\(\overline{a}\)  
or 
\(\underline{a}\)  
will hold, notation 
\(a|\).

Then 
\(a|\) 
is non-contradictory, but 
\(|a\) 
need not to be non\hyp{}contradictory. 
For
instance, let 
\(p\) 
be a point of the continuum in course of development, 
whose
continuation is free at this moment, 
but may be restricted at any moment in the
future; 
then 
\((p \text{ is rational})|\) 
is non-contradictory, 
but 
\(|(p \text{ is rational})\) 
is
contradictory, 
for the complete freedom which exists at this moment makes it
impossible to be sure that the rationality of 
\(p\) 
is contradictory, 
but also to be sure
that it is contradictory that the rationality of 
\(p\) 
is contradictory.
\elide%
\footnote{In this omitted passage,
Brouwer attempts an alternative 
to Proof~\ref{L074} in terms of \(|a\).}
However,
\elide
\(|a\)
does not seem admissible as a mathematical notion.~\citep[p.603–604]{Brouwer1975}
\end{quotation}
Niekus comments:
\begin{quote}
In struggling with his new notion of tensed objects he comes up with an explicit
distinction, 
which is the same as we make. 
For his 
\(a|\) 
and 
\(|a\)  
are the same as our
\(\neg\neg a \vee \neg a\)
and
\(G_n(\neg\neg a \vee \neg a)\). 
At the end of this note Brouwer expresses doubts about
introducing 
\(|a\) 
as a mathematical notion, 
without further argument. 
But we focus
here on the logical distinction. 
That Brouwer, given the distinction, 
would accept
(7) [=\eqref{L053}]
is of course out of question: stating that 
\(\neg\neg a \vee \neg a\) 
is contradictory and that
\(G_n(\neg\neg a \vee \neg a)\) 
is not refutes (7) 
[=\eqref{L053}]
in a very strong way. 
We conclude there is no base for
KS in Brouwer’s creating subject arguments.~\citep[p.38–39]{Niekus2010}
\end{quote}

Niekus abstracts from Brouwer’s unargued reservation about 
\(|a\) 
being mathematical, 
in order to concentrate on the logical distinction.
This cannot be done,
however,
because
on Brouwer’s conception of logic
(see subsection~\ref{L103} above),
on which it is but an application of mathematics,
as opposed to a prior foundation to it,
if a proposition is not mathematical 
then it has no logical significance either.

What may Brouwer’s reservation have consisted in?
The difference at hand
is that between proofs for which all information needed
is available now,
and proofs for which certain information is not yet present but will be generated
along the finite way.
The latter depend on future activity of the Creating Subject
and hence involves its essential freedom.
But 
for testability
the only mathematically relevant consideration 
is that a construction 
for either
\(\neg p\)
or
\(\neg\neg p\)
be finitely effectible,
and then
there is no
properly mathematical reason
to treat those
finite procesess for which
all
information is available before they begin
differently from the others.
A notion of testability that does just that,
such as
\(|a\), 
would for that reason not be 
‘admissible as a mathematical notion’.
With the notion
\(a|\),
on the other hand,
no problem arises as it is inclusive of both kinds.
So to the extent that an implication has to relate mathematical propositions,
and
\(G\)
is not acceptable as a mathematical notion,
neither~\eqref{L055},
which Niekus accepts,
nor~\eqref{L053},
which he rejects,
are acceptable.

This consideration obviously applies to provability in general.
A division among proofs is not mathematically motivated
if defined in terms of
\(G\).
It is for this reason that Brouwer can write,
when discussing the status of mathematical assertions with respect to truth
in a lecture manuscript from 1951:
\begin{quote}
An immediate consequence
[of the introduction of intuitionism]
was that for a mathematical assertion
α the two cases of truth and falsehood,
formerly exclusively admitted,
were replaced by the following three:
\begin{enumerate}
\item 
\(α\) has been proved to be true;
\item 
\(α\)  has been proved to be absurd;
\item 
\(α\)  has neither been proved to be true nor to be absurd,
nor do we know a finite algorithm leading to the statement
either
that 
\(α\)  is true
or that
\(α\)  is absurd.\(^{\dagger}\)
\end{enumerate}
\end{quote}
adding in the footnote 
\(\dagger\):
\begin{quote}
The case that 
\(α\)   has neither
been proved to be true nor to be absurd, 
but that we know a finite algorithm leading
to the statement either that 
\(α\)   is true, or that 
\(α\)  is absurd, obviously is
reducible to the first and second cases.~\citep[p.92]{Brouwer1981A}
\end{quote}
I do take it that a
‘finite algorithm’
may involve
making a specified finite number of choices;
to take an uncontroversial example,
Newton’s Method
for converging to an
\(x\)
such that
\(f(x)=0\)
begins
with the instruction
‘Choose a starting point
\(x_0\)’.
An example in
in Brouwer's writings is found
in ‘Points and Spaces’,
when
about the depth
(‘order’)
at which a choice sequence
(‘arrow’)
meets a bar
he writes:
\begin{quote}
The definition of a crude bar means that for every arrow 
\(α\) 
of 
\(K\) 
the order
\(n(α)\) 
of the postulated node of intersection with 
\(C(K)\) 
must be computable, however
complicated this calculation may be.
The algorithm in question may indicate the calculation of a maximal order 
\(n_1\) at
which will appear a finite method of calculation of a further maximal order 
\(n_2\) at
which will appear a finite method of calculation of a further maximal order 
\(n_3\) at
which will appear a finite method of calculation of a further maximal order 
\(n_4\) at
which the postulated node of intersection must have been passed.~\citep[p.12–13]{Brouwer1954A}
\end{quote}
The fact that the methods of calculation of the various orders
themselves come to appear at various orders,
that is,
at various points in the choice sequence, 
indicates the
possibility that these methods of calculation depend on the
choices made in between the orders in question.

In
2010
Niekus
acknowledges
that 
\bks{+}{}
is found in Brouwer:
\begin{quote}
There is an instance of KS 
in Brouwer’s work, 
from the last year 
in which he published, see Brouwer 1975, p. 525, 11th
line from below.~\citep[p.38n3]{Niekus2010}
\end{quote}
(This is the instance discussed in subsection~\ref{L072} above.)
But he continues by commenting that
\begin{quote}
Whether there are arguments for this specific instance of KS remains an interesting question. 
\end{quote}
It is not clear to me why Niekus does not say here,
on the basis of his own views,
that
there 
are no such arguments.
As regards his further comment that
\begin{quote}
Although there are one or more instances of KS for specific cases
in the work of Brouwer, 
he always carefully avoided its use as a general principle for
an unspecified formula.~\citep[p.38]{Niekus2010}
\end{quote}
However,
\bks{+}{}
as a fully general principle follows by
exactly the same reasoning as Brouwer
employs in his weak counterexample
(subsection~\ref{L075});
and 
for Brouwer there was no need 
to isolate either form of
\bks{}{}
as he could argue directly
from 
(in effect)
the three principles 
\cs{}{},
of which 
they are
immediate consequences
(subsection~\ref{L080}).

\subsection{Brouwer does not appeal to an ideal subject}\label{L065}

Niekus claims that to hold that
the Creating Subject is in some sense an ideal mathematician
is to hold that
the sequences defined in terms of its activity are
‘completely determined’:
\begin{quote}
The method of the creating subject characterizes Brouwer’s
papers after
1945,
when after a long delay he started to publish again.
The method has
always been supposed to be a radically new step in the work
of Brouwer.
The
expression 
‘creating subject’ 
was then interpreted as 
‘the
idealized mathematician’ 
and the generated sequences by the creating subject as
completely determined.

The notion of the idealized mathematician was formalized by
Kreisel
which resulted in the theory of the idealized mathematician.
This theory does
not reflect Brouwer’s reasoning well and it was struck by a
paradox,
discovered
by Troelstra,
that could not be resolved satisfactorily.

We propose a solution of the paradox in which Kreisel’s main
assumptions are
dropped.
A consequence of our solution is that the generated
sequences are no
longer completely determined,
they are choice sequences.
We will conclude that
the method of the creating subject is special,
not because of the introduction
of an idealized mathematician,
but by the systematic application of particular
choice sequences.~\citep[p.2]{Niekus2017}
\end{quote}
For the paradox,
see subsection~\ref{L077}
above,
where it is also mentioned that both the main solution proposed by Troelstra
and
that by Niekus
fail because of their dependence on Markov’s Principle.
Here the focus will be on the alleged property of complete determination.

That property is embodied,
according to Niekus,
in
\cs{-}{1},
which asserts the decidability of
\(\csop_{n}A\).
On the fact that 
\cs{-}{1}
is accepted in the Kreisel-Troelstra theory,
but not on his alternative,
he comments in an earlier paper that
\begin{quote}
For the reconstruction of Brouwer this has the consequence that, 
since A 2.1 
[=\cs{-}{1}]
is not valid anymore, 
we cannot define 
\({(a_n)}_n\) 
completely
\elide
This is what 
is to be expected if we interpret Brouwer’s method as above, 
because then 
\({(a_n)}_n\)
is a choice sequence, 
all of its values yet undetermined.~\citep[p.436]{Niekus1987}
\end{quote}
and
\begin{quotation}
An argument for taking 
\(β\) 
to be lawlike may be that in the TCS 
[i.e., The Creating Subject] 
the stages seem to
have a definite description, expressed by (1)
[=\cs{-}{1}]. 
But in the intuitionistic interpretation,
for a disjunction to hold we need a proof of one of the disjunctive parts. 
In the case of
\(β\) 
this seems to be not evident to us.%
\footnote{[Note MvA:
\(β\)  
is the Creating Subject sequence in Troelstra’s Paradox;
see subsection~\ref{L077}.]}

Let us return to Brouwer’s original use of creating subject, 
let us interpret it as
ourselves and let the stages cover our future. 
We can define 
\(β\)  
as above. 
Then its
values depend on our future results. 
We have no way to determine these values, 
other
than going in time to these stages, 
which are not specified at all. 
We think decidability
is questionable, 
and we do not want to call this sequence lawlike.~\citep[p.36–37]{Niekus2010}
\end{quotation}
More recently,
Niekus has commented
on Troelstra’s 
changed terminology in
1988 
(see subsection~\ref{L077} above):
\begin{quote}
In formulating the paradox Troelstra is now more cautious than he was
in 1969. 
He formulates the paradox with 
\(L(α)\) holds iff 
\(α\) 
‘fixed by a recipe’
instead of 
‘lawlike’. 
But Troelstra does not abandon his main argument from
Troelstra 1969 for calling a CS sequence, and thus the 
\(c(n)\), 
lawlike. 
That is
the decidability of 
 \(\csop_{n}φ\)
expressed by the TIM 
[The Ideal Mathematician]
Axiom 2: 
\(\csop_{n}φ \vee \neg\csop_{n}φ\).~\citep[p.8]{Niekus2017}
\end{quote}
However,
for 
\(\csop_{n}A\),
decidability is not questionable at all.
It is asserted on the ground 
that whether
\(\csop_{n}A\)
is true depends only on the
Creating Subject’s
own activity,
which,
for any specific
\(n\)
in the future,
it can simply carry on
for the finitely many  stages
required to pass
stage
\(n\),
after which simple inspection 
allows to determine whether 
in the preceeding acts
evidence of
\(A\)
was obtained.
This ground clearly is independent of these preceeding acts
being lawlike or to any extent free.
Niekus is, 
of course, 
right that the decision is made by actually
proceeding and making choices,
but he is wrong in rejecting this as the construction method
that proves either one of the disjuncts of
\cs{-}{1}
(where the number of choices required is finite).
After all, 
making a finite number of choices
can be part of
a genuine construction method~–
the theme of the end of the previous subsection~–
and so
\cs{-}{1}
does not entail that
Creating Subject sequences
are completely determined sequences,
as opposed to choice sequences.

Then Niekus continues and objects that Troelstra’s picture
is mistaken at an even deeper level:
\begin{quote}
Neither
does he question the conception underlying the axioms of the TIM: 
an idealized
mathematician, 
all his mathematical activity covered by a sequence of stages. 
This questioning is the key of the solution of Niekus 1987.~\citep[p.8]{Niekus2017}
\end{quote}
But
Niekus’ own conception
just as much entails that
all mathematical activity of the Creating Subject
is covered by a sequence of stages.
For,
as we saw above,
he writes
‘We imagine our future to be covered by a discrete sequence
of 
\(ω\)
stages,
starting with the present stage as stage
\(0\)’:
in particular,
then,
at the very
beginning of all our mathematical activity
we imagine our future to be covered by a discrete sequence
of 
\(ω\)
stages.
As we proceed,
elements in that sequence that at first corresponded to future stages
come to correspond to stages in our past.
But if we are intuitionistically entitled to imagine that sequence in the first place,
then the systematic change in these correspondences 
will not turn it into an unacceptable object.

Therefore also
in Niekus’s framework 
as he describes it 
the term
‘ideal(ised) mathematician’
is called for
(with a somewhat weaker meaning than he attaches to it,
because it does not include the idea that
the sequences it generates are
completely determined).
He writes that
\begin{quote}
According to Brouwer’s view, 
mathematics is a creation of the human mind and by
using the expression 
\textit{creating subject} 
Brouwer only made explicit his idealistic
position; 
it can be replaced by 
\textit{we} or 
\textit{I}.%
\footnote{[Note MvA: These are cases worth distinguishing~\citep{Rootselaar1970}.]}
Interpreted in this way, 
an idealized
mathematician is not needed at all for the reconstruction, 
a simple principle for
reasoning about the future is enough. 
\elide
We interpreted the expression
\textit{creating subject}
as
\textit{we},
and anybody else can interpret it as himself.~\citep[p.32 and p.39, original italics]{Niekus2010}%
\footnote{Also~[\citealp[p.435]{Niekus1987};~\citealp[p.226]{Niekus2005};~\citealp[p.9]{Niekus2017}].}
\end{quote}
At the same time,
as we just saw,
we are asked to imagine our future activity
as an 
\(ω\)-sequence of stages;
but unidealised human beings 
(Brouwer, Niekus, myself)
do not have such a long future.
(See also the beginning of subsection~\ref{L105}.)

\subsection{\texorpdfstring{\bks{-}{} and \cs{-}{} are incompatible with Brouwer’s notion of infinite proofs}{BKS- and CS- are incompatible with Brouwer’s notion of infinite proofs}}

In Brouwer’s Creating Subject arguments,
it is presupposed that evidence comes in an
ω-sequence.
On the other hand,
Brouwer also accepted infinite proofs,
i.e.,
proofs in which the conclusion is
made evident on the basis of
infinitely many elementary inferences.
The locus classicus for this
is footnote 8 in 
‘\german{Über Definitionsbereiche von Funktionen}’
from 1927,
the paper in which Brouwer gives his second demonstration of the Bar Theorem.
In this demonstration,
a well-ordered thin bar is obtained
by effecting a transformation
on the mental proof that a tree contains a
(decidable)
bar,
which in that footnote is claimed to be infinite:
\begin{quotation}\label{Brouwerfootnote}
Just as,
in general,
well-ordered species are produced by means of the two
generating operations from primitive species
\elide
so,
in particular,
mathematical proofs are
produced by means of the two generating operations from null
elements and elementary inferences 
[\german{Elementarschlüssen}]
that are immediately given in intuition
(albeit subject to the restriction that there always occurs
a last elementary inference).
These \textit{mental} mathematical proofs 
[\german{Beweisführungen}]
that in general contain
infinitely
many terms must not be confused with their linguistic
accompaniments,
which are finite and
necessarily inadequate,
hence do not belong to mathematics.

The preceding remark contains my main argument against the
claims of Hilbert’s metamathematics.~[\citealp[p.64]{Brouwer1927B};
trl.~\citealp[p.460n8]{Heijenoort1967}]
\end{quotation}
For
Brouwer’s constructive definition
of well-ordered species,
see
Definition~\ref{L097}
above.

It was Kreisel who noticed 
the contrast between 
the ω-sequence of the Creating Subject’s activities
and 
the transfinite length of canonical proofs,
and he came to see this as an incoherence in Brouwer’s thought.
Of the two contrasting ideas,
Kreisel found the analysis of proofs of the presence of a bar into an 
infinite canonical form to be
‘rather persuasive’~\citep[p.p.247]{Kreisel1967a},
so he located the problem squarely in the idea that the stages
form an ω-sequence,
and made the correspondingobjection to
\bks{}{}.%
\footnote{Similarly,~\citep[p.13]{Veldman1981}.}
It is amusing to hear,
over a number of years,
the crescendo
in Kreisel’s views:
\begin{enumerate}
\item 
1969:
\begin{quote}
there is no clear reason to
restrict oneself to 
\(ω\) stages
when the canonical proofs on p.59 consist
of a transfinite sequence.~\citep[p.61]{Kreisel.Newman1969}
\end{quote}

\item\label{L123} 
1970:
\begin{quote}
The assumptions used in deriving KS-,
namely thinking of the body of mathematical evidence as
arranged in an
\(ω\) order,
seem arbitrary (though not absurd) if,
as in the theory of ordinals,
one also thinks of individual proofs as consisting of
a transfinite sequence of steps
([3], footnote 8).%
\footnote{[Note MvA: The reference is to Brouwer’s footnote quoted on p.\pageref{Brouwerfootnote}.]}
Therefore the inconsistency of (KS) with Church’s thesis
does not,
I think,
refute the latter conclusively.~\citep[p.128]{Kreisel1970a}
\end{quote}

\item 1970:
\begin{quote}
[While validity in Kripke models implies validity in Heyting’s sense,]
the converse
is
dubious because (some of) the author’s counter models
allowed on pages
98–99 picture an essentially more elementary process of
treating
‘evidential situations’ than allowed in intuitionistic
mathematics.
Specifically, the author considers ω-series (in time)
of stages
of evidence while, at least occasionally, Brouwer considered
fully
analyzed proofs with a transfinite number of steps.~\citep[p.331]{Kreisel1970b}
\end{quote}

\item 1971:%
\footnote{This is from Kreisel’s ‘\german{Autorreferat}’ 
in
\textit{\german{Zentralblatt}}
of~\citep{Kreisel1970a};
Kreisel is here referring to his remark quoted in item~\ref{L123} in this list.}
\begin{quote}
At the end of \S4 (p.128) he considers the schema
KS which is inconsistent with CT\@.
(The schema KS was derived by Kripke from Brouwer’s
assertions about
the thinking
subject or, better, from the postulate of an
ω-ordering of
levels of proofs.)
The author’s objection to KS seems to the reviewer
much stronger than the author can have realized,
casting doubt on the interest of the papers
in the volume which are based on KS.%
\footnote{[Note MvA: 
In that volume~\citep{Kino.Myhill.Vesley1970},
\bks{}{} 
is discussed in the contributions by 
Kreisel~\citep{Kreisel1970a}, 
Myhill~\citep{Myhill1970}, 
Van Rootselaar~\citep{Rootselaar1970}, 
Vesley~\citep{Vesley1970}, 
and 
Scott~\citep{Scott1970}.]} 
\elide
In Section 5 the author apparently expects
an
(hypothetical)
abstract theory of functions and proofs to conflict with CT.\@
Without the kind of implausible
ω-ordering of proofs involved
in
KS,
there is no evidence for such a conflict.~\citep[p.301]{Kreisel1971b}
\end{quote}

\item 1972:
\begin{quote}
the contradiction pointed out at the bottom of p.128 of
[6]%
\footnote{[Note MvA: 
Kreisel’s reference ‘[6]’ is to~\citep{Kreisel1970a}; 
see item~\ref{L123} 
in this list.]},
between two well-known assertions of Brouwer;
one concerning the transfinite structure of
(fully analyzed) proofs,
the other concerning an ω-ordering
of the body of mathematical evidence
as the ‘thinking subject’ or,
equivalently,
the idealized mathematician proceeds in time.~\citep[p.325]{Kreisel1972}
\end{quote}
\end{enumerate}

Kreisel’s claim of a contradiction
can be countered by observing that there are two orderings in play,
and that once they are distinguished,
the perceived contradiction disappears.

For each element in the ordered species,
we distinguish the order it has according to the definition of the species
and the order in which it has been generated in time,
that is,
its genetic order.
We can then say,
in the case of the Bar Theorem,
that the elements of a canonical proof get ordered in two different orderings:
in the transfinite well-ordered species that is the
demonstration,
which order indicates where in the demonstration the
element fits in, 
and in the temporal order of the Creating Subject’s acts of mathematical construction,
which is an ω-order.
While a species may be of a greater order type,
our constructive access to it proceeds in
an ω-sequence of acts.

That the Creating Subject can indeed construct the elements of a
transfinite well-ordering
in an ω-sequence of acts is 
because the inductive definition
of well-ordered species
gives the freedom
to insert the elements in their species-order concurrently.
That is to say,
of each species that has been used in 
the construction of the whole species
–~Brouwer calls these its ‘constructive underspecies’~–,
the elements can be constructed
independently
of the construction of the elements
of the other constructive underspecies.
For example,
consider a well-ordering
of type 
\(ω + ω\),
say
the ordered sum
of the species 
\(X_0\)
of the even numbers in their natural order
and the species
\(X_1\)
of the odd numbers in their natural order.
Then the Creating Subject may first construct
\(0\) 
in
\(X_0\),
then 
\(1\)
in
\(X_1\),
then turn back to 
\(X_0\)
and construct
\(2\),
and so on.
More generally,
from the definition
of well-ordered species
one shows
by ordinary induction
that the elements of the species can be 
enumerated~\citep[p.7, p.30]{Brouwer1918B}.
That a conclusion can be drawn from infinitely many premises 
that at no point have all been constructed is because this infinity is governed by a
finite number of laws to construct it;
it is an insight into the construction processes that these laws 
describe that make the conclusion evident.
This is the same as in the case for induction on the natural numbers~\citep{Atten-forthcomingC}.

\section{Concluding remark}

The preceding considerations
indicate
that
the fact that Brouwer in 1954 was able to demonstrate 
\bks{+}{}
was highly dependent on his very specific views
on mathematical objects,
proofs,
truth,
and freedom.
Even slight changes in these notions or the role they are
assumed to play in mathematics may suffice to
make 
all or some versions of
\bks{}{}
or,
similarly,
\cs{}{}, 
implausible or false.
But in Brouwerian intuitionism,
these principles should be used freely.

\appendix

\renewcommand*{\thesection}{\Alph{section}}

\section{Brouwer’s implicit use of MP in 1918}\label{L131}

As Joan Moschovakis has observed in her review~\citep[p.274]{Moschovakis1979}
of vol.~1 of Brouwer’s
\textit{Collected Works},
in 1918 Brouwer once implicitly used MP
in the form
\begin{equation}\label{L133}
\neg\forall n\neg P(n) \rightarrow \exists n P(n)
\end{equation}
where  
\(P\)
is a decidable predicate~\citep[p.17]{Brouwer1918B};
this fact will now be presented in some detail.

The context is a proof that a certain species is closed.
Brouwer defines the species 
\(C\) 
of choice sequences of positive natural numbers
together with a bijective mapping
(that I will call \(f\))
from
\(C\)
to the dyadically expandable real numbers 
in
\((0,1)\).%
\footnote{The term
‘dyadically expandable’
translates the term Brouwer uses in 1925
for this class,
‘\german{dual entwickelbar}’~\citep[p.251]{Brouwer1925A}.}
That is,
\(f\) 
maps
the choice sequence
\(a_1, a_2, a_3, \dots\)
(\(a_i \in \numnat^+\))
to the real number
\begin{equation}
\frac{1}{2} + 
\frac{1}{4} + 
\dots + 
\frac{1}{2^{a_1-1}} +
\frac{1}{2^{a_1+a_2}} +
\frac{1}{2^{a_1+a_2+a_3}} +
\dots
\end{equation}
The intention behind Brouwer’s
somewhat unclear notation here is that
the sequence begins by summing the
first
\(a_1-1\)
values in the series
\(1/2^n\) 
(\(n \in \numnat^+\)),
which,
in case
\(a_1=1\),
is 
\(0\).%
\footnote{See Heyting’s note 6 in~\citep[p.591n6]{Brouwer1975};
that concerns the later presentation of the same mapping in~\citep[p.251]{Brouwer1925A},
but applies here, too.}

In the 1918 paper,
Brouwer
in fact  constructs real numbers in general
from
\(C\),
\(f\),
and a bijection 
from
\((0,1)\)
to
\((-\infty,\infty)\)~\citep[p.9]{Brouwer1918B};
but he quickly and definitively replaced this with 
a
much wider notion where
a real number may be any sequence of rationals
(alternatively, rational intervals)
as long as it converges.%
\footnote{On this replacement,
see~\citep[p.3–4, 4n1]{Brouwer1919A} 
and~\citep[p.955n1]{Brouwer1921A}.
}
It is the latter notion on which Brouwer’s weak and strong counterexamples
depend. 

Through the mapping
\(f\),
the natural order on the real numbers
in
\((0,1)\)
induces an order 
(which I notate as)
\(\leftabstord\)
on
\(C\).
For
\(x,y \in C\),
\begin{equation}
x \leftabstord y 
\equiv 
\exists v \in (0,1)\exists w \in (0,1)
(x = \inverse{f}(v) \wedge y = \inverse{f}(w) \wedge v < w)
\end{equation}
Thus,
increasing 
\(a_1\)
in a given element of
\(C\)
yields a greater element,
whereas increasing an
\(a_i\) 
with
\(i \geq 2\)
yields a smaller element.

Brouwer then sets out to prove that
\(C\),
thus ordered,
is a perfect species.
A species is perfect if it is dense in itself and closed.
He defines closedness as follows:
\begin{dfn}\mbox{}
\begin{quote}
The ordered species
\(M\)
is called 
\textit{closed},
if 
there can exist no
infinite sequence%
\footnote{[Note MvA: ‘\german{Fundamentalreihe}',
defined as any ordered species whose ordering is similar
to that of the natural numbers in their natural order
\citep[p.14]{Brouwer1918B}.
As Van Dalen~\citep[p.238n12]{Dalen2013} observes,
Brouwer’s actual use of the term wavers a bit.
He often but not always means a lawlike infinite sequence.
At times Brouwer makes an explicit distinction between
‘\german{Fundamentalreihe}'
and 
‘\german{unbegrenzt fortgesetzte Folge}',
such that the latter is the wider notion,
e.g.~\citep[p.202]{Brouwer1921A}.
But in the definition and proof under discussion here,
the wider notion is meant.]} 
of closed intervals
\mathlist{i_1, i_2, \dots}
in it 
of which
\(i_{ν+1}\)
is contained in
\(i_{ν}\)
for every 
\(ν\),
and which have no common element.~\citep[p.17, trl.~MvA]{Brouwer1918B}
\end{quote}
\end{dfn}
and then,
applying this definition to 
\(C\)
and
\(\leftabstord\),
presents this argument
(in which I have put the part relevant here in bold):
\begin{prf}\mbox{}
\begin{quote}
Let us now attempt to determine an infinite sequence
of closed intervals
\mathlist{i_1, i_2, \dots},
of which
\(i_{ν+1}\)
is contained in
\(i_{ν}\)
for every 
\(ν\),
and which have no common element. 
Let 
\mathlist{a_1, \dots, a_n, b_{n+1}, \dots} 
and
\mathlist{a_1, \dots, a_n, c_{n+1}, \dots}
(\(b_{n+1} > c_{n+1}\))
be the end elements
of
\(i_1\),
then the end elements
of an arbitrary
\(i_{ν}\)
have
the same initial segment
\mathlist{a_1, \dots, a_n},
whereas
for the later
\(i_{ν}\)
\(b_{n+1}\)
cannot increase
and
\(c_{n+1}\)
cannot decrease.%
\footnote{[Note MvA: 
An increase of 
\(b_{n+1}\)
would make the left end element smaller 
with respect to
\(\leftabstord\),
and
a decrease of 
\(c_{n+1}\)
would make the right end element greater;
but the sequence
\(i\)
is supposed to be decreasing.]}
{\boldmath\bfseries
Now,
as long as
\(b_{n+1}\)
and
\(c_{n+1}\)
retain the same distinct values,
the corresponding
\(i_v\)
contain
the element
\mathlist{a_1, \dots, a_n, c_{n+1}+1, 1, 1, 1, \dots};
in order to be sure that
this element does not belong to all
\(i_v\),
it must be possible to indicate a certain 
\(i_v\)
for which either
\(b_{n+1}\)
has decreased
or
\(c_{n+1}\)
has increased.}
As this reasoning can be repeated at will,
it must be possible to indicate
a later
\(i_v\)
for which 
\(b_{n+1} = c_{n+1} = a_{n+1}\)
will have come to hold,
and the end elements of which therefore
have the same
first
\(n+1\)
numbers.
Let these end elements be
\mathlist{a_1, \dots, a_{n+m}, b_{n+1}, \dots}
and
\mathlist{a_1, \dots, a_{n+m}, c_{n+1},} \dots. 
Then in the same way in which we derived 
from the sequence
\mathlist{a_1, \dots, a_n}
the sequence
\mathlist{a_1, \dots, a_{n+m}},
we can obtain from 
\mathlist{a_1, \dots, a_{n+m}}
a further sequence
\mathlist{a_1, \dots, a_{n+m+p}};
and,
continuing this way,
we can construct an infinitely proceeding sequence
\mathlist{a_1, a_2,} \dots. 
The element of 
\(C\)
that this sequence represents
belongs to 
\textit{all}
\(i_v\)
however,
by which we have reached a contradiction,
and have recognised that the ordered species in question is 
\textit{closed}.~\citep[p.17, original emphasis, trl.~MvA]{Brouwer1918B}
\end{quote}
\end{prf}
Write
\(b(i_v)\)
and
\(c(i_v)\)
for the values of
\(b_{n+1}\)
and
\(c_{n+1}\)
in the endpoints of
\(i_v\),
and assume that for a given
\(i_k\),
\(b(i_k) \neq c(i_k)\).
Then Brouwer in the bold passage concludes from
\begin{equation}
\neg\forall w\neg(b(i_{k+w}) \neq b(i_k)
\vee
c(i_{k+w}) \neq c(i_k))
\end{equation}
to
\begin{equation}
\exists w (b(i_{k+w}) \neq b(i_k)
\vee
c(i_{k+w}) \neq c(i_k))
\end{equation}
which inference corresponds to that licensed by MP in the form~\eqref{L133}.

That Brouwer came to see the problem with this reasoning
is strongly suggested by his next presentation of this proof~\citep[p.461–463]{Brouwer1926A}.
He there has changed the definition of closedness:

\begin{dfn}\mbox{}
\begin{quote}
In a virtually ordered species
\(M\)
an unbounded sequence of closed intervals
\mathlist{i_1, i_2, \dots},
where each
\(i_{ν+1}\)
is a subspecies of
\(i_{ν}\),
is called a 
\textit{hollow sequence of nested intervals}
[\german{hohle Intervallschachtelung}],
if
for each element
\(p\)
of 
\(M\)
a
\(ν_p\)
can be determined
such that
\(p\)
cannot belong to
\(i_{ν_p}\). 
\elide
If in
\(M\)
there can exist no
hollow sequence of nested intervals,
\(M\)
is called
\textit{closed}.~\citep[p.461, trl.~MvA]{Brouwer1926A}
\end{quote}
\end{dfn}
Thus,
the positive information, 
to produce which from the earlier definition of closedness required
MP,
has now become part of the definition of closedness itself.%
\footnote{The notion of a hollow sequence is related to that of a strong Specker double sequence~\citep[p.743]{ArdeshirRamezanian2010}.}
In the 1927 Berlin lectures,
Brouwer explicitly remarked on the 
greater strength of the new definition~\citep[p.40]{Brouwer1992}.

\section{Brouwer’s proof of the Negative Continuity Theorem}\label{L057}

Brouwer’s argument for the weak counterexample to
\(
\forall x\in\numreal(x \in \numrat \vee x \not\in \numrat)
\)
from 1954
(subsection~\ref{L069})
resembles that for his Negative Continuity Theorem from 1927~\citep{Brouwer1927B}.
Both arguments turn on the idea that,
given the definition of a real number
\(r\),
one can construct a real number
\(s\)
that starts out like
\(r\),
but that may come to diverge from it,
depending on an event of which it cannot 
be predicted if and when it will  occur.
In the case of the weak counterexample,
that is the event 
of proving
an as yet untestable
\(A\);
in the Negative Continuity Theorem,
simply the free choice to diverge.
The interpretation and correctness of Brouwer’s
proof of that theorem
has been the subject
of debate~[\citealp{Posy1976,  
Posy2008, 
Troelstra1982,
Martino1985,
Veldman1982,
Veldman1988}],
mostly concerned with the question whether a certain
negation occurring in the argument is weak or strong
and with the  exact grounds on which that negation is introduced.
I will not reconstruct that debate here in detail,
but want to present Brouwer’s argument
as I read it and
make two remarks:
one on Carl Posy’s criticism of that reading,
and one on the importance Brouwer attached to the argument.

The reading of Brouwer’s proof below agrees
with Heyting’s general suggestion 
to read it
in terms of the Creating Subject~[\citealp[p.131]{Heyting1981};~\citealp[p.479]{Troelstra1982}].
It also agrees with the
(in effect)
detailed elaboration of that
suggestion by Martino~\citep[p.383–384]{Martino1985},
and in particular
I agree with the latter~\citep[p.382]{Martino1985}
that
Brouwer’s argument is a proper proof
of Theorem~\ref{L106}
(by contradiction),
and not only a plausibility argument
(by constructing a weak counterexample to its antithesis)
as Veldman has suggested it is~\citep[p.291]{Veldman2001}.%
\footnote{Further on in his paper,
Martino remarks that justifications
of the continuity principle for lawless sequences
have not taken into account their givenness as individuals~\citep[p.386–390]{Martino1985}.
I later analysed the individuality of choice sequences,
and its relation to WC-N,
in my 1999 dissertation~\citep{Atten1999}, 
published in 2007~\citep{Atten2007};
part of it had found its way into
a paper with Van Dalen in 2002~\citep{Atten.Dalen2002}.
Martino’s paper was overlooked in 
(the work for)
each of these
three publications.
I regret that.}

Brouwer introduced
the notion of negative continuity  in
1924.
It was investigated further  by
Belinfante~\citep{Belinfante1929,Belinfante1930a,Belinfante1930b,Belinfante1931}
and Dijkman~\citep{Dijkman1948}.

\begin{dfn}[{{\citep[p.6]{Brouwer1924N}}}]
A sequence of real numbers
\(r_1, r_2, \dots\)
\textit{converges positively} 
to
a real number
\(r_0\)
if
\begin{equation}
\forall p \exists n\forall m(m > n \rightarrow \abs{r_0-r_m} < 1/p)
\end{equation}

A sequence of real numbers
\(r_1, r_2, \dots\)
\textit{converges negatively}  
to
a real number
\(r_0\)
if 
\begin{equation}
\forall p \neg\exists \underline{n} \forall m (\abs{r_0 - r_{\underline{n}(m)}} > 1/p)
\end{equation}
where 
\(\underline{n}\)
is a strictly increasing sequence of natural numbers.

A function 
\(f\)
is 
\textit{negatively continuous}  
at
\(r_0\)
if for every 
sequence of real numbers
\(r_1, r_2, \dots\)
that converges positively to
\(r_0\),
the sequence
\(f(r_1), f(r_2), \dots\)
converges negatively  to
\(f(r_0)\).

A function 
\(f\)
is 
\textit{negatively continuous}  
if it is
negatively continuous  
at
all points in its domain.
\end{dfn}

\begin{thm}[{{\citep[p.62]{Brouwer1927B}}}]\label{L106}
Let
\(f \colon [0,1] \rightarrow \numreal\)
be  a full  function.
Then 
\(f\)
is negatively continuous.  
\end{thm}

\begin{prf}\label{L109}
Let
\(r_0 \in [0,1]\)
be arbitrary,
and  let
\(r_1, r_2, \dots\)
be a sequence of real numbers
that 
converges positively to 
\(r_0\).
Assume,
towards a contradiction,
that
\begin{equation}\label{L107}
\exists p \exists \underline{n} \forall m (\abs{f(r_0) - f(r_{\underline{n}(m)})} > 1/p)
\end{equation}
where 
\(\underline{n}\)
is a strictly increasing sequence of natural numbers.

The Creating Subject constructs a choice sequence
\(r_{ω}\)
of rational numbers 
\(r_{ω}(i)\)
as follows.
\begin{itemize}
\item 
\(r_{ω}(i) = r_0(i)\)
if at the choice of
\(r_{ω}(k)\)
for some
\(k <  i\)
the Creating Subject 
made the free decision to
align all further choices in
\(r_{ω}\)
with
those of
\(r_0\).

\item 
\(r_{ω}(i) = r_{\underline{n}(m)}(i)\)
if at the choice of
\( r_{ω}(k)\)
for some
\(k < i\)
the Creating Subject 
made the free decision to
align all further choices in
\(r_{ω}\)
with
those of
\(r_{\underline{n}(m)}\),
for some
\(m\).

\item  
\( r_{ω}(i) = r_0(i)\) 
otherwise.
\end{itemize}
The decision 
with which sequence to align
the further choices in
\(r_{ω}\)
is free in respect to both its outcome
and the moment at which it is made,
if at all. 
It is a right that the Creating Subject reserves.%
\footnote{%
\label{L134}Brouwer writes:
‘we \elide reserve the right to determine,
at any time after the first,
second,
\dots,
\(m-1\)th,
and
\(m\)-th
intervals have been chosen,
the choice of all further intervals
(that is,
of the
\(m+1\)th,
\(m+2\)th,
and so on)
in such a way that either a point belonging to
\(ξ_{0}\)
or one belonging to a certain
\(ξ_{p_{ν}}\)
is generated.’~[\citealp[p.62]{Brouwer1927B};
trl.~\citealp[p.459]{Heijenoort1967}].
Brouwer’s use of
‘we’
here is conventional,
and does not imply that the Creating Subject is
not an idealised mathematician;
see subsections~\ref{L105} 
and~\ref{L065},
and compare footnote~\ref{L135}.}
As
assumption~\eqref{L107} 
guarantees that
\(f(r_0)\) 
is co-convergent
with no
\(f(r_{\underline{n}(m)})\),
no initial segment of 
\(r_{ω}\)
constructed prior to such a decision
provides sufficient 
information 
to allow a definition of
\(f(r_{ω})\);
This contradicts the hypothesis of the theorem
that
\(f\)
is a full function,
according to which a construction method
for
\(f(r_{ω})\)
is available from the outset.%
\footnote{Such a method may call for first ensuring that at least
\(k\) 
choices
have been made in 
\(r_{ω}\);
see subsection~\ref{L065} above.}
\end{prf}

The negation in the claim
‘%
\(f(r_{ω})\)
is not defined’
above
is a weak one:
it is not excluded that 
\(r_{ω}\)
will be defined later
(by making the required decision).
But no bound can be given on the stage
by which that would have happened.
The essential 
ingredient in Brouwer’s proof,
then, 
is the
contrast between
the condition that 
\(f\) 
be full
and
the fact that
the Creating Subject cannot be obliged to
exercise its right to fix
\(r_{ω}\)
within any specified finite time.

The reserved right amounts to a restriction on
the choices in 
\(r_{ω}\)
that is provisional
in that it can be lifted when the Creating Subject
chooses to do so~[\citealp{Atten.Dalen2002};~\citealp{Atten2007}].
In the latter
of these two references
I suggested,
as I had overlooked before,
that Brouwer here is exploiting just that notion~\citep[p.108–109]{Atten2007}.
It is a different aspect of the Creating Subject than those
described in the axioms of
\cs{}{}.

That suggestion was criticised by 
Posy~\citep{Posy2008} 
who objected that in that case there is something 
‘introspectively disingenuous’
about the Creating Subject’s behaviour in the interpretation
of the proof in terms of provisional restrictions:
‘We know full well that we want
\(f(r^\ast)\) [= \(f(r_{ω})\)]
to be undefined,
and won’t forget that fact in later
2nd and 1st order choices about
\(r^\ast\)’~\citep[p.30–31]{Posy2008}.
But 
the ground of
the undefinedness of
\(f(r_{ω})\)
is not 
the Subject’s
always choosing not to fix 
\(r_{ω}\);
the crux  is that
\(r_{ω}\) 
is a growing sequence
that
is
a genuine real number
from the outset,
even in absence of a bound
on the stage by which 
\(r_{ω}\)
would be fixed,
and that the 
hypothesis of the theorem 
implies that
also
\(r_{ω}\)
is
in the domain of
\(f\).
I agree,
then,
with Posy’s remark
that
‘We can say that
we don’t currently have a grasp of 
\(f(r^\ast)\); 
we cannot say that we can’t have
such a grasp’,
but add that the point of the argument is that
we cannot say that we
must have such a grasp by a specified stage.

It seems to me that 
the fact that the Creating Subject can at no point be obliged
to align the sequence
\(r_{ω}\) 
with one of the others
is
what Brouwer had in mind when he 
commented that the Negative Continuity Theorem is 
‘an immediate consequence of the intuitionistic point of view’~[\citealp[p.62]{Brouwer1927B};
trl.~\citealp[p.459]{Heijenoort1967}],
for it is an immediate
consequence of the constructive freedom of the Creating Subject.
The evidence for this property of the Creating Subject seems to be 
much easier accessible than that for the canonisability of proofs
appealed to in Brouwer’s argument for the Bar Theorem,
and used towards a demonstration of
the positive Continuity Theorem.
Indeed,
in the 1927 paper~\citep[p.62-63]{Brouwer1927B} 
Brouwer says that he had been mentioning the Negative Continuity Theorem
in lectures and conversations since 1918,
but that it was 
‘much later’
–~six years~–
that he could actually prove the Continuity Theorem that 
is made plausible by it,
and refers to his 1924 papers on the topic~\citep{Brouwer1924D2,Brouwer1924G2}.

Admittedly,
when Brouwer again presents a proof of
the Negative Continuity Theorem,
in the \textit{Cambridge Lectures},
and in a similar way announces it as
‘an immediate consequence of the fundamental thoughts of intuitionism
without using spread keys or well-ordered species’~\citep[p.80-81]{Brouwer1981A},
he has something different in mind,
as instead
of a Creating Subject sequence with a provisional restriction,
he
goes on to use a lawlike sequence and a fleeing property 
(see Definition~\ref{L113}):%
\footnote{The change is also remarked on
by Niekus~[\citealp[p.40-41]{Niekus2010};~\citealp[p.4]{Niekus2017}].}

\begin{prf}[of Theorem~\ref{L106}]\mbox{}
\begin{quote}
For,
let us suppose 
that \(y=f(x)\) 
is a full function of 
\(U\)
[the unit continuum];
\(ξ_{0}\) 
a real number belonging
to
\(U\);
\mathlist{ξ_{1}, ξ_{2}, \dots}
an infinite sequence of real numbers of
\(U\)
converging to
\(ξ_{0}\);
\(t\)
a natural number;
and that
\(\absval{f(ξ_{v})-f(ξ_{0})} > 1/t\)
for every
\(v\).

Let
\(g\) 
be a fleeing property
and
\(k_g\)
its critical number.
We define an infinite sequence of real numbers
\mathlist{q_1, q_2, \dots}
in the following way:
\(q_v=ξ_{v}\)
for
\(v \leq k_g\)
and
\(q_v=ξ_{k_g}\)
for
\(v \geq k_g\).
This sequence converges to a real number
\(q_0\),
to which no real number
\(f(q_0)\)
can be assigned.~\citep[p.80-81]{Brouwer1981A}
\end{quote}
\end{prf}
	
But the fact that Brouwer now presents
this simpler proof
is of course no indication  that he
had come to have second thoughts about 
the acceptability of that from 1927.

\section*{Acknowledgements} 

I am grateful to 
Dirk van Dalen,
Ulrich Kohlenbach,
Saul Kripke,
Per Martin-Löf,
Joan Moschovakis,
Joop Niekus,
Carl Posy,
Göran Sundholm,
Anne Troelstra,
Wim Veldman,
and Albert Visser
for discussion of these topics over the years,
and also to Kripke’s assistant Romina Padro.
Joan Moschovakis 
kindly sent corrections and valuable comments
to the penultimate version of this paper,
as  did a most helpful referee.
Earlier versions of some parts of this paper were presented at
‘Dirk van Dalen 80’,
Utrecht, 2013; 
‘Functions, Proofs, Constructions’,
Tübingen, 2014;
the Fifth Formal Topology Workshop, 
Mittag-Leffler Institute, Djursholm, 2015;
the
PhilMath Intersem, Paris, 2015;
‘Intuitionism, Computation, and Proof: Selected themes from the research of G.~Kreisel’,
Paris,
June 2016;
‘Constructive Semantics: Meaning in between Phenomenology and Constructivism’,
Friedrichshafen,
September-October 2016;
and 
‘Logique en question VII’,
Paris,
June 2017.
I thank the organisers for their invitations,
and the audiences for their questions and comments.


\bibliographystyle{plainnat}
\bibliography{Indagationes-VanAtten}

\end{document}